\documentclass[11pt,reqno]{amsart}

\usepackage[text={160mm,240mm},centering]{geometry}            
\usepackage{comment}
\usepackage{graphicx}
\usepackage{amssymb}
\usepackage{epstopdf}
\usepackage{dsfont}

\usepackage{amsmath,amsthm}
\usepackage{longtable}
\usepackage{array}
\usepackage{multirow}
\usepackage{multicol}
\usepackage{booktabs}

\usepackage[mathscr]{eucal}
\usepackage[all]{xy}
\usepackage{mathrsfs}
\usepackage{hyperref}
\usepackage{color}

\newtheorem{theorem}{Theorem}[section]
\newtheorem{lemma}[theorem]{Lemma}

\newtheorem{conjecture}[theorem]{Conjecture}
\newtheorem{corollary}[theorem]{Corollary}
\newtheorem{proposition}[theorem]{Proposition}

\newtheorem{definition-lemma}[theorem]{Definition}
\newtheorem{definition-theorem}[theorem]{Definition}

\newtheorem*{Main}{Main Theorem}

\theoremstyle{definition}

\newtheorem{definition}[theorem]{Definition}

\newcommand{\p}{\partial}

\newcommand{\Diff}{\operatorname{Diff}}
\newtheorem{remark}[theorem]{Remark}

\allowdisplaybreaks

\title[Nakai Conjecture]{The Nakai Conjecture for isolated hypersurface singularities of modality $\le 2$}

\author{Rui Li}
\address{
Zhili College,
	Tsinghua University,
	Beijing, 100084, P. R. China.}
\email{lirui22@mails.tsinghua.edu.cn}
\author{Zida Xiao}
\address{
Department of Mathematical Sciences,
	Tsinghua University,
	Beijing, 100084, P. R. China.}
\email{xiaozd21@mails.tsinghua.edu.cn}
\author{Huaiqing Zuo}
\address{Department of Mathematical Sciences,
	Tsinghua University,
	Beijing, 100084, P. R. China.}
\email{hqzuo@mail.tsinghua.edu.cn}

\thanks{Zuo  was  supported by NSFC Grant 12271280}

\begin{document}

\maketitle

\renewcommand{\abstractname}{Abstract}
\begin{abstract}
The well-known Nakai Conjecture concerns a very natural question: For an algebra of finite type over a characteristic zero field, if the ring of its differential operators is generated by the first order derivations, is the algebra regular? And it is natural to extend the Nakai Conjecture to local domains, in this paper, we verify it for isolated hypersurface singularities of modality $\le 2$, this extends the existing  works.

Keywords. high order derivations, isolated hypersurface singularity, Nakai Conjecture.
	
MSC(2020). 14B05, 32S05.
\end{abstract}

\section{Introduction}
In this paper, $k$ denotes a field of characteristic zero. Let $A$ be a $k$-algebra, $\Diff_k^q(A)$ be the set of $q$-th order differential operators on $A$ over $k$, which is defined inductively by $\Diff_k^q(A)=0$ for $q<0$, and $\Diff_k^q(A)=\{D\in Hom_k(A,A) \ | \ [D,a]\in \Diff_k^{q-1}(A), \forall a\in A\}$ for $q>0$. $Der^q_k(A)$ be the set of $q$-th order derivations on $A$ over $k$, which consists of $D\in Hom_k(A,A)$, satisfying 
$$D\left(a_0 a_1, \cdots, a_q\right)=\sum_{s=1}^{q+1}(-1)^{s-1} \sum_{i_1<\cdots<i_s} a_{i_1} \cdots a_{i_s} D\left(a_0 \cdots \check{a}_{i_1} \cdots \check{a}_{i_s} \cdots a_q\right), \forall a_0,a_1\cdots,a_q\in A. $$

It is known that $\Diff_k^q(A)=Der^q_k(A)\oplus A$, and there exist two natural filtrations $A=\Diff_k^0(A)\subset \Diff_k^1(A)\subset \Diff_k^2(A) \subset \cdots$, $0=Der_k^0(A)\subset Der_k^1(A)\subset Der_k^2(A)\subset \cdots$, denote $Der_k (A)=\mathop{\cup}\limits_{q\in \mathbb{N}}Der^q_k(A)$, $\Diff_k (A)=\mathop{\cup}\limits_{q\in \mathbb{N}}\Diff^q_k(A)$. The compositions of differential operators and derivations satisfy $\Diff_k^p(A)\Diff_k^q(A)\subset \Diff_k^{p+q}(A), Der_k^p(A) Der_k^q(A)\subset Der_k^{p+q}(A)$, the Lie brackets(commutators) of differential operators and derivations satisfy $[\Diff_k^p(A),\Diff_k^q(A)]\subset \Diff_k^{p+q-1}(A), [Der_k^p(A), Der_k^q(A)]\subset Der_k^{p+q-1}(A)$ (\cite{Nakai}). For simplicity, we omit the subscript $k$ in the notations of differential operators and derivations from now on.

 Denote by $der^q(A)$ the $A$-submodule of $Der^q(A)$ which consists of $A$-linear combinations of derivations of the form $\delta_1\delta_2\cdots\delta_j, 1\le j\le q, \delta_i\in Der^1(A), \forall i$, and denote by $der(A)$ the $A$-submodule of $Der(A)$ generated by the compositions of elements in $Der^1(A)$. 

\begin{remark}
    Note that $der(A)=\mathop{\cup}\limits_{q\in \mathbb{N}} der^q(A)$, but $der^q(A)=der(A)\mathop{\cap}Der^q(A)$ does not necessarily hold. Here we give an example to explain this, for more details about the example, one can refer to \cite{XYZ}.
    Let $A=\mathbb{C}\{x,y\}/(4x^3+y^2,xy)$, the Tjurina algebra of the $D_5$ singularity. We have $Der^1(A)=Span_{\mathbb{C}} \{y\p_y+\frac{2}{3}x\p_x, y^2\p_y, x^2\p_y+\frac{1}{4}y\p_x, y^2\p_x, x^2\p_x\}$, and $der^2(A)=Der^1(A)+Der^1(A)Der^1(A)=Der^1(A)\oplus Span_{\mathbb{C}} \{3y\p_y+9y^2\p_y^2+4x^2\p_x^2, y^2\p_x^2, y\p_x-y^2\p_x\p_y, y^2\p_x^2\}$. However, if we let $E=y\p_y+\frac{2}{3}x\p_x$, we see that $9E^4-36E^3+41E^2-14E=-4(y^2\p_y+\frac{1}{3}x^2\p_x)\notin der^2(A)$, and it is obvious that $9E^4-36E^3+41E^2-14E\in der(A)\cap Der^2(A)$. So $der^2(A)\ne der(A)\cap Der^2(A)$.
\end{remark}

Grothendieck \cite{Gr} showed that $Der(A)$ is generated by $Der^1(A)$ when $A$ is of finite type and regular. Nakai conjectured that the converse is also true, it seems that he had not quoted the conjecture rigorously in his papers, but the conjecture is often quoted in connection with his paper \cite{Nakai}. Nakai Conjecture was firstly stated in 1973 by K.R. Mount and O.E. Villamayor \cite{Mount} and they proved it in the case when $A$ is the affine ring of an irreducible algebraic curve.  In \cite{BLLS},  the statement of Nakai Conjecture is as follows.

\begin{conjecture}[\textbf{Nakai} \cite{BLLS} \cite{Nakai}]
	Let $k$ be a field of characteristic zero and $A$ be a $k$-algebra of finite type, if $der^q(A)=Der^q(A)$ for each integer $q\ge 1$, then $A$ is regular.
\end{conjecture}

An interesting result proved by Becker \cite{Be} and Rego \cite{Re} says that the Nakai Conjecture implies another well-known long-standing conjecture of Zariski-Lipman, which is still open in general cases. For recent progress of the Zariski-Lipman Conjecture, one can refer to  \cite{Berg}, \cite{Graf}, \cite{Kal}, \cite{Ste}, and so on.

\begin{conjecture}[\textbf{Zariski-Lipman} \cite{Lipman}]
    Let $k$ be a field of characteristic zero, $A$ be the affine ring or the local analytic ring of a variety over $k$, if $Der^1(A)$ is $A$-projective, then $A$ is regular.
\end{conjecture}

The Nakai Conjecture has been proved for several cases and is still widely open. It is known to be true for the case of algebraic curves \cite{Mount}, the case of monomial ideals \cite{Sch}, the case when $A$ is the invariant subring of $k[x_1,\cdots,x_n]$ acted by a finite subgroup of $GL(n, k)$ \cite{Yas}, the case of a cone on a Riemann surface of genus $>1$ \cite{Paul} and other special cases (\cite{Brown}, \cite{BLLS}, \cite{Traves}).

For the case that $A$ is the affine ring of a hypersurface, Singh \cite{Singh} presented the following conjecture in 1986, which is stronger than the Nakai Conjecture in this case:
if $A=k[x_1,\cdots,x_n]/(f)$ and $der^2(A)=Der^2(A)$, then $A$ is regular (i.e., Singh Conjecture), and he proved it for two variables case. In \cite{Brumatti}, the authors proved it for  isolated hypersurface homogeneous Brieskorn singularities, and Xiao-Yau-Zuo \cite{XYZ} generalized it to weighted homogeneous Brieskorn singularities. Recently, Yau-Zhu-Zuo \cite{YZZ} proved the Nakai Conjecture for the isolated homogeneous hypersurface singularities, they introduced new ideas to analyze the necessary conditions for $D\in Der^2(A)$ to be in $der^2(A)$ and completed the proof by concrete constructions. Furthermore, Xiao-Yau-Zhu-Zuo \cite{XYZZ} generalized it to   weighted homogeneous fewnomial singularities and  cusp singularities.

A natural question is to extend the class of algebra $A$ in the Nakai Conjecture : is it possible to replace the affine ring $A$ by a local domain? In \cite{BLLS}, the authors gave the negative answer by providing an example of one-dimensional pseudo-geometric ring. So the Nakai Conjecture does not hold for an arbitrary local domain $A$. However, in \cite{Be}, Becker showed that the Nakai Conjecture implies the Zariski-Lipman Conjecture both in algebraic and analytic version. In this paper, we consider the Nakai Conjecture for case that $A$ is an analytic local ring of  isolated hypersurface singularities.  

In section \ref{sec 2}, we recall some basic conceptions about high order derivations, and make some preparations for our proof. In section \ref{sec 3} and \ref{sec 4}, we do more concrete calculations case by case for unimodal and bimodal isolated hypersurface singularities classified by Arnold \cite{Arnold}. Our main theorem is stated as following: 

\begin{Main}
    Let $(V(f),0)\subset (\mathbb{C}^n,0)$ be an isolated hypersurface singularity of modality $\le 2$, and $A=\mathcal{O}_n/(f)$ be the algebra of holomorphic function germs of $V(f)$ at the origin, then $der^2(A)\ne Der^2(A)$, i.e., the Nakai Conjecture holds.
\end{Main}



\section{Preliminaries}\label{sec 2}

In this section, we make some preparations for high order derivations. First we fix some notations, let $P=k[x_1,\cdots,x_n]$ or $k\{x_1,\cdots,x_n\}$(the subring  of all convergent power series in $k[[x_1,\cdots,x_n]]$). Let $\mathbb{N}$ be the set of all non-negative integers and put $V=\mathbb{N}^n$. For
$\alpha=\left(\alpha_1, \cdots, \alpha_n\right) \in V$, we use the standard notation: $|\alpha|=\alpha_1+\cdots+\alpha_n$,\  $\alpha !=\alpha_{1} ! \cdots \alpha_{n} !, \ x^\alpha=x_1^{\alpha_1} \cdots x_n^{\alpha_n}$, etc. For $r \in \mathbb{Z}$, let $V_r=\{\alpha \in V|\ | \alpha | \leq r\}$ and $W_r=\{\alpha \in V|\ | \alpha | =r\}$. For $1 \leq i \leq n$, let $e_i=(0, \cdots, 1, \cdots, 0) \in W_1$ with $1$ at the $i$-th position.

When $A=P/I$ with $I$ a proper ideal of $P$, the higher derivations are presented as following:
\begin{theorem}\label{thm2.1} 
	Let  $P=k\left[x_1, x_2, \cdots, x_n\right]$ or $k\{x_1,\cdots,x_n\}$, $I$ be a proper ideal of $P$ and $A=P/I$. Then
	\begin{align*}
	Der^{m}(A)\cong \frac{\{D\in Der^{m}(P),D(I)\in I\}}{I Der^{m}(P)}\cong \{D\in Der^{m}(P,A) | D(I)=0\};
	\end{align*}
   \begin{align*}
	\Diff^m(A)\cong \frac{\{D\in \Diff^m(P),D(I)\in I\}}{I\Diff^m(P)}\cong \{D\in \Diff^m(P,A) | D(I)=0\}.
	\end{align*} 
    
\end{theorem}
From now on, we consider the cases of the $k$-algebra $A$ is of the form $A=P/I$, we will identify derivations(differential operators) in $Der(A)$($\Diff(A)$) with their lifts in $Der(P)$($\Diff(P)$) or in $Der(P,A)$($\Diff(P,A)$) throughout the later discussion. For $\alpha \in V$, denote $\p_x^{(\alpha)}$ the derivation $(1 / \alpha !) \partial^\alpha / \partial x^\alpha: P \rightarrow P$.  As the set of all high  order differential operators of $P$ is well-known as the Weyl algebra $\Diff(P)=P[\p_{x_1},\cdots,\p_{x_n}]$, then every $D \in \operatorname{Diff}(A)$ has a unique expression of the form in $\Diff(P,A)$:
$$
D=\sum_{\alpha \in V} c_{\alpha}(D) \p_x^{(\alpha)}
$$
with $c_\alpha(D) \in A$ for all $\alpha$ and $D(I)\in I$, and  $c_\alpha(D)=0$ for almost all $\alpha$.

\begin{definition}
For $D\in \Diff(P,A)$ and $\beta\in V$, define 
$$
\langle D,x^\beta\rangle=\sum_{\alpha\in V}c_{\alpha+\beta}(D)\p_x^{(\alpha)}.
$$
\end{definition}
Note that if $D\in \Diff^r(P,A)$ then $\langle D,x^\beta\rangle \in \Diff^{r-|\beta|}(P,A)$.

\begin{remark}\label{rmk2.3}
    In theorem \ref{thm2.1}, if $f_1,f_2,\cdots,f_r$ is a system of generators of the ideal $I$, then the condition $D(I)=0$ is equivalent to $\langle D,x^\beta\rangle (f_i)=0$, $\forall \beta \in V_{m-1}, 1\le i\le r$.
\end{remark}

\begin{definition}
Let $\Phi: \Diff(P, A) \times V \rightarrow \operatorname{Der}(P, A)$ be the pairing defined by $\Phi(D, \beta)=\left\langle D, x^\beta\right\rangle-\left(\left\langle D, x^\beta\right\rangle(1)\right)_P$. Note that $\Phi$ is the direct limit of the pairings
$$
\Phi_{m,r}: \operatorname{Diff}^{m}(P, A) \times W_{r} \rightarrow Der^{m-r}(P, A)
$$
given by
$$
\Phi_{m,r}(D, \beta)=\langle D, x^\beta\rangle-c_\beta(D) \p_x^{(0)} .
$$
\end{definition}
\begin{proposition}\cite{Singh}\label{prop2.5}
For  $r\le m$, we have an exact sequence
$$
0 \rightarrow \operatorname{Diff}^{r}(P, A) \longrightarrow \operatorname{Diff}^m(P, A) \stackrel{\theta_{m,r}}{\longrightarrow} \bigoplus_{\beta \in W_r} \operatorname{Der}^{m-r}(P, A),
$$
where $\theta_{m,r}(D)=\left(\Phi_{m,r}(D, \beta)\right)_{\beta \in W_r}$, and $\Diff^r(P, A) \hookrightarrow \Diff^m(P, A)$ is the natural inclusion.
\end{proposition}  

As for $D \in \operatorname{Diff}^m(P, A)$, the following three conditions are equivalent by theorem \ref{thm2.1} and remark \ref{rmk2.3}:\\
(i) $D \in \operatorname{Diff}^m(A)$;\\
(ii) $\left\langle D, x^\beta\right\rangle \in \operatorname{Diff}^{m-|\beta|}(A)$ for every $\beta \in V$;\\
(iii) $\left\langle D, x^\beta\right\rangle \in \operatorname{Diff}^{m-|\beta|}(A)$ for every $\beta \in V_{m-1}$.\\
The exact sequence in proposition \ref{prop2.5} induces the exact sequence 
$$
0 \rightarrow \operatorname{Der}^r(A) \to  \operatorname{Der}^m(A) \stackrel{\theta_{m,r}}{\longrightarrow} \bigoplus_{\beta \in W_r}\operatorname{Der}^{m-r}(A), 
$$
for the case of $r=m-1$, we write $\theta_{m,m-1}=\theta_m$ for simplicity.

\begin{definition}\label{def2.6}
     For $m \in \mathbb{Z}$, define
$$
\begin{aligned}
	\mathscr{D}^{m}(A)= & \{ \left(d_\beta\right)_{{\beta} \in W_{m-1}} \in \underset{\beta \in W_{m-1}}{\bigoplus} \operatorname{Der}^1(A) \mid d_\beta\left(x_i\right)=d_\gamma\left(x_j\right) \quad\text { whenever } \\
	& \beta+e_i=\gamma+e_j, \beta, \gamma \in W_{m-1}, 1 \leq i, j \leq n\}.
\end{aligned}
$$
It is easy to check that $\mathscr{D}^m(A)\subset \operatorname{Im}(\theta_m)$, and in particular, $$
\mathscr{D}^2(A)=\{(d_1,\cdots,d_n)\in\oplus_{i=1}^nDer^1(A)|d_i(x_j)=d_j(x_i) \text{ for all $i,j$}\}.
$$
\end{definition}

However, the expression of the generators of high order derivations on $A$ is rather difficult. The following proposition provides an idea to descent the first order derivations on $A$ to a quotient ring of $A$, which is an Artinian algebra.

\begin{theorem}\label{thm2.7}
    Let  $P=k\left[x_1, x_2, \cdots, x_n\right]$ or $k\{x_1,\cdots,x_n\}$, $I=(f_1,f_2,\cdots,f_r)$ which defines a complete intersection of dimension $n-r$, $A=P/I$. Then the natural projection $$A\to A/(\{\frac{\p(f_1,f_2,\cdots,f_r)}{\p(x_{i_1},x_{i_2},\cdots,x_{i_r})}\}_{1\le i_1<\cdots<i_r\le n})=P/(I,\{\frac{\p(f_1,f_2,\cdots,f_r)}{\p(x_{i_1},x_{i_2},\cdots,x_{i_r})}\}_{1\le i_1<\cdots<i_r\le n})$$ induces a natural map
    $$Der^1(A)\rightarrow Der^1(P/(I,\{\frac{\p(f_1,f_2,\cdots,f_r)}{\p(x_{i_1},x_{i_2},\cdots,x_{i_r})}\}_{1\le i_1<\cdots<i_r\le n})).$$
    Here the notation $\frac{\p(f_1,f_2,\cdots,f_r)}{\p(x_{i_1},x_{i_2},\cdots,x_{i_r})}$ denotes the Jacobian determinant.
\end{theorem}

\begin{proof}
    Denote the ideal $J:=(\{\frac{\p(f_1,f_2,\cdots,f_r)}{\p(x_{i_1},x_{i_2},\cdots,x_{i_r})}\}_{1\le i_1<\cdots<i_r\le n})$, from theorem \ref{thm2.1}, we need to show that for any $D=\sum_{s=1}^n \alpha_s \frac{\p}{\p x_s} \in Der^1(P)$ satisfying $D(I)\subset I$, then $\forall 1\le i_1<\cdots<i_r\le n$, we have $D(\frac{\p(f_1,f_2,\cdots,f_r)}{\p(x_{i_1},x_{i_2},\cdots,x_{i_r})})\in (I,J)$.

    As $D(f_i)\subset I$, write $D(f_i)=\sum_{k=1}^r c_i^{(k)}f_k, c_i^{(k)}\in P$, take the $j$-th partial derivative on both sides, then $\frac{\p}{\p x_j}(D(f_i))\equiv \sum_{k=1}^r c_i^{(k)}\frac{\p f_k}{\p x_j} \mod I$, and $$D(\frac{\p f_i}{\p x_j})=[D,\frac{\p}{\p x_j}](f_i)+\frac{\p}{\p x_j}(D(f_i))\equiv -\sum_{s=1}^n \frac{\p \alpha_s}{\p x_j}\frac{\p f_i}{\p x_s}+\sum_{k=1}^r c_i^{(k)}\frac{\p f_k}{\p x_j} \mod I. $$
    From the Leibniz rule, we have
    \begin{align*}
        &D(\frac{\p(f_1,f_2,\cdots,f_r)}{\p(x_{i_1},x_{i_2},\cdots,x_{i_r})})=
        \sum_{k=1}^r det\left(\begin{matrix}
            \frac{\p f_1}{\p x_{i_1}} & \frac{\p f_1}{\p x_{i_2}} & \cdots & \cdots & \frac{\p f_1}{\p x_{i_r}} \\
            \vdots & \vdots & \vdots & \vdots & \vdots  \\
            \frac{\p f_{k-1}}{\p x_{i_1}} & \frac{\p f_{k-1}}{\p x_{i_2}} & \cdots & \cdots & \frac{\p f_{k-1}}{\p x_{i_r}} \\
            D(\frac{\p f_{k}}{\p x_{i_1}}) & D(\frac{\p f_{k}}{\p x_{i_2}}) & \cdots & \cdots & D(\frac{\p f_{k}}{\p x_{i_r}}) \\
            \frac{\p f_{k+1}}{\p x_{i_1}} & \frac{\p f_{k+1}}{\p x_{i_2}} & \cdots & \cdots & \frac{\p f_{k+1}}{\p x_{i_r}} \\
            \vdots & \vdots & \vdots & \vdots & \vdots  \\
            \frac{\p f_{r}}{\p x_{i_1}} & \frac{\p f_{r}}{\p x_{i_2}} & \cdots & \cdots & \frac{\p f_{r}}{\p x_{i_n}} \\
        \end{matrix}\right)\\&\equiv 
        \sum_{k=1}^r det \left(\begin{matrix}
            \frac{\p f_1}{\p x_{i_1}} & \frac{\p f_1}{\p x_{i_2}} & \cdots & \cdots & \frac{\p f_1}{\p x_{i_r}} \\
            \vdots & \vdots & \vdots & \vdots & \vdots  \\
            \frac{\p f_{k-1}}{\p x_{i_1}} & \frac{\p f_{k-1}}{\p x_{i_2}} & \cdots & \cdots & \frac{\p f_{k-1}}{\p x_{i_r}} \\
            -\sum_{s=1}^n \frac{\p \alpha_s}{\p x_{i_1}}\frac{\p f_k}{\p x_s}+\sum_{t=1}^r c_k^{(t)}\frac{\p f_t}{\p x_{i_1}} & \cdots & \cdots & \cdots & -\sum_{s=1}^n \frac{\p \alpha_s}{\p x_{i_r}}\frac{\p f_k}{\p x_s}+\sum_{t=1}^r c_k^{(t)}\frac{\p f_t}{\p x_{i_r}} \\
            \frac{\p f_{k+1}}{\p x_{i_1}} & \frac{\p f_{k+1}}{\p x_{i_2}} & \cdots & \cdots & \frac{\p f_{k+1}}{\p x_{i_r}} \\
            \vdots & \vdots & \vdots & \vdots & \vdots  \\
            \frac{\p f_{r}}{\p x_{i_1}} & \frac{\p f_{r}}{\p x_{i_2}} & \cdots & \cdots & \frac{\p f_{r}}{\p x_{i_n}} \\
        \end{matrix}\right)\\&\equiv
        \sum_{k=1}^r c_k^{(k)} \frac{\p (f_1,\cdots,f_r)}{\p (x_{i_1},\cdots,x_{i_r})}-\sum_{k=1}^r \sum_{s=1}^n \frac{\p f_k}{\p x_s} det \left(\begin{matrix}
            \frac{\p f_1}{\p x_{i_1}} & \frac{\p f_1}{\p x_{i_2}} & \cdots & \cdots & \frac{\p f_1}{\p x_{i_r}} \\
            \vdots & \vdots & \vdots & \vdots & \vdots  \\
            \frac{\p f_{k-1}}{\p x_{i_1}} & \frac{\p f_{k-1}}{\p x_{i_2}} & \cdots & \cdots & \frac{\p f_{k-1}}{\p x_{i_r}} \\
             \frac{\p \alpha_s}{\p x_{i_1}} &  \frac{\p \alpha_s}{\p x_{i_2}} & \cdots & \cdots &  \frac{\p \alpha_s}{\p x_{i_r}} \\
            \frac{\p f_{k+1}}{\p x_{i_1}} & \frac{\p f_{k+1}}{\p x_{i_2}} & \cdots & \cdots & \frac{\p f_{k+1}}{\p x_{i_r}} \\
            \vdots & \vdots & \vdots & \vdots & \vdots  \\
            \frac{\p f_{r}}{\p x_{i_1}} & \frac{\p f_{r}}{\p x_{i_2}} & \cdots & \cdots & \frac{\p f_{r}}{\p x_{i_n}} \\
        \end{matrix}\right)\\&\equiv
        -\sum_{s=1}^n\sum_{k=1}^r \frac{\p f_k}{\p x_s} \sum_{t=1}^r (-1)^{k+t} \frac{\p \alpha_s}{\p x_{i_t}}det \left(\begin{matrix}
            \frac{\p f_1}{\p x_{i_1}} & \cdots & \frac{\p f_1}{\p x_{i_{t-1}}} & \frac{\p f_1}{\p x_{i_{t+1}}} & \cdots & \frac{\p f_1}{\p x_{i_r}} \\
            \vdots & \vdots & \vdots & \vdots & \vdots  \\
            \frac{\p f_{k-1}}{\p x_{i_1}} & \cdots & \frac{\p f_{k-1}}{\p x_{i_{t-1}}} & \frac{\p f_{k-1}}{\p x_{i_{t+1}}} & \cdots & \frac{\p f_{k-1}}{\p x_{i_r}} \\
            \frac{\p f_{k+1}}{\p x_{i_1}} & \cdots & \frac{\p f_{k+1}}{\p x_{i_{t-1}}} & \frac{\p f_{k+1}}{\p x_{i_{t+1}}} & \cdots & \frac{\p f_{k+1}}{\p x_{i_r}} \\
            \vdots & \vdots & \vdots & \vdots & \vdots  \\
            \frac{\p f_{r}}{\p x_{i_1}} & \cdots & \frac{\p f_{r}}{\p x_{i_{t-1}}} & \frac{\p f_{r}}{\p x_{i_{t+1}}} & \cdots & \frac{\p f_{r}}{\p x_{i_r}} \\
        \end{matrix}\right)\\&\equiv
        -\sum_{s=1}^n\sum_{t=1}^r\frac{\p \alpha_s}{\p x_{i_t}}\frac{\p (f_1,f_2,\cdots,f_r)}{\p (x_{i_1},\cdots,x_{i_{t-1}},x_{s},x_{i_{t+1}},x_{i_r})}\equiv 0  \ (mod (I,J)).
    \end{align*}
\end{proof}

In particular, for an isolated hypersurface singularity $(V(f),0)\subset (\mathbb{C}^n,0)$, if we take $P=\mathcal{O}_n$, $I=(f)$, then the above theorem tells that every first order derivation on $A$ induces a first order derivation on $\mathcal{O}_n/(f,f_{x_1},\cdots,f_{x_n})$, the Tjurina algebra of the singularity $(V(f),0)$. 

Now a natural question is, which elements in $Der^1(P/(I,\{\frac{\p(f_1,f_2,\cdots,f_r)}{\p(x_{i_1},x_{i_2},\cdots,x_{i_r})}\}_{1\le i_1<\cdots<i_r\le n})$ can be lifted to $Der^1(P/I)$ in theorem \ref{thm2.7} ? For the hypersurface case $I=(f)$, the answer is obvious.

\begin{proposition}\label{prop2.8}
    A derivation $D\in Der^1(P/(f,J(f)))$ is in the image of the map from $Der^1(P/(f))$ in theorem \ref{thm2.7} if and only $D(f)\in (f,J(f)^2)$ in $P$, where $J(f)=(f_{x_1},\cdots,f_{x_n})$ is the Jacobi ideal of $f$ (this notation will be used in the remaining part of the text).
\end{proposition}
\begin{proof}
    If $D=\sum_{i=1}^n {\alpha_i\frac{\p}{\p x_i}}\in Der^1(P/(f,J(f))$ and $D(f)\in (f,J(f)^2)$,  where $\alpha_i$ is a lift of $D(x_i)$ in $P$, $\forall i$. Write $D(f)=fg+\sum_{i=1}^nf_{x_i}g_i, g\in P, g_i\in J(f)$, then $D':=\sum_{i=1}^n (\alpha_i-g_i)\frac{\p}{\p x_i}$ satisfies $D'(f)=gf\in (f)$, $D'\in Der^1(P/(f))$, and the derivation $D'$ induced on $Der^1(P/(f,J(f))$ is $D$.

    The other direction is trivial.
\end{proof}

\begin{lemma}\label{lemma2.9}
   Denote $I_i$ the ideal of $A$ by  $\langle D(x_i): D\in Der^1(A)\rangle$(we will continue using this notation in the whole text). With notations in definition \ref{def2.6}, for $(d_1,d_2\cdots,d_n)\in \mathscr{D}^2(A)$, if $(d_1,d_2,\cdots,d_n)$ is in the image of $der^2(A)$ under $\theta_2$, then $d_i(x_i)\in I_i^2$.
\end{lemma}
\begin{proof}
    As $\theta_2(Der^1(A))=0$, $der^2(A)=Der^1(A)+Der^1(A)Der^1(A)$, we only need to show that, for $D_1,D_2\in Der^1(A)$, if we denote $\theta_2(D_1D_2)=(d_1,d_2,\cdots,d_n)$, then $d_i(x_i)\in I_i^2$.
    
    Let $D_1=\sum_{i=1}^n \alpha_i\frac{\p}{\p x_i}$, $D_2=\sum_{i=1}^n \beta_i \frac{\p}{\p x_i}$, then $D_1D_2=\sum_{i,j=1}^n (\alpha_i\beta_j \frac{\p^2}{\p x_i \p x_j}+\alpha_i\frac{\p \beta_j}{\p x_i}\frac{\p}{\p x_j})$, $\langle D,x_i\rangle=\sum_{j=1}^n (\alpha_i\beta_j+\alpha_j\beta_i)\frac{\p}{\p x_j}+\sum_{j=1}^n \alpha_j\frac{\p \beta_i}{\p x_j}$, thus $d_i=\sum_{j=1}^n (\alpha_i\beta_j+\alpha_j\beta_i)\frac{\p}{\p x_j}$, and $d_i(x_i)=2\alpha_i\beta_i\in I_i^2$
\end{proof}

From now on, we concentrate on dealing with the cases of isolated hypersurface singularities  $(V(f),0)$. 

\begin{theorem}\label{thm2.10}
Suppose $A=P/I$ and $I$ is principal. Then the sequence
$$
0 \rightarrow \operatorname{Der}^1(A) \to  \operatorname{Der}^2(A) \stackrel{\theta_2}{\longrightarrow} \mathscr{D}^2(A) \rightarrow 0
$$
is exact.
\end{theorem}
\begin{proof}
See Theorem 2.13 in \cite{Singh}. 
\end{proof}

\begin{proposition}\label{prop2.11}\cite{Kantor}
    For $(V(f),0)$ an isolated weighted homogeneous hypersurface singularity with weight type $(w_1,\cdots,w_n;1)$, let $A=P/(f)$, then $Der^1(A)$ is generated by the Euler derivation $E=\sum_{i=1}^n w_ix_i\frac{\p}{\p x_i}$ and the Hamiltonian derivations $D_{ij}:=f_{x_i}\frac{\p}{\p x_j}-f_{x_j}\frac{\p}{\p x_i}.$ Moreover, $I_i=(x_i,f_1,\cdots,f_{i-1},f_{i+1},\cdots,f_n)$.
\end{proposition}

\begin{theorem}\label{thm2.12}
	For $n$ derivations $d_1,\cdots,d_n\in Der^1(A)$, if $d_i(x_j)-d_{j}(x_i)\in J(f)\  \forall i,j$, then $\exists$ $(d'_1,\cdots,d'_n)\in \mathscr{D}^2(A)$, such that for each $i$, $d'_{i}-d_i$ is an $A$-linear combination of $D_{kl}$'s, $(k,l=1,\cdots,n)$.
\end{theorem}
\begin{proof}
    Let $d_i(x_j)-d_j(x_i)=\sum_{s=1}^n c_{ij}^{(s)}f_{x_s}, 1\le i<j \le n$, first we may assume that $c_{ij}^{(i)}=c_{ij}^{(j)}=0$, because we can replace $d_i$ by $d_i-c_{ij}^{(j)}D_{ij}$, and replace $d_j$ by $d_j-c_{ij}^{(i)}D_{ij}$.

    Now for each pair $(i,j)$, $1\le i<j\le n$, we define an operation $P_{ij}$ on the tuple $(d_1,d_2,\cdots,d_n)$ as following: define
    $P_{ij}(d_1,\cdots,d_n)=(d_1',\cdots,d_n')$ by 
    $$d_i':=d_i-\sum_{s=1}^{i-1}\frac{1}{2}c_{ij}^{(s)}D_{sj}-\sum_{s=i+1}^{j-1}\frac{1}{2}c_{ij}^{(s)}D_{sj}+\sum_{s=j+1}^n \frac{1}{2}c_{ij}^{(s)}D_{js} \ ,$$
    $$d_j':=d_j+\sum_{s=1}^{i-1}\frac{1}{2}c_{ij}^{(s)}D_{si}-\sum_{s=i+1}^{j-1} \frac{1}{2}c_{ij}^{(s)}D_{is}-\sum_{s=j+1}^n \frac{1}{2}c_{ij}^{(s)}D_{is} \ ,$$
    $$d_t':=d_t-\frac{1}{2}c_{ij}^{(t)}D_{ij} \ .$$
    Then we see that $d_s'(x_t)-d_t'(x_s)=d_s(x_t)-d_t(x_s), \forall s,t\ne i,j$, and $d_i'(x_s)-d_s'(x_i)=d_i(x_s)+\frac{1}{2}c_{ij}^{(s)}f_{x_j}-(d_s(x_i)+\frac{1}{2}c_{ij}^{(s)}f_{x_j})=d_i(x_s)-d_s(x_i)$, $d_j'(x_s)-d_s'(x_j)=d_j(x_s)-\frac{1}{2}c_{ij}^{(s)}f_{x_i}-(d_s(x_j)-\frac{1}{2}c_{ij}^{(s)}f_{x_i})=d_j(x_s)-d_s(x_j)$, $\forall s\ne i,j$, and $d_i'(x_j)-d_j'(x_i)=0$. The operation $P_{ij}$ only changes $d_i(x_j)-d_j(x_i)$ to $0$ and keeps the other $d_s(x_t)-d_t(x_s), \{s,t\}\ne \{i,j\}$ unchanged. Therefore, after taking the $\frac{n(n-1)}{2}$ times operations $P_{ij}$'s $1\le i< j\le n$(the order of the operations $P_{ij}$ is arbitrary), the tuple $(d_1',\cdots,d_n')$ we obtained satisfies the condition in the statement of the theorem.
\end{proof}

\begin{corollary}\label{cor2.13}
    If there exists $(d_1,\cdots,d_n)\in Der^1(A)^{\oplus n}$, satisfying $d_i(x_j)-d_j(x_i)\in J(f)$, and there exists some $i$, such that $d_i(x_i)\notin (I_i^2,f_{x_1},\cdots,f_{x_{i-1}},f_{x_{i+1}},\cdots,f_{x_n})$, then $(d_1,\cdots,d_n)$ does not lie in the image of $der^2(A)$ under $\theta_2$, and $der^2(A)\ne Der^2(A)$ follows from $\theta_2$ is surjective in theorem \ref{thm2.10}.
\end{corollary}

Combining with theorem \ref{thm2.7}, proposition \ref{prop2.8} and corollary \ref{cor2.13}, we obtain the following proposition for proving the Nakai Conjecture.

\begin{proposition}\label{prop2.14}
    Let $(V(f),0)$ be an isolated hypersurface singularity, if there exists a symmetric matrix $(\beta_{ij})_{n\times n}\in M_n(\mathcal{O}_n/(f,J(f))$, satisfying $(\widetilde{\beta_{ij}})\cdot (f_{x_1},\cdots,f_{x_n})^{T}\equiv \textbf{0} \mod(f,J(f)^2)$, where $(\widetilde{\beta_{ij}})$ is an arbitrary lift of $(\beta_{ij})$ in $M_n(\mathcal{O}_n)$, and there exists $1\le i\le n$, such that $\beta_{ii}\notin I_i^2 \mod J(f)$, then the Nakai Conjecture holds for the singularity $(V(f),0)$.
\end{proposition}
\begin{proof}
    From theorem \ref{thm2.7} and proposition \ref{prop2.8}, there exists derivations $d_1,\cdots,d_n\in Der^1(\mathcal{O}_n/(f))$, such that $d_i(x_j)\equiv \beta_{ij} \mod J(f)$, as $\beta_{ij}=\beta_{ji}$, $d_i(x_j)-d_j(x_i)\in J(f), \forall 1\le i,j\le n$. And there exists some $1\le i\le n$, $d_i(x_i)\notin I_i^2+J(f)$, so $d_i(x_i)\notin (I_i^2,f_{x_1},\cdots,f_{x_{i-1}},f_{x_{i+1}},\cdots,f_{x_n})$, the assertion follows from corollary \ref{cor2.13}.
\end{proof}
\textbf{Proof of the main theorem}

Now we begin to prove our main theorem. As the simple isolated hypersurface singularities (ADE singularities) are weighted homogeneous fewnomial, the Nakai Conjecture in which cases has been proved in \cite{XYZZ} as following. In the next two sections, we only consider the cases of modality $=1$ and $2$ isolated hypersurface singularities. 

\begin{theorem}\cite{XYZZ}\label{thm2.15}
    Let $A=P/I$ where $P=k[x_1,\cdots,x_n]$, $I=(f)$ with $(V(f),0)$ a weighted homogeneous fewnomial isolated hypersurface singularity, then $der^2(A)\ne Der^2(A)$. 
\end{theorem}

\section{The unimodal isolated hypersurface singularity case}\label{sec 3}
The unimodal isolated hypersurface singularities were classified into simple elliptic singularities, cusp hypersurface singularities, and 14 exceptional families in \cite{Arnold} as following.

\begin{center}
    \begin{longtable}{|p{1.2cm}|p{4.5cm}|p{3cm}|p{4cm}|}
    \caption{14 exceptional families of unimodal singularities }\label{table-0}\\
    \hline
    symbol & normal forms & Milnor number $\mu$ & Tjurina number $\tau$\\
    \hline 
    $E_{12}$ & $x^3+y^7+axy^5$ & $12$ & $12(a=0), 11(a\ne 0)$\\
    \hline
    $E_{13}$ & $x^3+xy^5+ay^8$ & $13$ & $13(a=0), 12(a\ne 0)$\\
    \hline
    $E_{14}$ & $x^3+y^8+axy^6$ & $14$ & $14(a=0), 13(a\ne 0)$\\
    \hline
    $Z_{11}$ & $x^3y+y^5+axy^4$ & $11$ & $11(a=0), 10(a\ne 0)$\\
    \hline
    $Z_{12}$ & $x^3y+xy^4+ax^2y^3$ & $12$ & $12(a=0), 11(a\ne 0)$\\
    \hline
    $Z_{13}$ & $x^3y+y^6+axy^5$ & $13$ & $13(a=0), 12(a\ne 0)$\\
    \hline
    $W_{12}$ & $x^4+y^5+ax^2y^3$ & $12$ & $12(a=0), 11(a\ne 0)$\\
    \hline
    $W_{13}$ & $x^4+xy^4+ay^6$ & $13$ & $13(a=0), 12(a\ne 0)$\\
    \hline
    $Q_{10}$ & $x^3+y^4+yz^2+axy^3$ & $10$ & $10(a=0), 9(a\ne 0)$\\
    \hline
    $Q_{11}$ & $x^3+y^2z+xz^3+az^5$ & $11$ & $11(a=0), 10(a\ne 0)$\\
    \hline
    $Q_{12}$ & $x^3+y^5+yz^2+axy^4$ & $12$ & $12(a=0), 11(a\ne 0)$\\
    \hline
    $S_{11}$ & $x^4+y^2z+xz^2+ax^3z$ & $11$ & $11(a=0), 10(a\ne 0)$\\
    \hline
    $S_{12}$ & $x^2y+y^2z+xz^3+az^5$ & $12$ & $12(a=0), 11(a\ne 0)$\\
    \hline
    $U_{12}$ & $x^3+y^3+z^4+axyz^2$ & $12$ & $12(a=0), 11(a\ne 0)$\\
    \hline
    \end{longtable}
\end{center}

In this section, for an isolated hypersurface singularity $(V(f),0)\subset (\mathbb{C}^n,0)$, $A=\mathcal{O}_n/(f)$, with notations as in section \ref{sec 2}, we will verify the Nakai Conjecture for unimodal isolated hypersurface singularities. As the cusp hypersurface singularities case was solved in \cite{XYZZ}. We only need to deal with simple elliptic singularities and the exceptional families of singularities listed as above.

\subsection{Cases of simple elliptic singularities}
\ 
\newline
\indent
Simple elliptic singularities can be classfied into three types: $\widetilde{E}_6, \widetilde{E}_7, \widetilde{E}_8$, they are all one parameter families defined by a weighted homogeneous polynomial $f_t$.

\begin{proposition}(\cite{XYZZ})\label{prop3.1}
    Let $f\in \mathbb{C}[x_1\cdots,x_n]$ be a weighted homogeneous polynomial with weight type $(w_1,\cdots,w_n;1)$, which defines an isolated hypersurface singularity at the origin, if there exists some hyperplane section $\{x_s=0\}\cap V(f)$ which has an isolated singularity at the origin, then for $A=\mathcal{O}_n/(f)$, we have $Der^2(A)\ne der^2(A)$.
\end{proposition}
\begin{proof}
    Denote $H_{ij}$ to be the algebra co-factor of $\frac{\p^2 f}{\p x_i \p x_j}$ in $Hess(f)$, from theorem 3.7 in \cite{XYZZ}, we have $w_ix_iH_{jk}-w_jx_jH_{ik}\in J(f), \forall 1\le i,j,k\le n$. Now define $d_i:=H_{is}\cdot E$, $1\le i\le n$, where $E=\sum_{k=1}^n w_k\frac{\p}{\p x_k}$ is the Euler derivation, then $\forall 1\le i,j\le n$, $d_i(x_j)-d_j(x_i)=w_jx_jH_{is}-w_ix_iH_{js}\in J(f)$.

    As $\{x_s=0\}\cap V(f)$ has an isolated singularity at the origin, the polynomial $f|_{x_s=0}$ defines an isolated hypersurface singularity in $\mathbb{C}^{n-1}$, so the $(n-1)$ polynomials $\frac{\p f}{\p x_i} (i\ne s)$ form a regular sequence in $\mathcal{O}_{n-1}=\mathcal{O}_n/(x_s)$, then $$H_{ss}=\frac{\p (f_{x_1},\cdots,f_{x_{s-1}},f_{x_{s+1}},\cdots,f_{x_n})}{\p (x_1,\cdots,x_{s-1},x_{s+1},\cdots,x_n)}\notin (f_{x_1},\cdots,f_{x_{s-1}},f_{x_{s+1}},\cdots,f_{x_n})\mathcal{O}_{n-1}.$$
    $x_s$ is regular in $\mathcal{O}_n$, so $x_sH_{ss}\notin (x_s^2, f_{x_1},\cdots,f_{x_{s-1}},f_{x_{s+1}},\cdots,f_{x_n})\mathcal{O}_{n}$, i.e. $d_s(x_s)=w_sx_sH_{ss}\notin (I_s^2,f_{x_1},\cdots,f_{x_{s-1}}, f_{x_{s+1}},\cdots,f_{x_n})$ from proposition \ref{prop2.11},  and $Der^2(A)\ne der^2(A)$ follows from theorem \ref{thm2.10} and corollary \ref{cor2.13}.
\end{proof}

For a $\widetilde{E}_6$ singularity defined by $f_t=x^3+y^3+z^3+txyz$ with $t^3+27\ne 0$, the hyperplane section $V(f_t)\cap \{x=0\}=\{(y,z) | y^3+z^3=0\}$ has an isolated singularity at the origin, the Nakai Conjecture holds by proposition \ref{prop3.1}.

For a $\widetilde{E}_7$ singularity defined by $f_t=x^4+y^4+tx^2y^2$ with $t^2\ne 4$, the hyperplane section $V(f_t)\cap \{x=0\}=\{y^4=0\}\subset \mathbb{C}$ has an isolated singularity at the origin, the Nakai Conjecture holds by proposition \ref{prop3.1}.

For a $\widetilde{E}_8$ singularity defined by $f_t=x^3+y^6+tx^2y^2$ with $t^3+\frac{27}{4}\ne 0$, the hyperplane section $V(f_t)\cap \{x=0\}=\{y^6=0\}\subset \mathbb{C}$ has an isolated singularity at the origin, the Nakai Conjecture holds by proposition \ref{prop3.1}.

\subsection{Cases of exceptional families of unimodal isolated hypersurface singularities}
\ 
\newline
\indent
Now we treat the exceptional families of unimodal isolated hypersurface singularities. In table \ref{table-0}, we see that if $a=0$, then they are reduced to weighted homogeneous fewnomial singularities (see \ref{thm2.15}), so we only consider the case of $a\ne 0$ in table \ref{table-0}.

1. $E_{12}$ singularity

For an $E_{12}$ singularity $(V(f),0)$ defined by $f=x^3+y^7+axy^5$ for some $a\ne 0$, $f_x=3x^2+ay^5$, $f_y=7y^6+5axy^4$. As $\frac{1}{3}xf_x+\frac{1}{7}yf_y-f=\frac{1}{21}axy^5$, $R=\mathcal{O}_2/(f,f_x,f_y)=\mathbb{C}\{x,y\}/(x^3,y^7,xy^5, 3x^2+ay^5,7y^6+5axy^4)$, it has a $\mathbb{C}$-basis: $\{1,y,y^2,y^3,y^4,y^5,y^6,x,xy,xy^2,xy^3\}$, $dim_{\mathbb{C}} R=11$. 

As a derivation $\alpha_1\p_x+\alpha_2\p_y\in Der^1(R)$ if and only if $Hess(f)\cdot (\alpha_1,\alpha_2)^T=0$ in $R$, by solving the equation $Hess(f)\cdot (\alpha_1,\alpha_2)^T=0$ in $R$, we obtain that 
\begin{align*}
\begin{pmatrix}
\alpha_1\\
\alpha_2
\end{pmatrix}&\in Span_{\mathbb{C}}\{
\begin{pmatrix} ay^4 \\ -\frac{6}{5}x \end{pmatrix},
\begin{pmatrix} y^5 \\ 0  \end{pmatrix},
\begin{pmatrix} y^6 \\ 0  \end{pmatrix},
\begin{pmatrix} xy \\ -\frac{2}{7}ax \end{pmatrix},
\begin{pmatrix} xy^2 \\ 0 \end{pmatrix},
\begin{pmatrix} xy^3 \\ 0 \end{pmatrix},
\begin{pmatrix} 0 \\ \frac{5}{7}ax+y^2 \end{pmatrix},
\begin{pmatrix} 0 \\ y^3  \end{pmatrix},
\begin{pmatrix} 0 \\ y^4  \end{pmatrix},\\&
\begin{pmatrix} 0 \\ y^5  \end{pmatrix},
\begin{pmatrix} 0 \\ y^6  \end{pmatrix},
\begin{pmatrix} 0 \\ xy  \end{pmatrix},
\begin{pmatrix} 0 \\ xy^2 \end{pmatrix},
\begin{pmatrix} 0 \\ xy^3 \end{pmatrix}
\}.
\end{align*}
So $Der^1(R)=Span_{\mathbb{C}}\{ay^4\p_x-\frac{6}{5}x\p_y, y^5\p_x, y^6\p_x, xy\p_x-\frac{2}{7}ax\p_y, xy^2\p_x, xy^3\p_x, (\frac{5}{7}ax+y^2)\p_y, y^3\p_y, \\y^4\p_y, y^5\p_y, y^6\p_y, xy\p_y, xy^2\p_y, xy^3\p_y\}$.
We see that the ideal $I_2\subset \mathfrak{m}+(f,J(f))$, $I_2^2\subset \mathfrak{m}^2+(f,J(f))$. Now we define $(\beta_{ij})\in M_2(R)$ by 
$$\beta_{11}= -\frac{1}{21}a^2 xy^3 + a y^5, $$
 $$\beta_{12}=\beta_{21}= -\frac{1}{1176}a^3 xy^2 - \frac{9}{8} xy -\frac{15}{392}a^2 y^4, $$
 $$\beta_{22}= -\frac{25}{697662}a^7 xy^2 + \frac{1}{294}a^4 xy + \frac{9}{392}ax - \frac{5}{6174}a^6 y^4 +
\frac{1}{196}a^3 y^3 - \frac{27}{56}y^2, $$
then the lifting $(\widetilde{\beta_{ij}})\in M_2(\mathcal{O}_2)$ satisfies$$
\begin{pmatrix}
\widetilde{\beta_{11}} & \widetilde{\beta_{12}} \\
\widetilde{\beta_{21}} & \widetilde{\beta_{22}}
\end{pmatrix}
\begin{pmatrix}
f_{x}\\
f_{y}
\end{pmatrix}
\equiv
\begin{pmatrix}
0\\
0
\end{pmatrix} \ \ \ (  \textup{mod} \ (f,J(f)^2)  ),$$ and $\widetilde{\beta_{22}}\notin \mathfrak{m}^2+(f,J(f))$, so $\beta_{22}\notin I_2^2 \mod J(f)$. By proposition \ref{prop2.14}, the Nakai Conjecture holds.
\\

Now we treat the other singularities in the same way.

2. $E_{13}$ singularity 

For an $E_{13}$ singularity $(V(f),0)$ defined by $f=x^3+xy^5+ay^8$ for some $a\ne 0$, $f_x=3x^2+y^5$, $f_y=5xy^4+8ay^7$. As $\frac{1}{3}xf_x+\frac{1}{8}yf_y-f=-\frac{1}{24}xy^5$, $R=\mathcal{O}_2/(f,f_x,f_y)=\mathbb{C}\{x,y\}/(x^3,y^8,xy^5, 3x^2+y^5,5xy^4+8ay^7)$, it has a $\mathbb{C}$-basis: $\{1,y,y^2,y^3,y^4,y^5,y^6,y^7,x,xy,xy^2,xy^3\}$, $dim_{\mathbb{C}} R=12$. 

By solving the equation $Hess(f)\cdot (\alpha_1,\alpha_2)^T=0$ in $R$, we obtain that 
\begin{align*}
\begin{pmatrix}
\alpha_1\\
\alpha_2
\end{pmatrix}&\in Span_{\mathbb{C}}\{
\begin{pmatrix} y^4 \\ -\frac{6}{5}x \end{pmatrix},
\begin{pmatrix} y^5 \\ 0  \end{pmatrix},
\begin{pmatrix} y^6 \\ 0  \end{pmatrix},
\begin{pmatrix} y^7 \\ 0  \end{pmatrix},
\begin{pmatrix} xy \\ \frac{2}{5}y^2 \end{pmatrix},
\begin{pmatrix} xy^2 \\ \frac{2}{5}y^3 \end{pmatrix},
\begin{pmatrix} xy^3 \\ 0 \end{pmatrix},
\begin{pmatrix} 0 \\ \frac{8}{5}x+ay^3 \end{pmatrix},
\begin{pmatrix} 0 \\ y^4  \end{pmatrix},\\&
\begin{pmatrix} 0 \\ y^5  \end{pmatrix},
\begin{pmatrix} 0 \\ y^6  \end{pmatrix},
\begin{pmatrix} 0 \\ y^7  \end{pmatrix},
\begin{pmatrix} 0 \\ xy  \end{pmatrix},
\begin{pmatrix} 0 \\ xy^2 \end{pmatrix},
\begin{pmatrix} 0 \\ xy^3 \end{pmatrix}
\}.
\end{align*}
We see that the ideal $I_2\subset \mathfrak{m}+(f,J(f))$, $I_2^2\subset \mathfrak{m}^2+(f,J(f))$. Now we define $(\beta_{ij})\in M_2(R)$ by
$$\beta_{11}=\frac{3}{10}a xy^3 + y^5, $$
$$\beta_{12}=\beta_{21}= -\frac{3213}{12500}a^4 xy^3 -\frac{153}{250}a^2 xy^2 - \frac{6}{5} xy + \frac{6}{25}a y^4, $$
$$\beta_{22}=\frac{45734787}{15625000}a^7 xy^3 + \frac{2177847}{312500}a^5 xy^2 -
\frac{6453}{6250}a^3 xy - \frac{18}{125}ax - \frac{24192}{15625}a^4 y^4 - \frac{9}{625}a^2 y^3 - \frac{12}{25}y^2, $$
then the lifting $(\widetilde{\beta_{ij}})\in M_2(\mathcal{O}_2)$ satisfies$$
\begin{pmatrix}
\widetilde{\beta_{11}} & \widetilde{\beta_{12}} \\
\widetilde{\beta_{21}} & \widetilde{\beta_{22}}
\end{pmatrix}
\begin{pmatrix}
f_{x}\\
f_{y}
\end{pmatrix}
\equiv
\begin{pmatrix}
0\\
0
\end{pmatrix} \ \ \ (  \textup{mod} \ (f,J(f)^2)  ),$$ and $\widetilde{\beta_{22}}\notin \mathfrak{m}^2+(f,J(f))$, so $\beta_{22}\notin I_2^2 \mod J(f)$. By proposition \ref{prop2.14}, the Nakai Conjecture holds.
\\

3. $E_{14}$ singularity 

For an $E_{14}$ singularity $(V(f),0)$ defined by $f=x^3+y^8+axy^6$ for some $a\ne 0$, $f_x=3x^2+ay^6$, $f_y=8y^7+6axy^5$. As $\frac{1}{3}xf_x+\frac{1}{8}yf_y-f=\frac{1}{12}axy^6$, $R=\mathcal{O}_2/(f,f_x,f_y)=\mathbb{C}\{x,y\}/(x^3,y^8,xy^6, 3x^2+ay^6,8y^7+6axy^5)$, it has a $\mathbb{C}$-basis: $\{1,y,y^2,y^3,y^4,y^5,y^6,y^7,x,xy,xy^2,xy^3,xy^4\}$, $dim_{\mathbb{C}} R=13$. 

By solving the equation $Hess(f)\cdot (\alpha_1,\alpha_2)^T=0$ in $R$, we obtain that 
\begin{align*}
\begin{pmatrix}
\alpha_1\\
\alpha_2
\end{pmatrix}&\in Span_{\mathbb{C}}\{
\begin{pmatrix} ay^5 \\ -x \end{pmatrix},
\begin{pmatrix} y^6 \\ 0  \end{pmatrix},
\begin{pmatrix} y^7 \\ 0  \end{pmatrix},
\begin{pmatrix} xy \\ -\frac{1}{4}ax \end{pmatrix},
\begin{pmatrix} xy^2 \\ 0 \end{pmatrix},
\begin{pmatrix} xy^3 \\ 0 \end{pmatrix},
\begin{pmatrix} xy^4 \\ 0 \end{pmatrix},
\begin{pmatrix} 0 \\ \frac{3}{4}ax+y^2 \end{pmatrix},
\begin{pmatrix} 0 \\ y^3  \end{pmatrix},\\&
\begin{pmatrix} 0 \\ y^4  \end{pmatrix},
\begin{pmatrix} 0 \\ y^5  \end{pmatrix},
\begin{pmatrix} 0 \\ y^6  \end{pmatrix},
\begin{pmatrix} 0 \\ y^7  \end{pmatrix},
\begin{pmatrix} 0 \\ xy  \end{pmatrix},
\begin{pmatrix} 0 \\ xy^2 \end{pmatrix},
\begin{pmatrix} 0 \\ xy^3 \end{pmatrix},
\begin{pmatrix} 0 \\ xy^4 \end{pmatrix}
\}.
\end{align*}
We see that the ideal $I_2\subset \mathfrak{m}+(f,J(f))$, $I_2^2\subset \mathfrak{m}^2+(f,J(f))$. Now we define $(\beta_{ij})\in M_2(R)$ by
 $$\beta_{11}= -\frac{1}{12}a^2 xy^4 + a y^6, $$
 $$\beta_{12}=\beta_{21}= -\frac{9}{10} xy - \frac{9}{160}a^2 y^5, $$
 $$\beta_{22}= \frac{27}{52480}a^4 xy^2 + \frac{9}{320}ax +\frac{9}{1280}a^3 y^4 - \frac{27}{80}y^2, $$
then the lifting $(\widetilde{\beta_{ij}})\in M_2(\mathcal{O}_2)$ satisfies$$
\begin{pmatrix}
\widetilde{\beta_{11}} & \widetilde{\beta_{12}} \\
\widetilde{\beta_{21}} & \widetilde{\beta_{22}}
\end{pmatrix}
\begin{pmatrix}
f_{x}\\
f_{y}
\end{pmatrix}
\equiv
\begin{pmatrix}
0\\
0
\end{pmatrix} \ \ \ (  \textup{mod} \ (f,J(f)^2)  ),$$ and $\widetilde{\beta_{22}}\notin \mathfrak{m}^2+(f,J(f))$, so $\beta_{22}\notin I_2^2 \mod J(f)$. By proposition \ref{prop2.14}, the Nakai Conjecture holds.
\\

4. $Z_{11}$ singularity

For a $Z_{11}$ singularity $(V(f),0)$ defined by $f=x^3y+y^5+axy^4$ for some $a\ne 0$, $f_x=3x^2y+ay^4$, $f_y=x^3+5y^4+4axy^3$. As $\frac{4}{15}xf_x+\frac{1}{5}yf_y-f=\frac{1}{15}axy^4$, $R=\mathcal{O}_2/(f,f_x,f_y)=\mathbb{C}\{x,y\}/(x^3y,y^5,xy^4, 3x^2y+ay^4,x^3+5y^4+4axy^3)$, it has a $\mathbb{C}$-basis: $\{1,y,y^2,y^3,y^4,x,xy,xy^2,xy^3,\\ x^2\}$, $dim_{\mathbb{C}} R=10$. 

By solving the equation $Hess(f)\cdot (\alpha_1,\alpha_2)^T=0$ in $R$, we obtain that 
\begin{align*}
\begin{pmatrix}
\alpha_1\\
\alpha_2
\end{pmatrix}&\in Span_{\mathbb{C}}\{
\begin{pmatrix} y^3 \\ 0  \end{pmatrix},
\begin{pmatrix} y^4 \\ 0  \end{pmatrix},
\begin{pmatrix} xy \\ 0 \end{pmatrix},
\begin{pmatrix} xy^2 \\ 0 \end{pmatrix},
\begin{pmatrix} xy^3 \\ 0 \end{pmatrix},
\begin{pmatrix} x^2 \\ 0 \end{pmatrix},
\begin{pmatrix} 0 \\ y^2  \end{pmatrix},
\begin{pmatrix} 0 \\ y^3  \end{pmatrix},
\begin{pmatrix} 0 \\ y^4  \end{pmatrix},
\begin{pmatrix} 0 \\ xy  \end{pmatrix},
\begin{pmatrix} 0 \\ xy^2 \end{pmatrix},\\&
\begin{pmatrix} 0 \\ xy^3 \end{pmatrix},
\begin{pmatrix} 0 \\ x^2 \end{pmatrix}
\}.
\end{align*}
We see that the ideal $I_1\subset \mathfrak{m}^2+(f,J(f))$, $I_1^2\subset \mathfrak{m}^4+(f,J(f))$. Now we define $(\beta_{ij})\in M_2(R)$ by
$$\beta_{11}=-12x^2 - \frac{5}{12}a^2xy^2 + ay^3, $$
$$\beta_{12}=\beta_{21}= -\frac{1}{48}a^3 xy^2 - 9 xy -\frac{9}{16}a^2 y^3, $$
$$\beta_{22}= \frac{5}{492}a^4 xy^2 +\frac{9}{16}a xy + \frac{1}{8}a^3 y^3 - \frac{27}{4}y^2, $$ 
then the lifting $(\widetilde{\beta_{ij}})\in M_2(\mathcal{O}_2)$ satisfies$$
\begin{pmatrix}
\widetilde{\beta_{11}} & \widetilde{\beta_{12}} \\
\widetilde{\beta_{21}} & \widetilde{\beta_{22}}
\end{pmatrix}
\begin{pmatrix}
f_{x}\\
f_{y}
\end{pmatrix}
\equiv
\begin{pmatrix}
0\\
0
\end{pmatrix} \ \ \ (  \textup{mod} \ (f,J(f)^2)  ),$$ and $\widetilde{\beta_{11}}\notin \mathfrak{m}^4+(f,J(f))$, so $\beta_{11}\notin I_1^2 \mod J(f)$.
By proposition \ref{prop2.14}, the Nakai Conjecture holds.\\

5. $Z_{12}$ singularity

For a $Z_{12}$ singularity $(V(f),0)$ defined by $f=x^3y+xy^4+ax^2y^3$ for some $a\ne 0$, $f_x=3x^2y+y^4+2axy^3$, $f_y=x^3+4xy^3+3ax^2y^2$. As $\frac{2}{7}xf_x+\frac{1}{7}yf_y-f=-\frac{1}{7}xy^4$, $\frac{1}{5}xf_x+\frac{1}{5}yf_y-f=-\frac{1}{5}x^3y$,  $R=\mathcal{O}_2/(f,f_x,f_y)=\mathbb{C}\{x,y\}/(x^3y, xy^4, x^2y^3,  3x^2y+y^4+2axy^3,x^3+4xy^3+3ax^2y^2)$, it has a $\mathbb{C}$-basis: $\{1,y,y^2,y^3,y^4,y^5,x,xy,xy^2,xy^3,x^2\}$, $dim_{\mathbb{C}} R=11$. 

By solving the equation $Hess(f)\cdot (\alpha_1,\alpha_2)^T=0$ in $R$, we obtain that 
\begin{align*}
\begin{pmatrix}
\alpha_1\\
\alpha_2
\end{pmatrix}&\in Span_{\mathbb{C}}\{
\begin{pmatrix} y^2 \\ \frac{3}{4}x-\frac{11}{12}ay^2 \end{pmatrix},
\begin{pmatrix} y^3 \\ 0  \end{pmatrix},
\begin{pmatrix} y^4 \\ 0  \end{pmatrix},
\begin{pmatrix} y^5 \\ 0  \end{pmatrix},
\begin{pmatrix} xy \\ \frac{2}{3}y^2 \end{pmatrix},
\begin{pmatrix} xy^2 \\ 0 \end{pmatrix},
\begin{pmatrix} xy^3 \\ 0 \end{pmatrix},
\begin{pmatrix} x^2 \\ 0 \end{pmatrix},
\begin{pmatrix} 0 \\ y^3  \end{pmatrix},\\&
\begin{pmatrix} 0 \\ y^4  \end{pmatrix},
\begin{pmatrix} 0 \\ y^5  \end{pmatrix},
\begin{pmatrix} 0 \\ xy  \end{pmatrix},
\begin{pmatrix} 0 \\ xy^2 \end{pmatrix},
\begin{pmatrix} 0 \\ xy^3 \end{pmatrix},
\begin{pmatrix} 0 \\ x^2 \end{pmatrix}
\}.
\end{align*}
We see that the ideal $I_1\subset \mathfrak{m}^2+(f,J(f))$, $I_1^2\subset \mathfrak{m}^4+(f,J(f))$. Now we define $(\beta_{ij})\in M_2(R)$ by
$$\beta_{11}= -6x^2 + axy^2, $$
$$\beta_{12}=\beta_{21}= \frac{1}{27}a^4xy^3 - \frac{1}{3}a^2xy^2 - 4xy,$$
 $$\beta_{22}= -\frac{23}{729}a^5xy^3 +\frac{11}{81}a^3xy^2 + \frac{4}{9}axy + \frac{2}{27}a^2y^3 - \frac{8}{3}y^2, $$
then the lifting $(\widetilde{\beta_{ij}})\in M_2(\mathcal{O}_2)$ satisfies$$
\begin{pmatrix}
\widetilde{\beta_{11}} & \widetilde{\beta_{12}} \\
\widetilde{\beta_{21}} & \widetilde{\beta_{22}}
\end{pmatrix}
\begin{pmatrix}
f_{x}\\
f_{y}
\end{pmatrix}
\equiv
\begin{pmatrix}
0\\
0
\end{pmatrix} \ \ \ (  \textup{mod} \ (f,J(f)^2)  ),$$ and $\widetilde{\beta_{11}}\notin \mathfrak{m}^4+(f,J(f))$, so $\beta_{11}\notin I_1^2 \mod J(f)$.
By proposition \ref{prop2.14}, the Nakai Conjecture holds.\\

6. $Z_{13}$ singularity

For a $Z_{13}$ singularity $(V(f),0)$ defined by $f=x^3y+y^6+axy^5$ for some $a\ne 0$, $f_x=3x^2y+ay^5$, $f_y=x^3+6y^5+5axy^4$. As $\frac{5}{18}xf_x+\frac{1}{6}yf_y-f=\frac{1}{9}axy^5$, $R=\mathcal{O}_2/(f,f_x,f_y)=\mathbb{C}\{x,y\}/(x^3y, xy^5, y^6,  3x^2y+ay^5,x^3+6y^5+5axy^4)$, it has a $\mathbb{C}$-basis: $\{1,y,y^2,y^3,y^4,y^5,x,xy,xy^2,\\xy^3,xy^4,x^2\}$, $dim_{\mathbb{C}} R=12$. 

By solving the equation $Hess(f)\cdot (\alpha_1,\alpha_2)^T=0$ in $R$, we obtain that 
\begin{align*}
\begin{pmatrix}
\alpha_1\\
\alpha_2
\end{pmatrix}&\in Span_{\mathbb{C}}\{
\begin{pmatrix} y^4 \\ 0  \end{pmatrix},
\begin{pmatrix} y^5 \\ 0  \end{pmatrix},
\begin{pmatrix} xy \\ 0 \end{pmatrix},
\begin{pmatrix} xy^2 \\ 0 \end{pmatrix},
\begin{pmatrix} xy^3 \\ 0 \end{pmatrix},
\begin{pmatrix} xy^4 \\ 0 \end{pmatrix},
\begin{pmatrix} x^2 \\ 0 \end{pmatrix},
\begin{pmatrix} 0 \\ y^2  \end{pmatrix},
\begin{pmatrix} 0 \\ y^3  \end{pmatrix},
\begin{pmatrix} 0 \\ y^4  \end{pmatrix},\\&
\begin{pmatrix} 0 \\ y^5  \end{pmatrix},
\begin{pmatrix} 0 \\ xy  \end{pmatrix},
\begin{pmatrix} 0 \\ xy^2 \end{pmatrix},
\begin{pmatrix} 0 \\ xy^3 \end{pmatrix},
\begin{pmatrix} 0 \\ xy^4 \end{pmatrix},
\begin{pmatrix} 0 \\ x^2 \end{pmatrix}
\}.
\end{align*}
We see that the ideal $I_1\subset \mathfrak{m}^2+(f,J(f))$, $I_1^2\subset \mathfrak{m}^4+(f,J(f))$. Now we define $(\beta_{ij})\in M_2(R)$ by
$$\beta_{11}=-\frac{15}{2}x^2 - \frac{7}{15}a^2xy^3 + ay^4, $$
$$\beta_{12}=\beta_{21}=-\frac{9}{2}xy - \frac{12}{25}a^2y^4, $$
$$\beta_{22}=\frac{9}{25}axy + \frac{12}{125}a^3y^4 -\frac{27}{10}y^2, $$
then the lifting $(\widetilde{\beta_{ij}})\in M_2(\mathcal{O}_2)$ satisfies$$
\begin{pmatrix}
\widetilde{\beta_{11}} & \widetilde{\beta_{12}} \\
\widetilde{\beta_{21}} & \widetilde{\beta_{22}}
\end{pmatrix}
\begin{pmatrix}
f_{x}\\
f_{y}
\end{pmatrix}
\equiv
\begin{pmatrix}
0\\
0
\end{pmatrix} \ \ \ (  \textup{mod} \ (f,J(f)^2)  ),$$ and $\widetilde{\beta_{11}}\notin \mathfrak{m}^4+(f,J(f))$, so $\beta_{11}\notin I_1^2 \mod J(f)$. By proposition \ref{prop2.14}, the Nakai Conjecture holds.
\\

7. $W_{12}$ singularity

For a $W_{12}$ singularity $(V(f),0)$ defined by $f=x^4+y^5+ax^2y^3$ for some $a\ne 0$, $f_x=4x^3+2axy^3$, $f_y=5y^4+3ax^2y^2$. As $\frac{1}{4}xf_x+\frac{1}{5}yf_y-f=\frac{1}{10}ax^2y^3$, $R=\mathcal{O}_2/(f,f_x,f_y)=\mathbb{C}\{x,y\}/(x^4,y^5,x^2y^3,  2x^3+axy^3,5y^4+3ax^2y^2)$, it has a $\mathbb{C}$-basis: $\{1,y,y^2,y^3,y^4,x,xy,xy^2,xy^3,x^2,\\x^2y\}$, $dim_{\mathbb{C}} R=11$. 

By solving the equation $Hess(f)\cdot (\alpha_1,\alpha_2)^T=0$ in $R$, we obtain that 
\begin{align*}
\begin{pmatrix}
\alpha_1\\
\alpha_2
\end{pmatrix}&\in Span_{\mathbb{C}}\{
\begin{pmatrix} y^3 \\ 0  \end{pmatrix},
\begin{pmatrix} y^4 \\ 0  \end{pmatrix},
\begin{pmatrix} xy \\ 0 \end{pmatrix},
\begin{pmatrix} xy^2 \\ 0 \end{pmatrix},
\begin{pmatrix} xy^3 \\ 0 \end{pmatrix},
\begin{pmatrix} x^2 \\ 0 \end{pmatrix},
\begin{pmatrix} x^2y \\ 0 \end{pmatrix},
\begin{pmatrix} 0 \\ y^2  \end{pmatrix},
\begin{pmatrix} 0 \\ y^3  \end{pmatrix},
\begin{pmatrix} 0 \\ y^4  \end{pmatrix},\\&
\begin{pmatrix} 0 \\ xy  \end{pmatrix},
\begin{pmatrix} 0 \\ xy^2 \end{pmatrix},
\begin{pmatrix} 0 \\ xy^3 \end{pmatrix},
\begin{pmatrix} 0 \\ x^2 \end{pmatrix},
\begin{pmatrix} 0 \\ x^2y \end{pmatrix}
\}.
\end{align*}
We see that the ideal $I_1\subset \mathfrak{m}^2+(f,J(f))$, $I_1^2\subset \mathfrak{m}^4+(f,J(f))$. Now we define $(\beta_{ij})\in M_2(R)$ by
$$\beta_{11}=-10x^2+ay^3, $$
$$\beta_{12}=\beta_{21}=\frac{3}{125}a^4xy^3 - \frac{2}{5}a^2xy^2 - 8xy, $$
$$\beta_{22}=\frac{16}{25}ax^2 +\frac{8}{125}a^2y^3 - \frac{32}{5}y^2, $$
then the lifting $(\widetilde{\beta_{ij}})\in M_2(\mathcal{O}_2)$ satisfies$$
\begin{pmatrix}
\widetilde{\beta_{11}} & \widetilde{\beta_{12}} \\
\widetilde{\beta_{21}} & \widetilde{\beta_{22}}
\end{pmatrix}
\begin{pmatrix}
f_{x}\\
f_{y}
\end{pmatrix}
\equiv
\begin{pmatrix}
0\\
0
\end{pmatrix} \ \ \ (  \textup{mod} \ (f,J(f)^2)  ),$$ and $\widetilde{\beta_{11}}\notin \mathfrak{m}^4+(f,J(f))$, so $\beta_{11}\notin I_1^2 \mod J(f)$. By proposition \ref{prop2.14}, the Nakai Conjecture holds.
\\

8. $W_{13}$ singularity

For a $W_{13}$ singularity $(V(f),0)$ defined by $f=x^4+xy^4+ay^6$ for some $a\ne 0$, $f_x=4x^3+y^4$, $f_y=4xy^3+6ay^5$. As $\frac{1}{4}xf_x+\frac{3}{16}yf_y-f=\frac{1}{8}ay^6$, $R=\mathcal{O}_2/(f,f_x,f_y)=\mathbb{C}\{x,y\}/(x^4,xy^4,y^6,  4x^3+y^4,2xy^3+3ay^5)$, it has a $\mathbb{C}$-basis: $\{1,y,y^2,y^3,y^4,y^5,x,xy,xy^2,x^2,x^2y,x^2y^2\}$, $dim_{\mathbb{C}} R=12$.

By solving the equation $Hess(f)\cdot (\alpha_1,\alpha_2)^T=0$ in $R$, we obtain that 
\begin{align*}
\begin{pmatrix}
\alpha_1\\
\alpha_2
\end{pmatrix}&\in Span_{\mathbb{C}}\{
\begin{pmatrix} y^3 \\ 0  \end{pmatrix},
\begin{pmatrix} y^4 \\ 0  \end{pmatrix},
\begin{pmatrix} y^5 \\ 0  \end{pmatrix},
\begin{pmatrix} xy \\ \frac{3}{4}y^2 \end{pmatrix},
\begin{pmatrix} xy^2 \\ 0 \end{pmatrix},
\begin{pmatrix} x^2 \\ 0 \end{pmatrix},
\begin{pmatrix} x^2y \\ 0 \end{pmatrix},
\begin{pmatrix} x^2y^2 \\ 0 \end{pmatrix},
\begin{pmatrix} 0 \\ y^3  \end{pmatrix},
\begin{pmatrix} 0 \\ y^4  \end{pmatrix},\\&
\begin{pmatrix} 0 \\ y^5  \end{pmatrix},
\begin{pmatrix} 0 \\ xy  \end{pmatrix},
\begin{pmatrix} 0 \\ xy^2 \end{pmatrix},
\begin{pmatrix} 0 \\ x^2 \end{pmatrix},
\begin{pmatrix} 0 \\ x^2y \end{pmatrix},
\begin{pmatrix} 0 \\ x^2y^2 \end{pmatrix}
\}.
\end{align*}
We see that the ideal $I_1\subset \mathfrak{m}^2+(f,J(f))$, $I_1^2\subset \mathfrak{m}^4+(f,J(f))$. Now we define $(\beta_{ij})\in M_2(R)$ by
$$\beta_{11}=-3x^2 + axy^2, $$
$$\beta_{12}=\beta_{21}=-\frac{11}{2}a^2x^2y - \frac{9}{4}xy + \frac{9}{8}ay^3, $$
$$\beta_{22}=\frac{413}{48}a^4x^2y^2 -\frac{9}{8}ax^2 - \frac{39}{16}a^2xy^2 - \frac{27}{16}y^2, $$
then the lifting $(\widetilde{\beta_{ij}})\in M_2(\mathcal{O}_2)$ satisfies$$
\begin{pmatrix}
\widetilde{\beta_{11}} & \widetilde{\beta_{12}} \\
\widetilde{\beta_{21}} & \widetilde{\beta_{22}}
\end{pmatrix}
\begin{pmatrix}
f_{x}\\
f_{y}
\end{pmatrix}
\equiv
\begin{pmatrix}
0\\
0
\end{pmatrix} \ \ \ (  \textup{mod} \ (f,J(f)^2)  ),$$ and $\widetilde{\beta_{11}}\notin \mathfrak{m}^4+(f,J(f))$, so $\beta_{11}\notin I_1^2 \mod J(f)$. By proposition \ref{prop2.14}, the Nakai Conjecture holds.
\\

9. $Q_{10}$ singularity

For a $Q_{10}$ singularity $(V(f),0)$ defined by $f=x^3+y^4+yz^2+axy^3$ for some $a\ne 0$, $f_x=3x^2+ay^3$, $f_y=4y^3+z^2+3axy^2$, $f_z=2yz$. As $\frac{1}{3}xf_x+\frac{2}{9}yf_y+\frac{7}{18}zf_z-f=-\frac{1}{9}y^4$, $\frac{1}{4}xf_x+\frac{1}{4}yf_y+\frac{3}{8}zf_z-f=-\frac{1}{4}x^3$, 
$R=\mathcal{O}_3/(f,f_x,f_y,f_z)=\mathbb{C}\{x,y,z\}/(x^3,y^4,yz,xy^3,3x^2+ay^3,4y^3+z^2+3axy^2)$, it has a $\mathbb{C}$-basis: $\{1,z,z^2,y,y^2,y^3,x,xz,xy\}$, $dim_{\mathbb{C}} R=9$.
By solving the equation $Hess(f)\cdot (\alpha_1,\alpha_2,\alpha_3)^T=0$ in $R$, we obtain that 
\begin{align*}
\begin{pmatrix}
\alpha_1\\
\alpha_2\\
\alpha_3
\end{pmatrix}&\in Span_{\mathbb{C}}\{
\begin{pmatrix} z^2 \\ 0 \\ 0  \end{pmatrix},
\begin{pmatrix} y^3 \\ 0 \\ 0  \end{pmatrix},
\begin{pmatrix} xz \\ 0 \\ 0  \end{pmatrix},
\begin{pmatrix} xy \\ 0 \\ 0  \end{pmatrix},
\begin{pmatrix} 0 \\ z \\ 3axy+4y^2  \end{pmatrix},
\begin{pmatrix} 0 \\ z^2 \\ 0  \end{pmatrix},
\begin{pmatrix} 0 \\ y^2 \\ 0  \end{pmatrix},
\begin{pmatrix} 0 \\ y^3 \\ 0  \end{pmatrix},
\begin{pmatrix} 0 \\ xz \\ 0  \end{pmatrix},\\&
\begin{pmatrix} 0 \\ xy \\ 0  \end{pmatrix},
\begin{pmatrix} 0 \\ 0 \\ z^2  \end{pmatrix},
\begin{pmatrix} 0 \\ 0 \\ y^3  \end{pmatrix},
\begin{pmatrix} 0 \\ 0 \\ xz  \end{pmatrix}
\}.
\end{align*}
We see that the ideal $I_1\subset (z^2,y^3,xz,xy)+(f,J(f))$, $I_1^2\subset (z^2,y^3,xz,xy)^2+(f,J(f))=(f,J(f))$. Now we define $(\beta_{ij})\in M_3(R)$ by
$$\beta_{11}=31ay^3 + az^2, $$
$$\beta_{12}=\beta_{21}= -\frac{243}{5}xy - \frac{18}{5}a^2y^3, $$
$$\beta_{13}=\beta_{31}= -\frac{729}{10}xz, $$ 
$$\beta_{22}= \frac{81}{20}axy + \frac{63}{80}a^3y^3 -\frac{729}{20}y^2, $$
$$\beta_{23}=\beta_{32}=0, $$
$$\beta_{33}=-\frac{243}{20}y^3 - \frac{1701}{20}z^2, $$
then the lifting $(\widetilde{\beta_{ij}})\in M_3(\mathcal{O}_3)$ satisfies$$
\begin{pmatrix}
\widetilde{\beta_{11}} & \widetilde{\beta_{12}} & \widetilde{\beta_{13}}\\
\widetilde{\beta_{21}} & \widetilde{\beta_{22}} & \widetilde{\beta_{23}}\\
\widetilde{\beta_{31}} & \widetilde{\beta_{32}} & \widetilde{\beta_{33}}
\end{pmatrix}
\begin{pmatrix}
f_{x}\\
f_{y}\\
f_{z}
\end{pmatrix}
\equiv
\begin{pmatrix}
0\\
0\\
0
\end{pmatrix} \ \ \ (  \textup{mod} \ (f,J(f)^2)  ),$$ and $\widetilde{\beta_{11}}\notin (f,J(f))$, so $\beta_{11}\notin I_1^2 \mod J(f)$. By proposition \ref{prop2.14}, the Nakai Conjecture holds.
\\

10. $Q_{11}$ singularity

For a $Q_{11}$ singularity $(V(f),0)$ defined by $f=x^3+y^2z+xz^3+az^5$ for some $a\ne 0$, $f_x=3x^2+z^3$, $f_y=2yz$, $f_z=y^2+3xz^2+5az^4$. As $\frac{1}{3}xf_x+\frac{7}{18}yf_y+\frac{2}{9}zf_z-f=\frac{1}{9}az^5$, $\frac{2}{5}xf_x+\frac{2}{5}yf_y+\frac{1}{5}zf_z-f=\frac{1}{5}x^3$, 
$R=\mathcal{O}_3/(f,f_x,f_y,f_z)=\mathbb{C}\{x,y,z\}/(x^3,z^5,xz^3, 3x^2+z^3,yz,y^2+3xz^2+5az^4)$, it has a $\mathbb{C}$-basis: $\{1,z,z^2,z^3,z^4,y,y^2,x,xz,xy\}$, $dim_{\mathbb{C}} R=10$.

By solving the equation $Hess(f)\cdot (\alpha_1,\alpha_2,\alpha_3)^T=0$ in $R$, we obtain that 
\begin{align*}
\begin{pmatrix}
\alpha_1\\
\alpha_2\\
\alpha_3
\end{pmatrix}&\in Span_{\mathbb{C}}\{
\begin{pmatrix} z^3 \\ 0 \\ 0  \end{pmatrix},
\begin{pmatrix} z^4 \\ 0 \\ 0  \end{pmatrix},
\begin{pmatrix} y^2 \\ 0 \\ 0  \end{pmatrix},
\begin{pmatrix} xz \\ 0 \\ \frac{2}{3}z^2  \end{pmatrix},
\begin{pmatrix} xy \\ 0 \\ 0  \end{pmatrix},
\begin{pmatrix} 0 \\ \frac{3}{5}xz+az^3 \\ \frac{1}{5}y  \end{pmatrix},
\begin{pmatrix} 0 \\ z^4 \\ 0  \end{pmatrix},
\begin{pmatrix} 0 \\ y^2 \\ 0  \end{pmatrix},
\begin{pmatrix} 0 \\ xy \\ 0  \end{pmatrix},\\&
\begin{pmatrix} 0 \\ 0 \\ z^3  \end{pmatrix},
\begin{pmatrix} 0 \\ 0 \\ z^4  \end{pmatrix},
\begin{pmatrix} 0 \\ 0 \\ y^2  \end{pmatrix},
\begin{pmatrix} 0 \\ 0 \\ xz  \end{pmatrix},
\begin{pmatrix} 0 \\ 0 \\ xy  \end{pmatrix}
\}.
\end{align*}
We see that the ideal $I_1\subset (z^3,y^2,xz,xy)+(f,J(f))$, $I_1^2\subset (z^3,y^2,xz,xy)^2+(f,J(f))=(f,J(f))$. Now we define $(\beta_{ij})\in M_3(R)$ by
$$\beta_{11}= -\frac{3}{14}a y^2 + z^3, $$
$$\beta_{12}=\beta_{21}= -\frac{7}{2}xy, $$
$$\beta_{13}=\beta_{31}=-2xz +\frac{169}{126}a^2y^2 + \frac{16}{21}az^3, $$ 
 $$\beta_{22}= -\frac{49}{12}y^2 + \frac{7}{12}az^4, $$
 $$\beta_{23}=\beta_{32}= 0, $$
 $$\beta_{33}=-\frac{6}{7}a xz - \frac{1139}{378}a^3 y^2 - \frac{97}{63}a^2 z^3 - \frac{4}{3}z^2, $$
then the lifting $(\widetilde{\beta_{ij}})\in M_3(\mathcal{O}_3)$ satisfies$$
\begin{pmatrix}
\widetilde{\beta_{11}} & \widetilde{\beta_{12}} & \widetilde{\beta_{13}}\\
\widetilde{\beta_{21}} & \widetilde{\beta_{22}} & \widetilde{\beta_{23}}\\
\widetilde{\beta_{31}} & \widetilde{\beta_{32}} & \widetilde{\beta_{33}}
\end{pmatrix}
\begin{pmatrix}
f_{x}\\
f_{y}\\
f_{z}
\end{pmatrix}
\equiv
\begin{pmatrix}
0\\
0\\
0
\end{pmatrix} \ \ \ (  \textup{mod} \ (f,J(f)^2)  ),$$ and $\widetilde{\beta_{11}}\notin (f,J(f))$, so $\beta_{11}\notin I_1^2 \mod J(f)$. By proposition \ref{prop2.14}, the Nakai Conjecture holds.
\\

11. $Q_{12}$ singularity

For a $Q_{12}$ singularity $(V(f),0)$ defined by $f=x^3+y^5+yz^2+axy^4$ for some $a\ne 0$, $f_x=3x^2+ay^4$, $f_y=5y^4+z^2+4axy^3$, $f_z=2yz$. As $\frac{1}{3}xf_x+\frac{1}{6}yf_y+\frac{5}{12}zf_z-f=-\frac{1}{6}y^5$, $\frac{1}{5}xf_x+\frac{1}{5}yf_y+\frac{2}{5}zf_z-f=-\frac{2}{5}x^3$, 
$R=\mathcal{O}_3/(f,f_x,f_y,f_z)=\mathbb{C}\{x,y,z\}/(x^3,y^5,xy^4,yz, 3x^2+ay^4, 5y^4+z^2+4axy^3)$, it has a $\mathbb{C}$-basis: $\{1,z,z^2,y,y^2,y^3,y^4,x,xz,xy,xy^2\}$, $dim_{\mathbb{C}} R=11$.

By solving the equation $Hess(f)\cdot (\alpha_1,\alpha_2,\alpha_3)^T=0$ in $R$, we obtain that 
\begin{align*}
\begin{pmatrix}
\alpha_1\\
\alpha_2\\
\alpha_3
\end{pmatrix}&\in Span_{\mathbb{C}}\{
\begin{pmatrix} z^2 \\ 0 \\ 0  \end{pmatrix},
\begin{pmatrix} y^4 \\ 0 \\ 0  \end{pmatrix},
\begin{pmatrix} xz \\ 0 \\ 0  \end{pmatrix},
\begin{pmatrix} xy \\ 0 \\ 0  \end{pmatrix},
\begin{pmatrix} xy^2 \\ 0 \\ 0  \end{pmatrix},
\begin{pmatrix} 0 \\ z \\ 4axy^2+5y^3  \end{pmatrix},
\begin{pmatrix} 0 \\ z^2 \\ 0  \end{pmatrix},
\begin{pmatrix} 0 \\ y^2 \\ 0  \end{pmatrix},
\begin{pmatrix} 0 \\ y^3 \\ 0  \end{pmatrix},\\&
\begin{pmatrix} 0 \\ y^4 \\ 0  \end{pmatrix},
\begin{pmatrix} 0 \\ xz \\ 0  \end{pmatrix},
\begin{pmatrix} 0 \\ xy \\ 0  \end{pmatrix},
\begin{pmatrix} 0 \\ xy^2 \\ 0  \end{pmatrix},
\begin{pmatrix} 0 \\ 0 \\ z^2  \end{pmatrix},
\begin{pmatrix} 0 \\ 0 \\ y^4  \end{pmatrix},
\begin{pmatrix} 0 \\ 0 \\ xz  \end{pmatrix}
\}.
\end{align*}
We see that the ideal $I_1\subset (z^2,y^4,xz,xy)+(f,J(f))$, $I_1^2\subset (z^2,y^4,xz,xy)^2+(f,J(f))=(f,J(f))$. Now we define $(\beta_{ij})\in M_3(R)$ by
$$\beta_{11}= 29ay^4 + az^2, $$
$$\beta_{12}=\beta_{21}=-\frac{216}{7}xy - \frac{132}{35}a^2y^4, $$
$$\beta_{13}=\beta_{31}=-\frac{432}{7}xz, $$
$$\beta_{22}= \frac{108}{35}axy + \frac{24}{35}a^3y^4 - \frac{648}{35}y^2, $$
$$\beta_{23}=\beta_{32}=0, $$
$$\beta_{33}= -\frac{108}{7}y^4 - \frac{540}{7}z^2, $$
then the lifting $(\widetilde{\beta_{ij}})\in M_3(\mathcal{O}_3)$ satisfies$$
\begin{pmatrix}
\widetilde{\beta_{11}} & \widetilde{\beta_{12}} & \widetilde{\beta_{13}}\\
\widetilde{\beta_{21}} & \widetilde{\beta_{22}} & \widetilde{\beta_{23}}\\
\widetilde{\beta_{31}} & \widetilde{\beta_{32}} & \widetilde{\beta_{33}}
\end{pmatrix}
\begin{pmatrix}
f_{x}\\
f_{y}\\
f_{z}
\end{pmatrix}
\equiv
\begin{pmatrix}
0\\
0\\
0
\end{pmatrix} \ \ \ (  \textup{mod} \ (f,J(f)^2)  ),$$ and $\widetilde{\beta_{11}}\notin (f,J(f))$, so $\beta_{11}\notin I_1^2 \mod J(f)$. By proposition \ref{prop2.14}, the Nakai Conjecture holds.
\\

12. $S_{11}$ singularity

For a $S_{11}$ singularity $(V(f),0)$ defined by $f=x^4+y^2z+xz^2+ax^3z$ for some $a\ne 0$, $f_x=4x^3+z^2+3ax^2z$, $f_y=2yz$, $f_z=y^2+2xz+ax^3$. As $\frac{1}{4}xf_x+\frac{3}{8}yf_y+\frac{1}{4}zf_z-f=-\frac{1}{4}xz^2$, $\frac{1}{5}xf_x+\frac{3}{10}yf_y+\frac{2}{5}zf_z-f=-\frac{1}{5}x^4$, 
$R=\mathcal{O}_3/(f,f_x,f_y,f_z)=\mathbb{C}\{x,y,z\}/(x^4,xz^2,yz,4x^3+z^2+3ax^2z,y^2+2xz+ax^3)$, it has a $\mathbb{C}$-basis: $\{1,z,z^2,y,y^2,x,xz,xy,x^2,x^2y\}$, $dim_{\mathbb{C}} R=10$.

By solving the equation $Hess(f)\cdot (\alpha_1,\alpha_2,\alpha_3)^T=0$ in $R$, we obtain that 
\begin{align*}
\begin{pmatrix}
\alpha_1\\
\alpha_2\\
\alpha_3
\end{pmatrix}&\in Span_{\mathbb{C}}\{
\begin{pmatrix} z^2 \\ 0 \\ 0  \end{pmatrix},
\begin{pmatrix} y^2 \\ 0 \\ 0  \end{pmatrix},
\begin{pmatrix} xz \\ 0 \\ 0  \end{pmatrix},
\begin{pmatrix} xy \\ 0 \\ 0  \end{pmatrix},
\begin{pmatrix} x^2 \\ 0 \\ -xz  \end{pmatrix},
\begin{pmatrix} x^2y \\ 0 \\ 0  \end{pmatrix},
\begin{pmatrix} 0 \\ z^2 \\ 0  \end{pmatrix},
\begin{pmatrix} 0 \\ y^2 \\ 0  \end{pmatrix},
\begin{pmatrix} 0 \\ xz \\ 0  \end{pmatrix},\\&
\begin{pmatrix} 0 \\ xy \\ 2xz  \end{pmatrix},
\begin{pmatrix} 0 \\ x^2y \\ 0  \end{pmatrix},
\begin{pmatrix} 0 \\ 0 \\ z^2  \end{pmatrix},
\begin{pmatrix} 0 \\ 0 \\ 2xz+y^2  \end{pmatrix},
\begin{pmatrix} 0 \\ 0 \\ x^2y  \end{pmatrix}
\}.
\end{align*}
We see that the ideal $I_1\subset (z^2,y^2,xz,xy,x^2)+(f,J(f))$, $I_1^2\subset (z^2,y^2,xz,xy,x^2)^2+(f,J(f))=(f,J(f))$. Now we define $(\beta_{ij})\in M_3(R)$ by
$$\beta_{11}= 2axz + ay^2, $$
$$\beta_{12}=\beta_{21}=-\frac{5}{4}a^2x^2y, $$ 
$$\beta_{13}=\beta_{31}= -4xz - 2y^2+ \frac{1}{2}az^2, $$
$$\beta_{22}=\frac{25}{3}xz + \frac{25}{6}y^2 - \frac{25}{24}az^2, $$
$$\beta_{23}=\beta_{32}=0, $$
$$\beta_{33}=0, $$
then the lifting $(\widetilde{\beta_{ij}})\in M_3(\mathcal{O}_3)$ satisfies$$
\begin{pmatrix}
\widetilde{\beta_{11}} & \widetilde{\beta_{12}} & \widetilde{\beta_{13}}\\
\widetilde{\beta_{21}} & \widetilde{\beta_{22}} & \widetilde{\beta_{23}}\\
\widetilde{\beta_{31}} & \widetilde{\beta_{32}} & \widetilde{\beta_{33}}
\end{pmatrix}
\begin{pmatrix}
f_{x}\\
f_{y}\\
f_{z}
\end{pmatrix}
\equiv
\begin{pmatrix}
0\\
0\\
0
\end{pmatrix} \ \ \ (  \textup{mod} \ (f,J(f)^2)  ),$$ and $\widetilde{\beta_{11}}\notin (f,J(f))$, so $\beta_{11}\notin I_1^2 \mod J(f)$. By proposition \ref{prop2.14}, the Nakai Conjecture holds.
\\

13. $S_{12}$ singularity

For a $S_{12}$ singularity $(V(f),0)$ defined by $f=x^2y+y^2z+xz^3+az^5$ for some $a\ne 0$, $f_x=2xy+z^3$, $f_y=x^2+2yz$, $f_z=y^2+3xz^2+5az^4$. As $\frac{2}{5}xf_x+\frac{2}{5}yf_y+\frac{1}{5}zf_z-f=\frac{1}{5}x^2y$, $\frac{4}{13}xf_x+\frac{5}{13}yf_y+\frac{3}{13}zf_z-f=\frac{2}{13}az^5$, 
$R=\mathcal{O}_3/(f,f_x,f_y,f_z)=\mathbb{C}\{x,y,z\}/(x^2y,z^5,xz^3,y^2z,2xy+z^3,x^2+2yz,y^2+3xz^2+5az^4)$. And in $R$, $xyz^2=0$, $y^3=-5ayz^4=\frac{5}{2}ax^2z^3=0$, $xy^2=-\frac{1}{2}yz^3=\frac{1}{4}x^2z^2=-\frac{1}{12}x(y^2+5az^4)=-\frac{1}{12}xy^2$, then $xy^2=0$ in $R$, $R$ has a $\mathbb{C}$-basis: $\{1,z,z^2,z^3,z^4,y,yz,yz^2,y^2,x,xz\}$, $dim_{\mathbb{C}} R=11$.

By solving the equation $Hess(f)\cdot (\alpha_1,\alpha_2,\alpha_3)^T=0$ in $R$, we obtain that 
\begin{align*}
\begin{pmatrix}
\alpha_1\\
\alpha_2\\
\alpha_3
\end{pmatrix}&\in Span_{\mathbb{C}}\{
\begin{pmatrix} z^2 \\ xz\\ -\frac{20}{3}axz+\frac{2}{3}y  \end{pmatrix},
\begin{pmatrix} z^3 \\ 0 \\ 0  \end{pmatrix},
\begin{pmatrix} z^4 \\ 0 \\ 0  \end{pmatrix},
\begin{pmatrix} yz \\ 0 \\ -xz  \end{pmatrix},
\begin{pmatrix} yz^2 \\ 0 \\ 0  \end{pmatrix},
\begin{pmatrix} y^2 \\ 0 \\ 0  \end{pmatrix},
\begin{pmatrix} xz \\ \frac{5}{4}yz \\ \frac{3}{4}z^2  \end{pmatrix},
\begin{pmatrix} 0 \\ z^3 \\ 2xz  \end{pmatrix},\\&
\begin{pmatrix} 0 \\ z^4 \\ 0  \end{pmatrix},
\begin{pmatrix} 0 \\ yz^2 \\ 0  \end{pmatrix},
\begin{pmatrix} 0 \\ y^2 \\ 0  \end{pmatrix},
\begin{pmatrix} 0 \\ 0 \\ z^3  \end{pmatrix},
\begin{pmatrix} 0 \\ 0 \\ z^4  \end{pmatrix},
\begin{pmatrix} 0 \\ 0 \\ yz  \end{pmatrix},
\begin{pmatrix} 0 \\ 0 \\ yz^2  \end{pmatrix},
\begin{pmatrix} 0 \\ 0 \\ y^2  \end{pmatrix}
\}.
\end{align*}
We see that the ideal $I_1\subset (z^2,yz,y^2,xz)+(f,J(f))$, $I_1^2\subset (z^2,yz,y^2,xz)^2+(f,J(f))=(z^4, f,J(f))$. Now we define $(\beta_{ij})\in M_3(R)$ by
$$\beta_{11}= yz^2, $$
$$\beta_{12}=\beta_{21}= \frac{5}{16}z^4, $$ 
$$\beta_{13}=\beta_{31}= \frac{1}{8}y^2, $$
$$\beta_{22}= 0, $$ 
$$\beta_{23}=\beta_{32}=-\frac{15}{32}yz^2, $$
$$\beta_{33}=\frac{9}{8}ayz^2 - \frac{9}{32}z^3, $$ 
then the lifting $(\widetilde{\beta_{ij}})\in M_3(\mathcal{O}_3)$ satisfies$$
\begin{pmatrix}
\widetilde{\beta_{11}} & \widetilde{\beta_{12}} & \widetilde{\beta_{13}}\\
\widetilde{\beta_{21}} & \widetilde{\beta_{22}} & \widetilde{\beta_{23}}\\
\widetilde{\beta_{31}} & \widetilde{\beta_{32}} & \widetilde{\beta_{33}}
\end{pmatrix}
\begin{pmatrix}
f_{x}\\
f_{y}\\
f_{z}
\end{pmatrix}
\equiv
\begin{pmatrix}
0\\
0\\
0
\end{pmatrix} \ \ \ (  \textup{mod} \ (f,J(f)^2)  ),$$ and $\widetilde{\beta_{11}}\notin (z^4, f,J(f))$, so $\beta_{11}\notin I_1^2 \mod J(f)$. By proposition \ref{prop2.14}, the Nakai Conjecture holds.
\\

14. $U_{12}$ singularity

For a $U_{12}$ singularity $(V(f),0)$ defined by $f=x^3+y^3+z^4+axyz^2$ for some $a\ne 0$, $f_x=3x^2+ayz^2$, $f_y=3y^2+axz^2$, $f_z=4z^3+2axyz$. As $\frac{1}{3}xf_x+\frac{1}{3}yf_y+\frac{1}{4}zf_z-f=\frac{1}{6}axyz^2$, 
$R=\mathcal{O}_3/(f,f_x,f_y,f_z)=\mathbb{C}\{x,y,z\}/(x^3,y^3,z^4,xyz^2,3x^2+ayz^2,3y^2+axz^2,4z^3+2axyz)$. it has a $\mathbb{C}$-basis: $\{1,z,z^2,z^3,y,yz,yz^2,y^2,x,xz,xy\}$, $dim_{\mathbb{C}} R=11$.

By solving the equation $Hess(f)\cdot (\alpha_1,\alpha_2,\alpha_3)^T=0$ in $R$, we obtain that 
\begin{align*}
\begin{pmatrix}
\alpha_1\\
\alpha_2\\
\alpha_3
\end{pmatrix}&\in Span_{\mathbb{C}}\{
\begin{pmatrix} z^3 \\ 0 \\ 0  \end{pmatrix},
\begin{pmatrix} yz^2 \\ 0 \\ 0  \end{pmatrix},
\begin{pmatrix} y^2 \\ 0 \\ 0  \end{pmatrix},
\begin{pmatrix} xz \\ 0 \\ 0  \end{pmatrix},
\begin{pmatrix} xy \\ 0 \\ 0  \end{pmatrix},
\begin{pmatrix} 0 \\ z^3 \\ 0  \end{pmatrix},
\begin{pmatrix} 0 \\ yz \\ 0  \end{pmatrix},
\begin{pmatrix} 0 \\ yz^2 \\ 0  \end{pmatrix},
\begin{pmatrix} 0 \\ y^2 \\ 0  \end{pmatrix},\\&
\begin{pmatrix} 0 \\ xy \\ 0  \end{pmatrix},
\begin{pmatrix} 0 \\ 0 \\ z^2  \end{pmatrix},
\begin{pmatrix} 0 \\ 0 \\ z^3  \end{pmatrix},
\begin{pmatrix} 0 \\ 0 \\ yz  \end{pmatrix},
\begin{pmatrix} 0 \\ 0 \\ yz^2  \end{pmatrix},
\begin{pmatrix} 0 \\ 0 \\ y^2  \end{pmatrix},
\begin{pmatrix} 0 \\ 0 \\ xz  \end{pmatrix},
\begin{pmatrix} 0 \\ 0 \\ xy  \end{pmatrix}
\}.
\end{align*}
We see that the ideal $I_1\subset (z^3,yz^2,y^2,xz,xy)+(f,J(f))$, $I_1^2\subset (z^3,yz^2,y^2,xz,xy)^2+(f,J(f))=(f,J(f))$. Now we define $(\beta_{ij})\in M_3(R)$ by
$$\beta_{11}=ayz^2, $$
$$\beta_{12}=\beta_{21}=-2xy, $$
$$\beta_{13}=\beta_{31}=-\frac{3}{2}xz, $$
$$\beta_{22}=-3y^2, $$
$$\beta_{23}=\beta_{32}=-\frac{3}{2}yz, $$
$$\beta_{33}=\frac{3}{16}axy-\frac{9}{8}z^2, $$
then the lifting $(\widetilde{\beta_{ij}})\in M_3(\mathcal{O}_3)$ satisfies$$
\begin{pmatrix}
\widetilde{\beta_{11}} & \widetilde{\beta_{12}} & \widetilde{\beta_{13}}\\
\widetilde{\beta_{21}} & \widetilde{\beta_{22}} & \widetilde{\beta_{23}}\\
\widetilde{\beta_{31}} & \widetilde{\beta_{32}} & \widetilde{\beta_{33}}
\end{pmatrix}
\begin{pmatrix}
f_{x}\\
f_{y}\\
f_{z}
\end{pmatrix}
\equiv
\begin{pmatrix}
0\\
0\\
0
\end{pmatrix} \ \ \ (  \textup{mod} \ (f,J(f)^2)  ),$$ and $\widetilde{\beta_{11}}\notin (f,J(f))$, so $\beta_{11}\notin I_1^2 \mod J(f)$. By proposition \ref{prop2.14}, the Nakai Conjecture holds.
\\

\section{The bimodal isolated hypersurfce singularity case}\label{sec 4}
The bimodal isolated hypersurface singularities were classified into 8 infinite series and 14 exceptional families in \cite{Arnold} as the following three tables.

\begin{center}
\begin{longtable}{|p{1.2cm}|p{5.5cm}|p{3cm}|p{3cm}|}
\caption{4 infinite series of bimodal singularities of corank 2}\label{table-1}\\
\hline
symbol & normal forms & Milnor number $\mu$ & Tjurina number $\tau$\\
\hline 
\multirow{2}{1cm}{$J_{3,0}$} &
$x^3+bx^2y^3+y^9+cxy^7$ & \multirow{2}{1cm}{$16$} & $16$ ($c=0$)\\
& ($4b^3+27\ne 0$) & & $15$ ($c\ne 0$)\\
\hline 
\multirow{2}{1cm}{$J_{3,p}$} & $x^3+x^2y^3+(a_0+a_1y)y^{9+p}$ & \multirow{2}{1cm}{$16+p$} & \multirow{2}{1cm}{$14+p$}\\
& ($p>0, a_0\ne 0$) & &\\
\hline 
\multirow{2}{1cm}{$Z_{1,0}$} & $y(x^3+dx^2y^2+cxy^5+y^6)$ & \multirow{2}{1cm}{$15$} & $15$ ($c=0$)\\
& ($4d^3+27\ne 0$) & & $14$ ($c\ne 0$)\\
\hline
\multirow{2}{1cm}{$Z_{1,p}$} & $x^3y+x^2y^3+(a_0+a_1y)y^{7+p}$ & \multirow{2}{1cm}{$15+p$} & \multirow{2}{1cm}{$13+p$}\\
& ($p>0, a_0\ne 0$) & &\\
\hline
\multirow{2}{1cm}{$W_{1,0}$} & $x^4+(a_0+a_1y)x^2y^3+y^6$ & \multirow{2}{1cm}{$15$} & $15$ ($a_1=0$)\\
& ($a_0^2\ne 4$) & & $14$ ($a_1\ne 0$)\\
\hline
\multirow{2}{1cm}{$W_{1,p}$} & $x^4+x^2y^3+(a_0+a_1y)y^{6+p}$ & \multirow{2}{1cm}{$15+p$} & \multirow{2}{1cm}{$13+p$}\\
& ($p>0, a_0\ne 0$) & &\\
\hline
\multirow{2}{1cm}{$W_{1,2q-1}^{\#}$} & $(x^2+y^3)^2+(a_0+a_1y)xy^{4+q}$ & \multirow{2}{3cm}{$15+2q-1$} & \multirow{2}{3cm}{$13+2q-1$}\\
& ($q>0, a_0\ne 0$) & & \\
\hline
\multirow{2}{1cm}{$W_{1,2q}^{\#}$} & $(x^2+y^3)^2+(a_0+a_1y)x^2y^{3+q}$ & 
\multirow{2}{3cm}{$15+2q$} & \multirow{2}{3cm}{$13+2q$}\\
& ($q>0, a_0\ne 0$) & & \\
\hline
\end{longtable}
\end{center}

\begin{center}
\begin{longtable}{|p{1.2cm}|p{5.5cm}|p{3cm}|p{3cm}|}
\caption{4 infinite series of bimodal singularities of corank 3}\label{table-2}\\
\hline
symbol & normal forms & Milnor number $\mu$ & Tjurina number $\tau$\\
\hline 
\multirow{2}{1cm}{$Q_{2,0}$} & $x^3+yz^2+(a_0+a_1y)x^2y^2+xy^4$ & \multirow{2}{1cm}{$14$} & $14$ ($a_1=0$)\\
& ($a_0^2\ne 4$) & & $13$ ($a_1\ne 0$)\\
\hline
\multirow{2}{1cm}{$Q_{2,p}$} & $x^3+yz^2+x^2y^2+(a_0+a_1y)y^{6+p}$ & \multirow{2}{1cm}{$14+p$} & \multirow{2}{1cm}{$12+p$}\\
& ($p>0, a_0\ne 0$) & & \\
\hline
\multirow{2}{1cm}{$S_{1,0}$} & $x^2z+yz^2+y^5+(a_0+a_1y)zy^3$ & \multirow{2}{1cm}{$14$} & $14$ ($a_1=0$)\\
& ($a_0^2\ne 4$) & & $13$ ($a_1\ne 0$)\\
\hline
\multirow{2}{1cm}{$S_{1,p}$} & $x^2z+yz^2+x^2y^2+(a_0+a_1y)y^{5+p}$ & 
\multirow{2}{1cm}{$14+p$} & \multirow{2}{1cm}{$12+p$}\\
& ($p>0, a_0\ne 0$) & & \\
\hline
\multirow{2}{1cm}{$S_{1,2q-1}^{\#}$} & $x^2z+yz^2+zy^3+(a_0+a_1y)xy^{3+q}$ & \multirow{2}{3cm}{$14+2q-1$} & \multirow{2}{3cm}{$12+2q-1$}\\
& ($q>0, a_0\ne 0$) & &\\
\hline
\multirow{2}{1cm}{$S_{1,2q}^{\#}$} & $x^2z+yz^2+zy^3+(a_0+a_1y)x^2y^{2+q}$ & \multirow{2}{3cm}{$14+2q$} & \multirow{2}{3cm}{$12+2q$}\\
& ($q>0, a_0\ne 0$) & &\\
\hline
\multirow{2}{1cm}{$U_{1,0}$} & $x^3+xz^2+xy^3+(a_0+a_1y)y^3z$ & \multirow{2}{1cm}{$14$} & $14$ ($a_1=0$)\\
& ($a_0(a_0^2+1)\ne 0$) & & $13$ ($a_1\ne 0$)\\
\hline
\multirow{2}{1cm}{$U_{1,2q-1}$} & $x^3+xz^2+xy^3+(a_0+a_1y)y^{1+q}z^2$ &
\multirow{2}{3cm}{$14+2q-1$} & \multirow{2}{3cm}{$12+2q-1$}\\
& ($q>0, a_0\ne 0$) & &\\
\hline
\multirow{2}{1cm}{$U_{1,2q}$} & $x^3+xz^2+xy^3+(a_0+a_1y)y^{3+q}z$ &
\multirow{2}{3cm}{$14+2q$} & \multirow{2}{3cm}{$12+2q$}\\
& ($q>0, a_0\ne 0$) & &\\
\hline
\end{longtable}
\end{center}

\begin{center}
\begin{longtable}{|p{1.2cm}|p{5.5cm}|p{3cm}|p{3cm}|}
\caption{14 exceptional families of bimodal singularities (the three numbers in the column of Tjurina number corresponds to the cases of 
 $a_0=a_1=0$, $a_0=0  \ \& \ a_1\ne 0$, and $a_0\ne 0$)}\label{table-3}\\
\hline
symbol & normal forms & Milnor number $\mu$ & Tjurina number $\tau$\\
\hline 
$E_{18}$ & $x^3+y^{10}+(a_0+a_1y)xy^7$ & $18$ & $18,17,16$\\
\hline
$E_{19}$ & $x^3+xy^7+(a_0+a_1y)y^{11}$ & $19$ & $19,18,17$\\
\hline
$E_{20}$ & $x^3+y^{11}+(a_0+a_1y)xy^8$ & $20$ & $20,19,18$\\
\hline
$Z_{17}$ & $x^3y+y^8+(a_0+a_1y)xy^6$ & $17$ & $17,16,15$\\
\hline
$Z_{18}$ & $x^3y+xy^6+(a_0+a_1y)y^9$ & $18$ & $18,17,16$\\
\hline
$Z_{19}$ & $x^3y+y^9+(a_0+a_1y)xy^7$ & $19$ & $19,18,17$\\
\hline
$W_{17}$ & $x^4+xy^5+(a_0+a_1y)y^7$ & $17$ & $17,16,15$\\
\hline
$W_{18}$ & $x^4+y^7+(a_0+a_1y)x^2y^4$ & $18$ & $18,17,16$\\
\hline
$Q_{16}$ & $x^3+yz^2+y^7+(a_0+a_1y)xy^5$ & $16$ & $16,15,14$\\
\hline
$Q_{17}$ & $x^3+yz^2+xy^5+(a_0+a_1y)y^8$ & $17$ & $17,16,15$\\
\hline
$Q_{18}$ & $x^3+yz^2+y^8+(a_0+a_1y)xy^6$ & $18$ & $18,17,16$\\
\hline
$S_{16}$ & $x^2z+yz^2+xy^4+(a_0+a_1y)y^6$ & $16$ & $16,15,14$\\
\hline
$S_{17}$ & $x^2z+yz^2+y^6+(a_0+a_1y)zy^4$ & $17$ & $17,16,15$\\
\hline
$U_{16}$ & $x^3+xz^2+y^5+(a_0+a_1y)x^2y^2$ & $16$ & $16,15,14$\\
\hline
\end{longtable}
\end{center}
In this section, for an isolated hypersurface singularity $(V(f),0)\subset (\mathbb{C}^n,0)$, $A=\mathcal{O}_n/(f)$, with notations as in section \ref{sec 2}, we will verify the Nakai Conjecture for these classified families of singularities in this section.

\subsection{Cases of corank 2 bimodal singularities}
\ 
\newline
\indent
1. $J_{3,0}$ singularity
    
    For a $J_{3,0}$ singularity $(V(f),0)$ defined by $f=x^3+bx^2 y^3+y^9+cx y^7$ for some $b,c$ with $4 b^3+27\ne 0$, 
$f_x=3x^2+2bx y^3+cy^7$, $f_y=3bx^2y^2+9y^8+7cx y^6$. If $c=0$, then it reduces to the weighted homogeneous binomial case, and $(V(f)\cap \{x=0\},0)$ is an isolated singularity in $\mathbb{A}^1_{y}$, then the Nakai Conjecture holds by proposition \ref{prop3.1}, now we only need to consider the case of $c\ne 0$. 

As $4xf_x+yf_y-11f=x^3-2y^9$, $3xf_x+yf_y-9f=cxy^7$, 
then $R=\mathbb{C}\{x,y\}/(f,f_x,f_y)\simeq \mathbb{C}\{x,y\}/(x^3-2y^9,xy^7,3x^2+2bx y^3+cy^7,3bx^2y^2+9y^8+7cx y^6)$. In $R$, $y^{10}=-\frac{y^2}{9}(3bx^2y^2+7cxy^6)=-\frac{b}{3}x^2y^4=\frac{b}{9}y^4(2bxy^3+cy^7)=\frac{bc}{9}y^{11}$, as $1-\frac{bc}{9}y$ is invertible in $R$, $y^{10}=0$ in $R$. Then $R$ has a $\mathbb{C}$-basis: $\{1,y,y^2,\cdots,y^7,x,xy,xy^2,\cdots,xy^6\}$, $dim_{\mathbb{C}} R=15$.

(1) If $b\ne 0$, by solving the equation $Hess(f)\cdot (\alpha_1,\alpha_2)^{T}=0$ in $R$, we obtain that 
\begin{align*}
\begin{pmatrix}
\alpha_1\\
\alpha_2
\end{pmatrix}&\in Span_{\mathbb{C}}\{
\begin{pmatrix}
    y^5\\ (-\frac{13c}{36b}+\frac{21c}{8b^4})x y +(\frac{b}{9}+\frac{3}{2b^2})x 
\end{pmatrix}, 
\begin{pmatrix}
    y^6\\
    (\frac{b}{9}+\frac{3}{2b^2})xy
\end{pmatrix},
\begin{pmatrix}
    y^7\\
    0
\end{pmatrix},
\begin{pmatrix}
    x y\\
    -\frac{7c^2}{216b^2}xy+\frac{c}{54}x+\frac{1}{3}y^2
\end{pmatrix},
\\&
\begin{pmatrix}
    xy^2\\
    -(\frac{c}{18}+\frac{7c}{8b^3})x y-\frac{1}{2b}x
\end{pmatrix},
\begin{pmatrix}
    xy^3\\
    -\frac{1}{2b}xy
\end{pmatrix},
\begin{pmatrix}
    xy^4\\
    0
\end{pmatrix},
\begin{pmatrix}
    xy^5\\
    0
\end{pmatrix},
\begin{pmatrix}
    xy^6\\
    0
\end{pmatrix},
\begin{pmatrix}
    0\\
    (\frac{2c}{9}+\frac{21c}{8b^3})xy+\frac{3}{2b}x+y^3
\end{pmatrix},\\&  
\begin{pmatrix}
    0\\
    \frac{3}{2b}xy+y^4
\end{pmatrix},
\begin{pmatrix}
    0\\
    y^5
\end{pmatrix},
\begin{pmatrix}
    0\\
    y^6
\end{pmatrix},
\begin{pmatrix}
    0\\
    y^7
\end{pmatrix},
\begin{pmatrix}
    0\\
    xy^2
\end{pmatrix},
\begin{pmatrix}
    0\\
    xy^3
\end{pmatrix},
\begin{pmatrix}
    0\\
    xy^4
\end{pmatrix},
\begin{pmatrix}
    0\\
    xy^5
\end{pmatrix},
\begin{pmatrix}
    0\\
    xy^6
\end{pmatrix}
\}.
\end{align*}
We see that the ideal $I_1\subset (y^5,xy)+(f,J(f))$, $I_1^2\subset (y^{10},xy^6,x^2y^2)+(f,J(f))=(xy^5,f,J(f))$. Now we define 
$(\beta_{ij})\in M_2(R)$ by $$
\beta_{11}=\frac{111}{2}c^2xy^6-81b^2cxy^5+(18b^4+\frac{243}{2}b)xy^4,$$
$$
\beta_{12}=\beta_{21}=\frac{2347}{96}c^3xy^5-\frac{97}{12}b^2c^2xy^4+(b^4c+\frac{99}{2}bc)xy^3-(9b^3+\frac{243}{4})xy^2+(\frac{1}{2}b^3c^2+\frac{149}{4}c^2)y^7+\frac{3}{2}b^2cy^6,$$
$$
\beta_{22}=\frac{79}{96}bc^2xy^2-(\frac{1}{2}b^3c+\frac{237}{8}c)xy+\frac{725}{96}c^3y^6-\frac{103}{48}b^2c^2y^5-4bcy^4-(3b^3+\frac{81}{4})y^3, $$then the lifting $(\widetilde{\beta_{ij}})\in M_2(\mathcal{O}_2)$ satisfies$$
\begin{pmatrix}
\widetilde{\beta_{11}} & \widetilde{\beta_{12}} \\
\widetilde{\beta_{21}} & \widetilde{\beta_{22}}
\end{pmatrix}
\begin{pmatrix}
f_{x}\\
f_{y}
\end{pmatrix}
\equiv
\begin{pmatrix}
0\\
0
\end{pmatrix} \ \ \ (  \textup{mod} \ (f,J(f)^2)  ),$$ and from $18b^4+\frac{243}{2}b=18b(b^3+\frac{27}{4})\ne 0$, $\widetilde{\beta_{11}}\notin (x y^5,f,J(f))$, so $\beta_{11}\notin I_1^2 \mod J(f)$. By proposition \ref{prop2.14}, the Nakai Conjecture holds.

(2) If $b=0$, similarly by solving the equation $Hess(f)\cdot (\alpha_1,\alpha_2)^{T}=0$ in $R$, we obtain that 
\begin{align*}
    \begin{pmatrix}
\alpha_1\\
\alpha_2
\end{pmatrix}&\in Span_{\mathbb{C}}\{
\begin{pmatrix}
    y^6\\
    -\frac{6}{7c}x
\end{pmatrix},
\begin{pmatrix}
    y^7\\
    0
\end{pmatrix},
\begin{pmatrix}
    xy\\
    -\frac{2c}{9}x
\end{pmatrix},
\begin{pmatrix}
    xy^2\\
    0
\end{pmatrix},
\begin{pmatrix}
    xy^3\\ 0
\end{pmatrix},
\begin{pmatrix}
    xy^4\\ 0
\end{pmatrix},
\begin{pmatrix}
    xy^5\\ 0
\end{pmatrix},
\begin{pmatrix}
    xy^6\\ 0
\end{pmatrix},
\begin{pmatrix}
    0\\ \frac{7c}{9}x+y^2
\end{pmatrix},\\&
\begin{pmatrix}
    0\\ y^3
\end{pmatrix},
\begin{pmatrix}
    0\\ y^4
\end{pmatrix},
\begin{pmatrix}
    0\\ y^5
\end{pmatrix},
\begin{pmatrix}
    0\\ y^6
\end{pmatrix},
\begin{pmatrix}
    0\\ y^7
\end{pmatrix},
\begin{pmatrix}
    0\\xy
\end{pmatrix},
\begin{pmatrix}
    0\\xy^2
\end{pmatrix},
\begin{pmatrix}
    0\\xy^3
\end{pmatrix},
\begin{pmatrix}
    0\\xy^4
\end{pmatrix},
\begin{pmatrix}
    0\\xy^5
\end{pmatrix},
\begin{pmatrix}
    0\\xy^6
\end{pmatrix}
\}.
\end{align*}
We see that the ideal $I_1\subset (y^6,xy)+(f,J(f))$, $I_1^2\subset (f,J(f))$. Now we define $(\beta_{ij})\in M_2(R)$ by $$
\beta_{11}=cxy^6,$$
$$\beta_{12}=\beta_{21}=\frac{1}{3}cy^7,$$
$$\beta_{22}= -\frac{1}{4}xy-\frac{1}{108}c^2y^6,$$then the lifting $(\widetilde{\beta_{ij}})\in M_2(\mathcal{O}_2)$ satisfies$$
\begin{pmatrix}
\widetilde{\beta_{11}} & \widetilde{\beta_{12}} \\
\widetilde{\beta_{21}} & \widetilde{\beta_{22}}
\end{pmatrix}
\begin{pmatrix}
f_{x}\\
f_{y}
\end{pmatrix}
\equiv
\begin{pmatrix}
0\\
0
\end{pmatrix} \ \ \ (  \textup{mod} \ (f,J(f)^2)  ),$$ and $\widetilde{\beta_{11}}\notin (f,J(f))$, so $\beta_{11}\notin I_1^2 \mod J(f)$. By proposition \ref{prop2.14}, the Nakai Conjecture holds.
\\

2.$J_{3,p}$ singularity

For a $J_{3,p}$ singularity $(V(f),0)$ defined by $f=x^3+x^2y^3+(a_0+a_1y)y^{9+p}$ for some $p>0,a_0\ne 0$, 
$f_x=3x^2+2x y^3$, $f_y=3x^2y^2+(9+p)a_0y^{8+p}+(10+p)a_1y^{9+p}$.
As $\frac{1}{3}xf_x+\frac{1}{9}yf_y-f=y^{9+p}(\frac{p}{9}(a_0+a_1y)+\frac{a_1}{9}y)$ and $\frac{p}{9}(a_0+a_1y)+\frac{a_1}{9}y$ is invertible in $R=\mathbb{C}\{x,y\}/(f,f_x,f_y)$, $y^{9+p}=0$ in $R$. So $R=\mathbb{C}\{x,y\}/(x^3+x^2y^3, 3x^2+2xy^3,3x^2y^2+(9+p)a_0y^{8+p},y^{9+p})=\mathbb{C}\{x,y\}/(x^3,x^2y^3,3x^2+2xy^3,3x^2y^2+(9+p)a_0y^{8+p},y^{9+p})$, it has a $\mathbb{C}$-basis: $\{1,y,y^2,\cdots,y^{8+p},x,xy,xy^2,xy^3,xy^4\}$, $dim_{\mathbb{C}}R=14+p$.

By solving the equation $Hess(f)\cdot (\alpha_1,\alpha_2)^{T}=0$ in $R$, we obtain that 
\begin{align*}
\begin{pmatrix}
\alpha_1\\
\alpha_2
\end{pmatrix}&\in Span_{\mathbb{C}}\{
\begin{pmatrix} y^{p+6}\\ 0 \end{pmatrix},
\begin{pmatrix} y^{p+7}\\ 0 \end{pmatrix},
\begin{pmatrix} y^{p+8}\\ 0 \end{pmatrix},
\begin{pmatrix} y^{p+5}\\ \frac{1}{(9+p)a_0}x \end{pmatrix},
\begin{pmatrix} xy\\ \frac{1}{3}y^2 \end{pmatrix},
\begin{pmatrix} xy^2\\ -\frac{1}{2}x \end{pmatrix},
\begin{pmatrix} xy^3\\ 0 \end{pmatrix},
\begin{pmatrix} xy^4\\ 0 \end{pmatrix},\\
&\begin{pmatrix} 0\\ \frac{3}{2}x+y^3 \end{pmatrix},
\begin{pmatrix} 0\\ y^4 \end{pmatrix},
\begin{pmatrix} 0\\ y^5 \end{pmatrix}, \cdots
\begin{pmatrix} 0\\ y^{8+p} \end{pmatrix},
\begin{pmatrix} 0\\ xy \end{pmatrix},
\begin{pmatrix} 0\\ xy^2 \end{pmatrix},
\begin{pmatrix} 0\\ xy^3 \end{pmatrix},
\begin{pmatrix} 0\\ xy^4 \end{pmatrix}
\}.
\end{align*}
We see that the ideal $I_1\subset (xy,y^{p+6})+(f,J(f))$, $I_1^2\subset (x^2y^2, xy^{p+7},y^{2p+12})+(f,J(f))=(y^{8+p},f,J(f))$.

Now for $p=1$, we define $(\beta_{ij})\in M_2(R)$ by
$$
\beta_{11}=\frac{4}{a_0}xy^4+y^8,$$
$$\beta_{12}=\beta_{21}=(-\frac{1801}{2646}+\frac{4a_1}{21a_0^2})xy^3-\frac{2}{a_0}xy^2-\frac{10}{21}y^6,
$$
\begin{align*}
    \beta_{22}=&(-\frac{494180a_1^2}{27783a_0^2}+\frac{320a_1^3}{49a_0^4})xy^3+(-\frac{70019a_1}{294a_0}+\frac{2168a_1^2}{21a_0^3})xy^2+(\frac{35999}{294}-\frac{706a_1}{21a_0^2})xy+\frac{2}{21a_0}x\\&+(-\frac{68129a_1}{441a_0}+\frac{464a_1^2}{7a_0^3})y^5+(\frac{34109}{441}-\frac{148a_1}{7a_0^2})y^4-\frac{4}{7a_0}y^3;
\end{align*}
for $p=2$ we define $(\beta_{ij})\in M_2(R)$ by
$$
\beta_{11}=\frac{2}{a_0}xy^4+y^9,
$$
$$
\beta_{12}=\beta_{21}=\frac{a_1}{8a_0^2}xy^3-\frac{1}{a_0}xy^2-\frac{11}{24}y^7,$$
$$
\beta_{22}=(\frac{1111a_1}{768a_0}+\frac{187a_1^3}{256a_0^4})xy^3+(\frac{1089}{32}+\frac{411a_1^2}{32a_0^3})xy^2-\frac{49a_1}{8a_0^2}xy+\frac{1}{12a_0}x+(\frac{165}{8}+\frac{63a_1^2}{8a_0^3})y^5-\frac{29a_1}{8a_0^2}y^4-\frac{1}{4a_0}y^3;
$$
for $p=3$ we define $(\beta_{ij})\in M_2(R)$ by
$$
\beta_{11}=\frac{4}{3a_0}xy^4+y^{10},
$$
$$
\beta_{12}=\beta_{21}=\frac{8a_1}{81a_0^2}xy^3-\frac{2}{3a_0}xy^2-\frac{4}{9}y^8,
$$
$$
\beta_{22}=(-\frac{2}{27}+\frac{1280a_1^3}{6561a_0^4})xy^3+\frac{1040a_1^2}{243a_0^3}xy^2-\frac{194a_1}{81a_0^2}xy+\frac{2}{27a_0}x+\frac{608a_1^2}{243a_0^3}y^5-\frac{4a_1}{3a_0^2}y^4-\frac{4}{27a_0}y^3;
$$
and for $p\ge 4$, we define $(\beta_{ij})\in M_2(R)$ by
$$
\beta_{11}=\frac{4}{pa_0}xy^4+y^{p+7},
$$
$$
\beta_{12}=\beta_{21}=\frac{2(p+1)a_1}{3p(p+6)a_0^2}xy^3-\frac{2}{pa_0}xy^2-\frac{p+9}{3(p+6)}y^{p+5},
$$
\begin{align*}
    \beta_{22}=&-\frac{2(p+1)^2(p+9)(p^2-p-36)a_1^3}{9p^4(p+6)^2a_0^4}xy^3+\frac{2(p+1)(p^3+28p^2+189p+324)a_1^2}{3p^4(p+6)a_0^3}xy^2\\&-\frac{2(p^3+19p^2+171p+162)a_1}{3p^4(p+6)a_0^2}xy+\frac{2}{3(p+6)a_0}x+\frac{8(p+1)(p^2+10p+18)a_1^2}{p^4(p+6)a_1^3}y^5\\&
    -\frac{4(p^2+18p+18)a_1}{p^3(p+6)a_0^2}y^4-\frac{4}{p(p+6)a_0}y^3.
\end{align*}
Then the lifting $(\widetilde{\beta_{ij}})\in M_2(\mathcal{O}_2)$ satisfies$$
\begin{pmatrix}
\widetilde{\beta_{11}} & \widetilde{\beta_{12}} \\
\widetilde{\beta_{21}} & \widetilde{\beta_{22}}
\end{pmatrix}
\begin{pmatrix}
f_{x}\\
f_{y}
\end{pmatrix}
\equiv
\begin{pmatrix}
0\\
0
\end{pmatrix} \ \ \ (  \textup{mod} \ (f,J(f)^2)  ),$$ and $\widetilde{\beta_{11}}\notin (y^{8+p},f,J(f))$, so $\beta_{11}\notin I_1^2 \mod J(f)$.
By proposition \ref{prop2.14}, the Nakai Conjecture holds.
\\

3.$Z_{1,0}$ singularity

 For a $Z_{1,0}$ singularity $(V(f),0)$ defined by $f=y(x^3+dx^2y^2+cxy^5+y^6)$ for some $c,d$ with $4 d^3+27\ne 0$, 
$f_x=3x^2y+2dxy^3+cy^6$, $f_y=x^3+3dx^2y^2+6cxy^5+7y^6$. If $c=0$, then it reduces to the weighted homogeneous binomial case, and $(V(f)\cap \{x=0\},0)$ is an isolated singularity in $\mathbb{A}^1_{y}$, then the Nakai Conjecture holds by proposition \ref{prop3.1}, now we only need to consider the case of $c\ne 0$. 

As $2xf_x+yf_y-7f=cxy^6$, $3xf_x+yf_y-9f=x^3y-2y^7$, then $R=\mathbb{C}\{x,y\}/(f,f_x,f_y)=\mathbb{C}\{x,y\}/(xy^6,x^3y-2y^7,3x^2y+2dxy^3+cy^6,x^3+3dx^2y^2+6cxy^5+7y^6)$. In $R$, $y^8=-\frac{1}{7}y^2(x^3+3dx^2y^2+6cxy^5)=-\frac{1}{7}x^3y^2-\frac{3d}{7}x^2y^4=-\frac{2}{7}y^8-\frac{3d}{7}y^3(-\frac{1}{3})(2dxy^3+cy^6)=-\frac{2}{7}y^8+\frac{cd}{7}y^9$, i.e. $y^8(\frac{9}{7}-\frac{cd}{7}y)=0$ in $R$, as $\frac{9}{7}-\frac{cd}{7}y$ is invertible in $R$, $y^8=0$ in $R$. Then $R$ has a $\mathbb{C}$-basis: $\{1,y,y^2,\cdots,y^6,x,xy,xy^2,\cdots,xy^5,x^2\}$, $dim_{\mathbb{C}} R=14$.

(1) If $d\ne 0$, by solving the equation $Hess(f)\cdot (\alpha_1,\alpha_2)^{T}=0$ in $R$, we obtain that 
\begin{align*}
\begin{pmatrix}
\alpha_1\\
\alpha_2
\end{pmatrix}&\in Span_{\mathbb{C}}\{
\begin{pmatrix} y^3 \\ -(\frac{cd^2}{84}+\frac{59c}{56d})xy+\frac{2}{21}dx-(\frac{d^2}{21}+\frac{3}{2d})y^2 \end{pmatrix},
\begin{pmatrix}
    y^4 \\ (\frac{d}{6}+\frac{9}{4d^2})xy
\end{pmatrix},
\begin{pmatrix}   y^5 \\ 0 \end{pmatrix},
\begin{pmatrix}   y^6 \\ 0 \end{pmatrix},\\&
\begin{pmatrix}   xy \\ \frac{c}{24}xy+\frac{1}{2}y^2 \end{pmatrix},
\begin{pmatrix}   xy^2 \\ -\frac{3}{4d}xy \end{pmatrix},
\begin{pmatrix}   xy^3 \\ 0 \end{pmatrix},
\begin{pmatrix}   xy^4 \\ 0 \end{pmatrix},
\begin{pmatrix}   xy^5 \\ 0 \end{pmatrix},
\begin{pmatrix}   x^2 \\ \frac{1}{2}xy \end{pmatrix},
\begin{pmatrix}   0 \\ \frac{3}{2d}xy+y^3 \end{pmatrix},
\begin{pmatrix}   0 \\ y^4 \end{pmatrix},\\&
\begin{pmatrix}   0 \\ y^5 \end{pmatrix},
\begin{pmatrix}   0 \\ y^6 \end{pmatrix},
\begin{pmatrix}   0 \\ xy^2 \end{pmatrix},
\begin{pmatrix}   0 \\ xy^3 \end{pmatrix},
\begin{pmatrix}   0 \\ xy^4 \end{pmatrix},
\begin{pmatrix}   0 \\ xy^5 \end{pmatrix},
\begin{pmatrix}   0 \\ x^2 \end{pmatrix}
\}.
\end{align*}
We see that the ideal $I_1\subset (y^3,xy,x^2)+(f,J(f))$, $I_1^2\subset (y^3,xy,x^2)^2+(f,J(f))=(y^6,xy^4,f,J(f))$.
Now we define $(\beta_{ij})\in M_2(R)$ by $$
\beta_{11}=-(\frac{8d^3}{9c}+\frac{6}{c})x^2+y^5,
$$
$$
\beta_{12}=\beta_{21}=\frac{\frac{5}{288}c^2d^3+\frac{27}{64}c^2}{d^3+\frac{27}{4}}xy^4-\frac{\frac{7}{108}cd^6+\frac{81}{32}c}{d(d^3+\frac{27}{4})}xy^3-\frac{4d^3+27}{9c}xy+\frac{\frac{1}{8}cd^3-\frac{45}{16}c}{d^3+\frac{27}{4}}y^5+\frac{d^2}{9}y^4,
$$
\begin{align*}
    \beta_{22}=&\frac{(\frac{109}{3456}d^9+\frac{7}{16}d^6+\frac{9423}{2048}d^3+\frac{28431}{4096})c^3}{(d^3+\frac{27}{4})^2}xy^4+\frac{(-\frac{7}{1296}d^6+\frac{59}{384}d^3-\frac{315}{256})c^2}{d(d^3+\frac{27}{4})}xy^3+\frac{(\frac{1}{16}d^3+\frac{9}{32})c}{d^2}xy^2\\&+\frac{1}{4}xy+\frac{d}{12}y^3-\frac{2(d^3+\frac{27}{4})}{9c}y^2,
\end{align*}
then the lifting $(\widetilde{\beta_{ij}})\in M_2(\mathcal{O}_2)$ satisfies$$
\begin{pmatrix}
\widetilde{\beta_{11}} & \widetilde{\beta_{12}} \\
\widetilde{\beta_{21}} & \widetilde{\beta_{22}}
\end{pmatrix}
\begin{pmatrix}
f_{x}\\
f_{y}
\end{pmatrix}
\equiv
\begin{pmatrix}
0\\
0
\end{pmatrix} \ \ \ (  \textup{mod} \ (f,J(f)^2)  ),$$ and $\widetilde{\beta_{11}}\notin (y^6,xy^4,f,J(f))$, so $\beta_{11}\notin I_1^2 \mod J(f)$.
By proposition \ref{prop2.14}, the Nakai Conjecture holds.

(2) If $d=0$,  similarly by solving the equation $Hess(f)\cdot (\alpha_1,\alpha_2)^{T}=0$ in $R$, we obtain that 
\begin{align*}
\begin{pmatrix}
\alpha_1\\
\alpha_2
\end{pmatrix}&\in Span_{\mathbb{C}}\{
\begin{pmatrix}   y^5 \\ 0 \end{pmatrix},
\begin{pmatrix}   y^6 \\ 0 \end{pmatrix},
\begin{pmatrix}   xy \\ 0 \end{pmatrix},
\begin{pmatrix}   xy^2 \\ 0 \end{pmatrix},
\begin{pmatrix}   xy^3 \\ 0 \end{pmatrix},
\begin{pmatrix}   xy^4 \\ 0 \end{pmatrix},
\begin{pmatrix}   xy^5 \\ 0 \end{pmatrix},
\begin{pmatrix}   x^2 \\ 0 \end{pmatrix},
\begin{pmatrix}   0 \\ y^2 \end{pmatrix},
\begin{pmatrix}   0 \\ y^3 \end{pmatrix},\\&
\begin{pmatrix}   0 \\ y^4 \end{pmatrix},
\begin{pmatrix}   0 \\ y^5 \end{pmatrix},
\begin{pmatrix}   0 \\ y^6 \end{pmatrix},
\begin{pmatrix}   0 \\ xy \end{pmatrix},
\begin{pmatrix}   0 \\ xy^2 \end{pmatrix},
\begin{pmatrix}   0 \\ xy^3 \end{pmatrix},
\begin{pmatrix}   0 \\ xy^4 \end{pmatrix},
\begin{pmatrix}   0 \\ xy^5 \end{pmatrix},
\begin{pmatrix}   0 \\ x^2 \end{pmatrix}
\}.
\end{align*}
We see that the ideal $I_1\subset (y^5,xy,x^2)+(f,J(f))$, $I_1^2\subset (y^5,xy,x^2)^2+(f,J(f))=(f,J(f))$.
Now we define $(\beta_{ij})\in M_2(R)$ by $$
\beta_{11}=y^6,
$$
$$
\beta_{12}=\beta_{21}=-\frac{15}{13c}xy^2,
$$
$$
\beta_{22}=-\frac{5}{338}xy^2-\frac{15}{26c}y^3,$$
then the lifting $(\widetilde{\beta_{ij}})\in M_2(\mathcal{O}_2)$ satisfies$$
\begin{pmatrix}
\widetilde{\beta_{11}} & \widetilde{\beta_{12}} \\
\widetilde{\beta_{21}} & \widetilde{\beta_{22}}
\end{pmatrix}
\begin{pmatrix}
f_{x}\\
f_{y}
\end{pmatrix}
\equiv
\begin{pmatrix}
0\\
0
\end{pmatrix} \ \ \ (  \textup{mod} \ (f,J(f)^2)  ),$$ and $\widetilde{\beta_{11}}\notin (f,J(f))$, so $\beta_{11}\notin I_1^2 \mod J(f)$.
By proposition \ref{prop2.14}, the Nakai Conjecture holds.
\\

4.$Z_{1,p}$ singularity

For a $Z_{1,p}$ singularity $(V(f),0)$ defined by $f=x^3y+x^2y^3+(a_0+a_1y)y^{7+p}$ for some $p>0,a_0\ne 0$, 
$f_x=3x^2y+2x y^3$, $f_y=x^3+3x^2y^2+(7+p)a_0y^{6+p}+(8+p)a_1y^{7+p}$.
As $2xf_x+yf_y-7f=(pa_0+(p+1)a_1y)y^{7+p}$ and $pa_0+(p+1)a_1y$ is invertible in $R=\mathbb{C}\{x,y\}/(f,f_x,f_y)$, $y^{7+p}=0$ in $R$. So $R=\mathbb{C}\{x,y\}/(y^{7+p},3x^2y+2x y^3, x^3+3x^2y^2+(7+p)a_0y^{6+p})=\mathbb{C}\{x,y\}/(x^3y,x^2y^3,y^{7+p},3x^2y+2x y^3, x^3+3x^2y^2+(7+p)a_0y^{6+p})=\mathbb{C}\{x,y\}/(xy^5,y^{7+p},3x^2y+2x y^3, x^3+3x^2y^2+(7+p)a_0y^{6+p})$, it has a $\mathbb{C}$-basis: $\{1,y,y^2,\cdots,y^{6+p},x,xy,xy^2,xy^3,xy^4,x^2\}$, $dim_{\mathbb{C}}R=13+p$.

By solving the equation $Hess(f)\cdot (\alpha_1,\alpha_2)^{T}=0$ in $R$, we obtain that 
\begin{align*}
\begin{pmatrix}
\alpha_1\\
\alpha_2
\end{pmatrix}&\in Span_{\mathbb{C}}\{
\begin{pmatrix} y^{p+3} \\ \frac{2}{3(p+7)a_0}x-\frac{1}{3(p+7)a_0}y^2 \end{pmatrix},
\begin{pmatrix} y^{p+4} \\ 0 \end{pmatrix},
\begin{pmatrix} y^{p+5} \\ 0 \end{pmatrix},
\begin{pmatrix} y^{p+6} \\ 0 \end{pmatrix},
\begin{pmatrix} xy \\ \frac{1}{2}y^2 \end{pmatrix},
\begin{pmatrix} x y^2 \\ 0 \end{pmatrix},
\begin{pmatrix} x y^3 \\ 0 \end{pmatrix},\\&
\begin{pmatrix} x y^4 \\ 0 \end{pmatrix},
\begin{pmatrix} x^2 \\ 0 \end{pmatrix},
\begin{pmatrix} 0 \\ y^3 \end{pmatrix},
\begin{pmatrix} 0 \\ y^4 \end{pmatrix},
\begin{pmatrix} 0 \\ y^5 \end{pmatrix},
\cdots \begin{pmatrix} 0 \\ y^{p+6} \end{pmatrix},
\begin{pmatrix} 0 \\ xy \end{pmatrix},
\begin{pmatrix} 0 \\ xy^2 \end{pmatrix},
\begin{pmatrix} 0 \\ xy^3 \end{pmatrix},
\begin{pmatrix} 0 \\ xy^4 \end{pmatrix},
\begin{pmatrix} 0 \\ x^2 \end{pmatrix}
\}
\end{align*}
We see that the ideal $I_1\subset (x^2,xy,y^{p+3})+(f,J(f))$, $I_1^2\subset (x^2,xy,y^{p+3})^2+(f,J(f))=(xy^4,f,J(f))$.

Now for $p=1$, we define $(\beta_{ij})\in M_2(R)$ by
$$\beta_{11}=\frac{8}{3a_0}xy^3+y^6,$$
$$\beta_{12}=\beta_{21}=(-\frac{589}{600}+\frac{4a_1}{15a_0^2})xy^3-\frac{2}{a_0}xy^2-\frac{7}{10}y^5,
$$
$$\beta_{22}=(\frac{49427a_1}{6000a_0}-\frac{28a_1^2}{15a_0^3})xy^3+(\frac{24693}{200}-\frac{164a_1}{5a_0^2})xy^2+\frac{1}{5a_0}xy+(\frac{7601}{100}-\frac{20a_1}{a_0^2})y^4-\frac{4}{5a_0}y^3;
$$
for $p=2$ we define $(\beta_{ij})\in M_2(R)$ by
$$\beta_{11}=\frac{4}{3a_0}xy^3+y^7,$$
$$\beta_{12}=\beta_{21}=\frac{a_1}{6a_0^2}xy^3-\frac{1}{a_0}xy^2-\frac{2}{3}y^6,
$$
$$\beta_{22}=-(\frac{1}{6}+\frac{a_1^2}{6a_0^3})xy^3-\frac{6a_1}{a_0^2}xy^2+\frac{1}{6a_0}xy-\frac{10a_1}{3a_0^2}y^4-\frac{1}{3a_0}y^3;
$$
and for $p\ge 3$ we define $(\beta_{ij})\in M_2(R)$ by
$$\beta_{11}=\frac{8}{3pa_0}xy^3+y^{p+5},$$
$$\beta_{12}=\beta_{21}=\frac{2(p+1)a_1}{3p(p+4)a_0^2}xy^3-\frac{2}{pa_0}xy^2-\frac{p+6}{2(p+4)}y^{p+4},
$$
\begin{align*}
    \beta_{22}=&\frac{(p+1)(p+6)(p-3)a_1^2}{3p^3(p+4)a_0^3}xy^3-\frac{(p^3+13p^2+78p+72)a_1}{p^3(p+4)a_0^2}xy^2+\frac{1}{(p+4)a_0}xy\\&-\frac{4(p^2+12p+12)a_1}{p^3(p+4)a_0^2}y^4-\frac{4}{p(p+4)a_0}y^3.
\end{align*}
Then the lifting $(\widetilde{\beta_{ij}})\in M_2(\mathcal{O}_2)$ satisfies$$
\begin{pmatrix}
\widetilde{\beta_{11}} & \widetilde{\beta_{12}} \\
\widetilde{\beta_{21}} & \widetilde{\beta_{22}}
\end{pmatrix}
\begin{pmatrix}
f_{x}\\
f_{y}
\end{pmatrix}
\equiv
\begin{pmatrix}
0\\
0
\end{pmatrix} \ \ \ (  \textup{mod} \ (f,J(f)^2)  ),$$ and $\widetilde{\beta_{11}}\notin (xy^4,f,J(f))$, so $\beta_{11}\notin I_1^2 \mod J(f)$.
By proposition \ref{prop2.14}, the Nakai Conjecture holds.
\\

5.$W_{1,0}$ singularity

 For a $W_{1,0}$ singularity $(V(f),0)$ defined by $f=x^4+(a_0+a_1y)x^2y^3+y^6$ for some $a_0,a_1$ with $a_0^2\ne 4$, 
$f_x=4x^3+2xy^3(a_0+a_1y)$, $f_y=3a_0x^2y^2+4a_1x^2y^3+6y^5$. If $a_1=0$, then it reduces to the weighted homogeneous binomial case, and $(V(f)\cap \{x=0\},0)$ is an isolated singularity in $\mathbb{A}^1_{y}$, then the Nakai Conjecture holds by proposition \ref{prop3.1}, now we only need to consider the case of $a_1\ne 0$. 

As $3xf_x+2yf_y-12f=2a_1x^2y^4$, $xf_x-2f=2(x^4-y^6)$, $8f-2xf_x-yf_y=a_0x^2y^3+2y^6$, then $R=\mathbb{C}\{x,y\}/(f,f_x,f_y)=\mathbb{C}\{x,y\}/(x^2y^4,x^4-y^6,a_0x^2y^3+2y^6, 4x^3+2xy^3(a_0+a_1y),3a_0x^2y^2+4a_1x^2y^3+6y^5)$. In $R$, $y^7=-\frac{a_0}{2}x^2y^4=0$, $xy^5=-\frac{x}{6}(3a_0x^2y^2+4a_1x^2y^3)=-\frac{a_0}{2}x^3y^2-\frac{2a_1}{3}x^3y^3=-(\frac{a_0}{2}y^2+\frac{2a_1}{3}y^3)\cdot(-\frac{1}{2}xy^3)(a_0+a_1y)=xy^5\cdot(a_0+a_1y)(\frac{a_0}{4}+\frac{a_1}{3}y)$, as $a_0^2\ne 4$, then $1-(a_0+a_1y)(\frac{a_0}{4}+\frac{a_1}{3}y)=(1-\frac{a_0^2}{4})-(\frac{7}{12}a_0a_1+\frac{a_1^2}{3}y)\cdot y$ is invertible in $R$, $xy^5=0$ in $R$. So $R=\mathbb{C}\{x,y\}/(y^7,x y^5,x^2y^4, x^4-y^6, 4x^3+2xy^3(a_0+a_1y),3a_0x^2y^2+4a_1x^2y^3+6y^5)$, it has a $\mathbb{C}$-basis: $\{1,y, y^2,y^3,y^4,y^5,x,xy,xy^2,xy^3,xy^4,x^2,x^2y,x^2y^2\}$, $dim_{\mathbb{C}} R=14$.

(1)If $a_0\ne 0, a_0^2\ne 30$, by solving the equation $Hess(f)\cdot (\alpha_1,\alpha_2)^T=0$ in $R$, we obtain that
\begin{align*}
    \begin{pmatrix}
\alpha_1\\
\alpha_2
\end{pmatrix}&\in Span_{\mathbb{C}}\{
\begin{pmatrix}
    y^3 \\ (\frac{a_0}{6}-\frac{2}{a_0})xy
\end{pmatrix},
\begin{pmatrix}  y^4 \\ 0 \end{pmatrix},
\begin{pmatrix}  y^5 \\ 0 \end{pmatrix},
\begin{pmatrix}  xy \\ \frac{2}{3}y^2 \end{pmatrix},
\begin{pmatrix}  xy^2 \\ 0 \end{pmatrix},
\begin{pmatrix}  xy^3 \\ 0 \end{pmatrix},
\begin{pmatrix}  xy^4 \\ 0 \end{pmatrix},
\begin{pmatrix}  x^2 \\ \frac{2}{3}xy \end{pmatrix},
\begin{pmatrix}  x^2y \\ 0 \end{pmatrix},\\&
\begin{pmatrix}  x^2y^2 \\ 0 \end{pmatrix},
\begin{pmatrix}  0 \\ y^3 \end{pmatrix},
\begin{pmatrix}  0 \\ y^4 \end{pmatrix},
\begin{pmatrix}  0 \\ y^5 \end{pmatrix},
\begin{pmatrix}  0 \\ xy^2 \end{pmatrix},
\begin{pmatrix}  0 \\ xy^3 \end{pmatrix},
\begin{pmatrix}  0 \\ xy^4 \end{pmatrix},
\begin{pmatrix}  0 \\ x^2 \end{pmatrix},
\begin{pmatrix}  0 \\ x^2y \end{pmatrix},
\begin{pmatrix}  0 \\ x^2y^2 \end{pmatrix}
\}:=W.
\end{align*}

(2)If $a_0^2=30$, similarly by solving the equation $Hess(f)\cdot (\alpha_1,\alpha_2)^T=0$ in $R$, we obtain that
\begin{align*}
    \begin{pmatrix}
\alpha_1\\
\alpha_2
\end{pmatrix}&\in W+\mathbb{C}<\begin{pmatrix}
    y^2 \\ \frac{23}{90}a_1xy+\frac{1}{10}a_0x
\end{pmatrix}>.
\end{align*}

(3)If $a_0=0$, similarly by solving the equation $Hess(f)\cdot (\alpha_1,\alpha_2)^T=0$ in $R$,
 we obtain that
\begin{align*}
    \begin{pmatrix}
\alpha_1\\
\alpha_2
\end{pmatrix}&\in Span_{\mathbb{C}}\{
\begin{pmatrix}  y^4 \\ 0 \end{pmatrix},
\begin{pmatrix}  y^5 \\ 0 \end{pmatrix},
\begin{pmatrix}  xy \\ 0 \end{pmatrix},
\begin{pmatrix}  xy^2 \\ 0 \end{pmatrix},
\begin{pmatrix}  xy^3 \\ 0 \end{pmatrix},
\begin{pmatrix}  xy^4 \\ 0 \end{pmatrix},
\begin{pmatrix}  x^2 \\ 0 \end{pmatrix},
\begin{pmatrix}  x^2y \\ 0 \end{pmatrix},
\begin{pmatrix}  x^2y^2 \\ 0 \end{pmatrix},\\&
\begin{pmatrix}  0 \\ y^2 \end{pmatrix},
\begin{pmatrix}  0 \\ y^3 \end{pmatrix},
\begin{pmatrix}  0 \\ y^4 \end{pmatrix},
\begin{pmatrix}  0 \\ y^5 \end{pmatrix},
\begin{pmatrix}  0 \\ xy \end{pmatrix},
\begin{pmatrix}  0 \\ xy^2 \end{pmatrix},
\begin{pmatrix}  0 \\ xy^3 \end{pmatrix},
\begin{pmatrix}  0 \\ xy^4 \end{pmatrix},
\begin{pmatrix}  0 \\ x^2 \end{pmatrix},
\begin{pmatrix}  0 \\ x^2y \end{pmatrix},
\begin{pmatrix}  0 \\ x^2y^2 \end{pmatrix}
\}.
\end{align*}

We see that in all cases, $I_1\subset (x^2,xy,y^2)+(f,J(f))$, $I_1^2\subset (x^2,xy,y^2)^2+(f,J(f))=(y^4,xy^3,x^2y^2,f,J(f))$.
Now we define $(\beta_{ij})\in M_2(R)$ by $$
\beta_{11}=\frac{7}{3}a_0a_1x^2y+(a_0^2-4)x^2+\frac{2}{3}a_1y^4,
$$
$$
\beta_{12}=\beta_{21}=-\frac{2}{9}a_1^2xy^3+\frac{4}{3}a_0a_1xy^2+\frac{2(a_0^2-4)}{3}xy,
$$
$$
\beta_{22}=\frac{8}{27}a_1x^2+\frac{4}{81}a_1^2y^4+\frac{28}{27}a_0a_1y^3+\frac{4}{9}(a_0^2-4)y^2,$$
then the lifting $(\widetilde{\beta_{ij}})\in M_2(\mathcal{O}_2)$ satisfies$$
\begin{pmatrix}
\widetilde{\beta_{11}} & \widetilde{\beta_{12}} \\
\widetilde{\beta_{21}} & \widetilde{\beta_{22}}
\end{pmatrix}
\begin{pmatrix}
f_{x}\\
f_{y}
\end{pmatrix}
\equiv
\begin{pmatrix}
0\\
0
\end{pmatrix} \ \ \ (  \textup{mod} \ (f,J(f)^2)  ),$$ and $\widetilde{\beta_{11}}\notin (y^4,xy^3,x^2y^2,f,J(f))$, so $\beta_{11}\notin I_1^2 \mod J(f)$.
By proposition \ref{prop2.14}, the Nakai Conjecture holds.
\\

6.$W_{1,p}$ singularity

For a $W_{1,p}$ singularity $(V(f),0)$ defined by $f=x^4+x^2y^3+(a_0+a_1y)y^{6+p}$ for some $p>0,a_0\ne 0$, 
$f_x=4x^3+2x y^3$, $f_y=3x^2y^2+(6+p)a_0y^{5+p}+(7+p)a_1y^{6+p}$. 
As $\frac{1}{4}xf_x+\frac{1}{6}yf_y-f=(\frac{p}{6}a_0+\frac{p+1}{6}a_1y)y^{p+6}$, $\frac{p}{6}a_0+\frac{p+1}{6}a_1y$ is invertible in $R=\mathbb{C}\{x,y\}/(f,f_x,f_y)$, $y^{p+6}=0$ in $R$. And in $R$, $xy^5=-2x^3y^2=\frac{2}{3}(6+p)a_0xy^{p+5}$, as $1-\frac{2}{3}(6+p)y^p$ is invertible in $R$, $xy^5=0$ in $R$. So $R=\mathbb{C}\{x,y\}/(y^{p+6},xy^5,x^4,x^2y^3,4x^3+2x y^3,3x^2y^2+(6+p)a_0y^{5+p})$, it has a $\mathbb{C}$-basis: $\{1,y,y^2,\cdots,y^{p+5},x,xy,xy^2,xy^3,xy^4,x^2,x^2y\}$, $dim_{\mathbb{C}} R=p+13$.

By solving the equation $Hess(f)\cdot (\alpha_1,\alpha_2)^{T}=0$ in $R$, we obtain that 
\begin{align*}
\begin{pmatrix}
\alpha_1\\
\alpha_2
\end{pmatrix}&\in Span_{\mathbb{C}}\{
\begin{pmatrix} y^{p+3} \\ 0  \end{pmatrix},
\begin{pmatrix} y^{p+4} \\ 0  \end{pmatrix},
\begin{pmatrix} y^{p+5} \\ 0  \end{pmatrix},
\begin{pmatrix} xy \\ \frac{2}{3}y^2  \end{pmatrix},
\begin{pmatrix} xy^2 \\ 0  \end{pmatrix},
\begin{pmatrix} xy^3 \\ 0  \end{pmatrix},
\begin{pmatrix} xy^4 \\ 0  \end{pmatrix},
\begin{pmatrix} x^2 \\ 0  \end{pmatrix},
\begin{pmatrix} x^2y \\ 0  \end{pmatrix},\\&
\begin{pmatrix} 0 \\ y^3  \end{pmatrix},
\begin{pmatrix} 0 \\ y^4  \end{pmatrix},\cdots,
\begin{pmatrix} 0 \\ y^{p+5}  \end{pmatrix},
\begin{pmatrix} 0 \\ xy  \end{pmatrix},
\begin{pmatrix} 0 \\ xy^2  \end{pmatrix},
\begin{pmatrix} 0 \\ xy^3  \end{pmatrix},
\begin{pmatrix} 0 \\ xy^4  \end{pmatrix},
\begin{pmatrix} 0 \\ x^2  \end{pmatrix},
\begin{pmatrix} 0 \\ x^2y  \end{pmatrix}
\}.
\end{align*}
We see that $I_1\subset (x^2,xy,y^{p+3})+(f,J(f))$, $I_1^2\subset (x^2,xy,y^{p+3})^2+(f,J(f))=(xy^4,y^{p+5},f,J(f))$.

Now for $p=1$, we define $(\beta_{ij})\in M_2(R)$ by
$$\beta_{11}=-\frac{3}{a_0}x^2y+y^5,$$
$$\beta_{12}=\beta_{21}=(-\frac{1}{6}+\frac{a_1}{2a_0^2})xy^3-\frac{2}{a_0}xy^2,$$
$$\beta_{22}=(\frac{2}{9}-\frac{a_1}{a_0^2})x^2y+\frac{1}{3a_0}x^2-\frac{1}{a_0}y^3;$$
and for $p\ge 2$, we define $(\beta_{ij})\in M_2(R)$ by
$$\beta_{11}=-\frac{3}{pa_0}x^2y+y^{p+4},$$
$$\beta_{12}=\beta_{21}=\frac{(p+1)a_1}{p(p+3)a_0^2}xy^3-\frac{2}{pa_0}xy^2,$$
$$\beta_{22}=-\frac{4(p+2)a_1}{3p(p+3)a_0^2}x^2y+\frac{4}{3(p+3)a_0}x^2-\frac{4}{p(p+3)a_0}y^3.$$
Then the lifting $(\widetilde{\beta_{ij}})\in M_2(\mathcal{O}_2)$ satisfies$$
\begin{pmatrix}
\widetilde{\beta_{11}} & \widetilde{\beta_{12}} \\
\widetilde{\beta_{21}} & \widetilde{\beta_{22}}
\end{pmatrix}
\begin{pmatrix}
f_{x}\\
f_{y}
\end{pmatrix}
\equiv
\begin{pmatrix}
0\\
0
\end{pmatrix} \ \ \ (  \textup{mod} \ (f,J(f)^2)  ),$$ and $\widetilde{\beta_{11}}\notin (xy^4,y^{p+5},f,J(f))$, so $\beta_{11}\notin I_1^2 \mod J(f)$.
By proposition \ref{prop2.14}, the Nakai Conjecture holds.
\\

7.$W^\#_{1,2q-1}$ singularity

For a $W^\#_{1,2q-1}$ singularity $(V(f),0)$ defined by $f=(x^2+y^3)^2+(a_0+a_1y)xy^{4+q}$ for some $q>0$, $a_0\ne 0$, $f_x=4x(x^2+y^3)+(a_0+a_1y)y^{4+q}$, $f_y=6y^2(x^2+y^3)+(4+q)a_0xy^{3+q}+(5+q)a_1xy^{4+q}$.
As $\frac{1}{4}xf_x+\frac{1}{6}yf_2-f=((\frac{q}{6}-\frac{1}{12})a_0+(\frac{q}{6}+\frac{1}{12})a_1y)xy^{4+q}$, $(\frac{q}{6}-\frac{1}{12})a_0+(\frac{q}{6}+\frac{1}{12})a_1y$ is invertible in $R=\mathbb{C}\{x,y\}/(f,f_x,f_y)$, $xy^{4+q}=0$ in $R$. And we have $x^4= -x^2y^3= y^6$ in $R$, then in $R$, $x^2y^{q+3}=-y^{q+6}$, $0=3 y^2 f_1-2 x f_2=(3(a_0+a_1 y)+2(4+q)a_0)y^{q+6}$, as $3(a_0+a_1 y)+2(4+q)a_0$ is invertible in $R$, $y^{q+6}=0$ in $R$. So $R=\mathbb{C}\{x,y\}/(x^4+x^2y^3,x^2y^3+y^6,x y^{4+q},y^{q+6},4x(x^2+y^3)+(a_0+a_1y)y^{4+q},6y^2(x^2+y^3)+(4+q)a_0xy^{3+q})$, it has a $\mathbb{C}$-basis: $\{1,y,y^2,\cdots,y^{q+5},x,xy,xy^2,\cdots,xy^{q+3},x^2,x^2y\}$, $dim_{\mathbb{C}} R=2q+12$.

By solving the equation $Hess(f)\cdot (\alpha_1,\alpha_2)^{T}=0$ in $R$, we obtain that 
\begin{align*}
\begin{pmatrix}
\alpha_1\\
\alpha_2
\end{pmatrix}&\in Span_{\mathbb{C}}\{
\begin{pmatrix} y^3\\ -\frac{2}{3}xy \end{pmatrix},
\begin{pmatrix} y^4\\ -\frac{2}{3}xy^2 \end{pmatrix},\cdots,
\begin{pmatrix} y^{q+2}\\ -\frac{2}{3}xy^q \end{pmatrix},
\begin{pmatrix} y^{q+3}\\ 0 \end{pmatrix},
\begin{pmatrix} y^{q+4}\\ 0 \end{pmatrix},
\begin{pmatrix} y^{q+5}\\ 0 \end{pmatrix},\\&
\begin{pmatrix} xy\\ \frac{2q-1}{36}a_0xy^q+\frac{2}{3}y^2 \end{pmatrix},
\begin{pmatrix} xy^2\\ -\frac{2}{3}x^2 \end{pmatrix},
\begin{pmatrix} xy^3\\ -\frac{2}{3}x^2y \end{pmatrix},
\begin{pmatrix} xy^4\\ \frac{2}{3}y^5 \end{pmatrix},
\begin{pmatrix} xy^5\\ \frac{2}{3}y^6 \end{pmatrix},\cdots,
\begin{pmatrix} xy^q\\ \frac{2}{3}y^{q+1} \end{pmatrix},\\&
\begin{pmatrix} xy^{q+1}\\ 0 \end{pmatrix},
\begin{pmatrix} xy^{q+2}\\ 0 \end{pmatrix},
\begin{pmatrix} xy^{q+3}\\ 0 \end{pmatrix},
\begin{pmatrix} x^2\\ \frac{2}{3}xy \end{pmatrix},
\begin{pmatrix} x^2y\\ \frac{2}{3}xy^2 \end{pmatrix},
\begin{pmatrix} 0\\ x^2+y^3 \end{pmatrix},
\begin{pmatrix} 0\\ x^2y+y^4 \end{pmatrix},\\&
\begin{pmatrix} 0\\ y^{q+2} \end{pmatrix},
\begin{pmatrix} 0\\ y^{q+3} \end{pmatrix},
\begin{pmatrix} 0\\ y^{q+4} \end{pmatrix},
\begin{pmatrix} 0\\ y^{q+5} \end{pmatrix},
\begin{pmatrix} 0\\ xy^{q+1} \end{pmatrix},
\begin{pmatrix} 0\\ xy^{q+2} \end{pmatrix},
\begin{pmatrix} 0\\ xy^{q+3} \end{pmatrix}
\}.
\end{align*}
We see that $I_1\subset (y^3,xy,x^2)+(f,J(f))$, $I_1^2\subset (y^3,xy,x^2)^2+(f,J(f))\subset (y^5,xy^4,f,J(f))$.

Now for $q=1$, we define $(\beta_{ij})\in M_2(R)$ by
$$\beta_{11}=13x^2y+y^4, $$
$$\beta_{12}=\beta_{21}=-\frac{61}{216}a_0x^2y+\frac{53a_0^2}{324}xy^3+8xy^2+(\frac{34981}{373248}a_0^3+\frac{11}{18}a_1)y^5-\frac{97}{216}a_0y^4, $$
$$\beta_{22}=-\frac{715a_0^2}{3888}x^2y+\frac{4}{9}x^2-(\frac{20665}{559872}a_0^3+\frac{35}{27}a_1)xy^3-\frac{7}{27}a_0xy^2-(\frac{1403173}{80621568}a_0^4+\frac{1271}{1944}a_0a_1)y^5+\frac{52}{9}y^3;$$
for $q=2$, we define $(\beta_{ij})\in M_2(R)$ by
$$\beta_{11}=5x^2y+y^4 ,$$
$$\beta_{12}=\beta_{21}=\frac{8}{3}xy^2+\frac{a_1}{6}y^6-\frac{a_0}{6}y^5,$$
$$\beta_{22}=\frac{4}{9}x^2-\frac{7a_1}{9}xy^4-\frac{a_0}{3}xy^3-\frac{29}{144}a_0^2y^6+\frac{20}{9}y^3;$$
and for $q\ge 3$, we define $(\beta_{ij})\in M_2(R)$ by
$$\beta_{11}=\frac{2q+11}{2q-1}x^2y+y^4,$$
$$\beta_{12}=\beta_{21}=\frac{8}{2q-1}xy^2-\frac{(2q-13)a_1}{18(2q-1)}y^{q+4}-\frac{a_0}{6}y^{q+3},$$
$$\beta_{22}=\frac{4}{9}x^2-\frac{(2q+3)(2q+5)a_1}{27(2q-1)}xy^{q+2}-\frac{2q+5}{27}a_0xy^{q+1}+\frac{4(2q+11)}{9(2q-1)}y^3.$$
Then the lifting $(\widetilde{\beta_{ij}})\in M_2(\mathcal{O}_2)$ satisfies$$
\begin{pmatrix}
\widetilde{\beta_{11}} & \widetilde{\beta_{12}} \\
\widetilde{\beta_{21}} & \widetilde{\beta_{22}}
\end{pmatrix}
\begin{pmatrix}
f_{x}\\
f_{y}
\end{pmatrix}
\equiv
\begin{pmatrix}
0\\
0
\end{pmatrix} \ \ \ (  \textup{mod} \ (f,J(f)^2)  ),$$ and $\widetilde{\beta_{11}}\notin (y^5,xy^4,f,J(f))$, so $\beta_{11}\notin I_1^2 \mod J(f)$.
By proposition \ref{prop2.14}, the Nakai Conjecture holds.
\\

8.$W^\#_{1,2q}$ singularity

For a $W^\#_{1,2q}$ singularity $(V(f),0)$ defined by $f=(x^2+y^3)^2+(a_0+a_1y)x^2y^{3+q}$ for some $q>0$, $a_0\ne 0$, $f_x=4x(x^2+y^3)+2(a_0+a_1y)xy^{3+q}$, $f_y=6y^2(x^2+y^3)+(3+q)a_0x^2y^{2+q}+(4+q)a_1x^2y^{3+q}$. As $\frac{1}{4}xf_x+\frac{1}{6}yf_y-f=(\frac{qa_0}{6}+\frac{(q+1)a_1}{6}y)x^2y^{3+q}$, $\frac{qa_0}{6}+\frac{(q+1)a_1}{6}y$ is invertible in $R=\mathbb{C}\{x,y\}/(f,f_x,f_y)$, $x^2y^{3+q}=0$ in $R$. And we have $x^4=-x^2y^3=y^6$ in $R$, then $0=y^{q+1}f_y=6y^{q+3}(x^2+y^3)=6y^{q+6}$ in $R$, $y^{q+6}=0$ in $R$. Since $0=3y^2f_x-2xf_y=4 a_0x y^{q+5}-2(3+q)a_0 x^3 y^{q+2}=(4a_0+2(3+q)a_0)xy^{q+5}$ in $R$, $xy^{q+5}=0$ in $R$.
Therefore, $R=\mathbb{C}\{x,y\}/(y^{q+6},xy^{q+5},x^4+x^2y^3,x^2y^3-y^6,4x(x^2+y^3)+2(a_0+a_1y)xy^{3+q},6y^2(x^2+y^3)+(3+q)a_0x^2y^{2+q})$, and $x^2y^2=-y^5-\frac{3+q}{6}a_0x^2y^{q+2}=-y^5+\frac{3+q}{6}a_0y^{q+5}$ in $R$, $R$ has a $\mathbb{C}$-basis: $\{1,y,y^2,\cdots,y^{q+5},x,xy,xy^2,\cdots,xy^{q+4},x^2,x^2y\}$, $dim_{\mathbb{C}}R=13+2q$.

By solving the equation $Hess(f)\cdot (\alpha_1,\alpha_2)^{T}=0$ in $R$, we obtain that 
\begin{align*}
\begin{pmatrix}
\alpha_1\\
\alpha_2
\end{pmatrix}&\in Span_{\mathbb{C}}\{
\begin{pmatrix}  y^3 \\ -\frac{2}{3}xy \end{pmatrix},
\begin{pmatrix}  y^4 \\ -\frac{2}{3}xy^2 \end{pmatrix},
\begin{pmatrix}  y^5 \\ -\frac{2}{3}xy^3 \end{pmatrix},\cdots,
\begin{pmatrix}  y^{q+2} \\ -\frac{2}{3}xy^{q} \end{pmatrix},
\begin{pmatrix}  y^{q+3} \\ 0 \end{pmatrix},
\begin{pmatrix}  y^{q+4} \\ 0 \end{pmatrix},
\begin{pmatrix}  y^{q+5} \\ 0 \end{pmatrix},\\&
\begin{pmatrix}  xy \\ -\frac{q}{9}a_0y^{q+2}+\frac{2}{3}y^2\end{pmatrix},
\begin{pmatrix}  xy^2 \\ -\frac{2}{3}x^2 \end{pmatrix},
\begin{pmatrix}  xy^3 \\ -\frac{2}{3}x^2y \end{pmatrix},
\begin{pmatrix}  xy^4 \\ \frac{2}{3}xy^5 \end{pmatrix},
\begin{pmatrix}  xy^5 \\ \frac{2}{3}xy^6 \end{pmatrix},\cdots,
\begin{pmatrix}  xy^{q+1} \\ \frac{2}{3}y^{q+2} \end{pmatrix},
\begin{pmatrix}  xy^{q+2} \\ 0 \end{pmatrix},\\&
\begin{pmatrix}  xy^{q+3} \\ 0 \end{pmatrix},
\begin{pmatrix}  xy^{q+4} \\ 0 \end{pmatrix},
\begin{pmatrix}  x^2 \\ \frac{2}{3}xy \end{pmatrix},
\begin{pmatrix}  x^2y \\ \frac{2}{3}xy^2 \end{pmatrix},
\begin{pmatrix}  0 \\ x^2+y^3 \end{pmatrix},
\begin{pmatrix}  0 \\ x^2y+y^4 \end{pmatrix},
\begin{pmatrix}  0 \\ y^{q+3} \end{pmatrix},
\begin{pmatrix}  0 \\ y^{q+4} \end{pmatrix},\\&
\begin{pmatrix}  0 \\ y^{q+5} \end{pmatrix},
\begin{pmatrix}  0 \\ xy^{q+1} \end{pmatrix},
\begin{pmatrix}  0 \\ xy^{q+2} \end{pmatrix},
\begin{pmatrix}  0 \\ xy^{q+3} \end{pmatrix},
\begin{pmatrix}  0 \\ xy^{q+4} \end{pmatrix}
\}.
\end{align*}
We see that $I_1\subset (y^3,xy,x^2)+(f,J(f))$, $I_1^2\subset (y^3,xy,x^2)^2+(f,J(f))\subset (y^5,xy^4,f,J(f))$.

Now for $q=1$, we define $(\beta_{ij})\in M_2(R)$ by
$$\beta_{11}=7 x^2 y+y^4,$$
$$\beta_{12}=\beta_{21}=-(\frac{a_0^2}{27}+\frac{2 a_1}{9})x y^4-\frac{a_0}{3}x y^3+4xy^2,$$
$$\beta_{22}=-\frac{2a_0}{27}x^2y+\frac{4}{9}x^2+(-\frac{2a_0^2}{27}+\frac{2a_1}{9})y^5+\frac{28}{9}y^3; $$
and for $q\ge 2$, we define $(\beta_{ij})\in M_2(R)$ by
$$\beta_{11}=\frac{q+6}{q}x^2y+y^4,$$
$$\beta_{12}=\beta_{21}=-\frac{2}{9}a_1xy^{q+3}-\frac{a_0}{3}xy^{q+2}+\frac{4}{q}xy^2,$$
$$\beta_{22}=\frac{4}{9}x^2+\frac{2(q+2)a_1}{27}y^{q+4}+\frac{2qa_0}{27}y^{q+3}+\frac{4(q+6)}{9q}y^3.$$
Then the lifting $(\widetilde{\beta_{ij}})\in M_2(\mathcal{O}_2)$ satisfies$$
\begin{pmatrix}
\widetilde{\beta_{11}} & \widetilde{\beta_{12}} \\
\widetilde{\beta_{21}} & \widetilde{\beta_{22}}
\end{pmatrix}
\begin{pmatrix}
f_{x}\\
f_{y}
\end{pmatrix}
\equiv
\begin{pmatrix}
0\\
0
\end{pmatrix} \ \ \ (  \textup{mod} \ (f,J(f)^2)  ),$$ and $\widetilde{\beta_{11}}\notin (y^5,xy^4,f,J(f))$, so $\beta_{11}\notin I_1^2 \mod J(f)$.
By proposition \ref{prop2.14}, the Nakai Conjecture holds.
\\

\subsection{Cases of corank 3 bimodal singularities}
\ 
\newline
\indent
1. $Q_{2,0}$ singularity

For a $Q_{2,0}$ singularity $(V(f),0)$ defined by $f=x^3+yz^2+(a_0+a_1y)x^2y^2+xy^4$ for some $a_0,a_1$ with $a_0^2\ne 4$, $f_x=3x^2+2(a_0+a_1y)xy^2+y^4$, $f_y=z^2+2a_0x^2y+3a_1x^2y^2+4xy^3$, $f_z=2yz$. If $a_1=0$, then it reduces to the weighted homogeneous trinomial case, and $(V(f)\cap \{z=0\},0)$ is an isolated singularity in $\mathbb{A}^2_{x,y}$, then the Nakai Conjecture holds by proposition \ref{prop3.1}, now we only need to consider the case of $a_1\ne 0$.
By calculation, $R=\mathbb{C}\{x,y,z\}/(f,J(f))=\mathbb{C}\{x,y,z\}/(yz,x^2y^3,x^3-xy^4,a_0x^2y^2+2xy^4,3x^2+2(a_0+a_1y)xy^2+y^4, z^2+2a_0x^2y+3a_1x^2y^2+4xy^3)$, 
it has a $\mathbb{C}$-basis: $\{1,z,z^2,xz,y,xy,x^2y,y^2,xy^2,x^2y^2,y^3,x,x^2\}$, $dim_{\mathbb{C}}R=13$.

If $a_0^2\ne 3$, by solving the equation $Hess(f)\cdot (\alpha_1,\alpha_2,\alpha_3)^{T}=0$ in $R$, we obtain that 
\begin{align*}
\begin{pmatrix}
\alpha_1\\
\alpha_2\\
\alpha_3
\end{pmatrix}&\in Span_{\mathbb{C}}\{
\begin{pmatrix}  z^2 \\ 0 \\ 0 \end{pmatrix},
\begin{pmatrix}  xz \\ 0 \\ 0 \end{pmatrix},
\begin{pmatrix}  x^2 \\ \frac{1}{2}xy \\ 0\end{pmatrix},
\begin{pmatrix}  xy \\ \frac{1}{4}a_1xy+\frac{1}{2}y^2 \\ 0 \end{pmatrix},
\begin{pmatrix}  x^2y \\ 0 \\ 0 \end{pmatrix},
\begin{pmatrix}  xy^2 \\ \frac{1}{2}y^3 \\ 0 \end{pmatrix},
\begin{pmatrix}  x^2y^2 \\ 0 \\ 0 \end{pmatrix},\\&
\begin{pmatrix}  0 \\ z \\ 3a_1x^2y + 2a_0x^2 + 4xy^2\end{pmatrix},
\begin{pmatrix}  0 \\ z^2 \\ 0 \end{pmatrix},
\begin{pmatrix}  0 \\ xz \\ 0 \end{pmatrix},
\begin{pmatrix}  0 \\ x^2 \\ 0\end{pmatrix},
\begin{pmatrix}  0 \\ (a_0^2 - 6)xy - a_0y^3 \\ 0 \end{pmatrix},
\begin{pmatrix}  0 \\ x^2y \\ 0 \end{pmatrix},\\&
\begin{pmatrix}  0 \\ xy^2 \\ 0 \end{pmatrix},
\begin{pmatrix}  0 \\ x^2y^2 \\ 0 \end{pmatrix},
\begin{pmatrix}  0 \\ 0 \\ z^2 \end{pmatrix},
\begin{pmatrix}  0 \\ 0 \\ xz \end{pmatrix},
\begin{pmatrix}  0 \\ 0 \\ x^2y^2 \end{pmatrix}
\}:=W.
\end{align*}
We see that $I_1\subset (z^2,xz,x^2,xy)+(f,J(f))$, $I_1^2\subset (z^2,xz,x^2,xy)^2+(f,J(f))=(x^2y^2,f,J(f))$.
If $a_0^2=3$, similarly by solving the equation $Hess(f)\cdot (\alpha_1,\alpha_2,\alpha_3)^{T}=0$ in $R$, we obtain that 
\begin{align*}
\begin{pmatrix}
\alpha_1\\
\alpha_2\\
\alpha_3
\end{pmatrix}&\in W+\mathbb{C}<\begin{pmatrix}
    y^3 \\ \frac{3a_1}{4}y^3-\frac{a_0}{2}y^2 \\ 0
\end{pmatrix}>.
\end{align*}
We see that $I_1\subset (z^2,xz,x^2,xy,y^3)+(f,J(f))$, $I_1^2\subset (z^2,xz,x^2,xy,y^3)^2+(f,J(f))=(x^2y^2,f,J(f))$.

Now in both cases, we define $(\beta_{ij})\in M_3(R)$ by $$
\beta_{11}=x^2y,$$
$$\beta_{12}=\beta_{21}=\frac{1}{2}xy^2,$$
$$\beta_{22}=\frac{7a_1^2}{16(a_0^2-4)}x^2y-\frac{3a_0a_1}{8(a_0^2-4)}x^2-\frac{(a_0^2-1)a_1}{4(a_0^2-4)}xy^2+\frac{1}{4}y^3,$$
$$\beta_{13}=\beta_{31}=\beta_{23}=\beta_{32}=\beta_{33}=0.$$
Then the lifting $(\widetilde{\beta_{ij}})\in M_3(\mathcal{O}_3)$ satisfies$$
\begin{pmatrix}
\widetilde{\beta_{11}} & \widetilde{\beta_{12}} & \widetilde{\beta_{13}}\\
\widetilde{\beta_{21}} & \widetilde{\beta_{22}} & \widetilde{\beta_{23}}\\
\widetilde{\beta_{31}} & \widetilde{\beta_{32}} & \widetilde{\beta_{33}}
\end{pmatrix}
\begin{pmatrix}
f_{x}\\
f_{y}\\
f_{z}
\end{pmatrix}
\equiv
\begin{pmatrix}
0\\
0\\
0
\end{pmatrix} \ \ \ (  \textup{mod} \ (f,J(f)^2)  ),$$ and $\widetilde{\beta_{11}}\notin (x^2y^2,f,J(f))$, so $\beta_{11}\notin I_1^2 \mod J(f)$.
By proposition \ref{prop2.14}, the Nakai Conjecture holds.
\\

2.$Q_{2,p}$ singularity

For a $Q_{2,p}$ singularity $(V(f),0)$ defined by $f=x^3+yz^2+x^2y^2+(a_0+a_1y)y^{6+p}$ for some $p>0,a_0\ne 0$, $f_x=3x^2+2xy^2$, $f_y=z^2+2x^2y+(7+p)a_1y^{6+p}+(6+p)a_0y^{5+p}$, $f_z=2yz$. As $\frac{1}{3}xf_x+\frac{1}{6}yf_y+\frac{5}{12}zf_z-f=(\frac{p}{6}a_0+\frac{p+1}{6}a_1y)y^{p+6}$, and $\frac{p}{6}a_0+\frac{p+1}{6}a_1y$ is invertible in $R=\mathbb{C}\{x,y,z\}/(f,J(f))$, $y^{p+6}=0$ in $R$. Then $R=\mathbb{C}\{x,y,z\}/(f,J(f))=\mathbb{C}\{x,y,z\}/(y^{p+6},yz,x^3,x^2y^2,3x^2+2xy^2,z^2+2x^2y+(6+p)a_0y^{5+p})$, it has a $\mathbb{C}$-basis: $\{1,z,z^2,y,y^2,\cdots,y^{p+5},x,xz,xy,xy^2\}$, $dim_{\mathbb{C}}R=12+p$.

 By solving the equation $Hess(f)\cdot (\alpha_1,\alpha_2,\alpha_3)^{T}=0$ in $R$, we obtain that 
\begin{align*}
\begin{pmatrix}
\alpha_1\\
\alpha_2\\
\alpha_3
\end{pmatrix}&\in Span_{\mathbb{C}}\{
\begin{pmatrix}  z^2 \\ 0 \\ 0\end{pmatrix},
\begin{pmatrix}  y^{p+4} \\ 0 \\ 0\end{pmatrix},
\begin{pmatrix}  y^{p+5} \\ 0 \\ 0\end{pmatrix},
\begin{pmatrix}  xz \\ 0 \\ 0\end{pmatrix},
\begin{pmatrix}  xy \\ \frac{1}{2}y^2 \\ 0\end{pmatrix},
\begin{pmatrix}  xy^2 \\ 0 \\ 0\end{pmatrix},
\begin{pmatrix}  0 \\ z \\ -\frac{4}{3}xy^2+(p+6)a_0y^{p+5}\end{pmatrix},
\\&
\begin{pmatrix}  0 \\ z^2 \\ 0\end{pmatrix},
\begin{pmatrix}  0 \\ y^3 \\ 0\end{pmatrix},
\begin{pmatrix}  0 \\ y^4 \\ 0\end{pmatrix},\cdots,
\begin{pmatrix}  0 \\ y^{p+5} \\ 0\end{pmatrix},
\begin{pmatrix}  0 \\ xz \\ 0\end{pmatrix},
\begin{pmatrix}  0 \\ xy \\ 0\end{pmatrix},
\begin{pmatrix}  0 \\ xy^2 \\ 0\end{pmatrix},
\begin{pmatrix}  0 \\ 0 \\ z^2\end{pmatrix},
\begin{pmatrix}  0 \\ 0 \\ y^{p+5}\end{pmatrix},
\begin{pmatrix}  0 \\ 0 \\ xz\end{pmatrix}
\}.
\end{align*}
We see that $I_1\subset (z^2,y^{p+4},xz,xy)+(f,J(f))$, $I_1^2\subset (z^2,y^{p+4},xz,xy)^2+(f,J(f))=(f,J(f))$.

Now for $p=1$, we define $(\beta_{ij})\in M_3(R)$ by $$\beta_{11}=\frac{15a_0}{2}y^6+z^2 ,$$
$$\beta_{12}=\beta_{21}=-xy^2-\frac{27a_0}{77}y^5 ,$$
$$\beta_{22}=(\frac{729}{11}a_0-\frac{1359a_1}{77a_0})xy^2+\frac{9}{77}xy+(\frac{3159a_0}{77}-\frac{834a_1}{77a_0})y^4-\frac{30}{77}y^3 ,$$
$$\beta_{13}=\beta_{31}=\beta_{23}=\beta_{32}=\beta_{33}=0 ;$$
and for $p\ge 2$, we define $(\beta_{ij})\in M_3(R)$ by
$$\beta_{11}=(\frac{3}{2}p+6)a_0y^{p+5}+z^2,$$
$$\beta_{12}=\beta_{21}=-xy^2-\frac{9p(p+5)}{2(17p+60)}a_0y^{p+4},$$
$$\beta_{22}=-\frac{9(p^3+18p^2+72p+60)a_1}{p^2(17p+60)a_0}xy^2+\frac{9p}{17p+60}xy-\frac{6(12p^2+67p+60)a_1}{p^2(17p+60)a_0}y^4-\frac{30}{17p+60}y^3,$$
$$\beta_{13}=\beta_{31}=\beta_{23}=\beta_{32}=\beta_{33}=0.$$
Then the lifting $(\widetilde{\beta_{ij}})\in M_3(\mathcal{O}_3)$ satisfies$$
\begin{pmatrix}
\widetilde{\beta_{11}} & \widetilde{\beta_{12}} & \widetilde{\beta_{13}}\\
\widetilde{\beta_{21}} & \widetilde{\beta_{22}} & \widetilde{\beta_{23}}\\
\widetilde{\beta_{31}} & \widetilde{\beta_{32}} & \widetilde{\beta_{33}}
\end{pmatrix}
\begin{pmatrix}
f_{x}\\
f_{y}\\
f_{z}
\end{pmatrix}
\equiv
\begin{pmatrix}
0\\
0\\
0
\end{pmatrix} \ \ \ (  \textup{mod} \ (f,J(f)^2)  ),$$ and $\widetilde{\beta_{11}}\notin (f,J(f))$, so $\beta_{11}\notin I_1^2 \mod J(f)$.
By proposition \ref{prop2.14}, the Nakai Conjecture holds.
\\

3.$S_{1,0}$ singularity

For a $S_{1,0}$ singularity $(V(f),0)$ defined by $f=x^2z+yz^2+y^5+(a_0+a_1y)zy^3$ for some $a_0,a_1$ with $a_0^2\ne 4$, $f_x=2xz$, $f_y=z^2+5y^4+4a_1y^3z+3a_0y^2z$, $f_z=x^2+2yz+(a_0+a_1y)y^3$. If $a_1=0$, then it reduces to the weighted homogeneous trinomial case, and $(V(f)\cap \{x=0\},0)$ is an isolated singularity in $\mathbb{A}^2_{y,z}$, then the Nakai Conjecture holds by proposition \ref{prop3.1}, now we only need to consider the case of $a_1\ne 0$. As $\frac{3}{2}xf_x+yf_y+2zf_z-5f=a_1y^4z$, by calculation, $R=\mathbb{C}\{x,y,z\}/(f,J(f))=\mathbb{C}\{x,y,z\}/(xz,y^4z,yz^2-y^5,2yz^2+a_0y^3z, z^2+5y^4+4a_1y^3z+3a_0y^2z, x^2+2yz+(a_0+a_1y)y^3)$, it has a $\mathbb{C}$-basis: $\{1,z,z^2,y,yz,y^2,y^2z,y^3,y^3z,x,xy,xy^2,xy^3\}$, $dim_{\mathbb{C}}R=13$.
    
If $a_0^2\ne \frac{40}{11}$, by solving the equation $Hess(f)\cdot (\alpha_1,\alpha_2,\alpha_3)^T=0$ in $R$, we obtain that 
\begin{align*}
\begin{pmatrix}
\alpha_1\\
\alpha_2\\
\alpha_3
\end{pmatrix}&\in Span_{\mathbb{C}}\{
\begin{pmatrix} z^2 \\ 0 \\ 0\end{pmatrix},
\begin{pmatrix}  \frac{1}{2}a_0y^3+yz \\ -\frac{1}{3}xy \\ 0  \end{pmatrix},
\begin{pmatrix}  y^2z \\ 0 \\ 0  \end{pmatrix},
\begin{pmatrix}  (a_0^3-4a_0)y^3 \\ (-\frac{2}{3}a_0^2+4)xy \\ -2a_0xy^2  \end{pmatrix},
\begin{pmatrix}   y^3z \\ 0 \\ 0 \end{pmatrix},
\begin{pmatrix}   xy \\ \frac{2}{3}y^2 \\ \frac{1}{6}a_0a_1y^2z+\frac{4}{3}yz \end{pmatrix},\\&
\begin{pmatrix}  xy^2 \\ 0 \\ (-\frac{1}{2}a_0^2+2)y^2z  \end{pmatrix},
\begin{pmatrix}  xy^3 \\ 0 \\ 0  \end{pmatrix},
\begin{pmatrix}   0 \\ \frac{10}{3}y^2+a_0z \\ \frac{16}{3}a_0a_1y^2z+(\frac{3}{2}a_0^2-\frac{10}{3})yz  \end{pmatrix},
\begin{pmatrix}   0 \\ z^2 \\ 0 \end{pmatrix},
\begin{pmatrix}   0 \\ yz \\ -a_0y^2z \end{pmatrix},
\begin{pmatrix}   0 \\ y^2z \\ 0 \end{pmatrix},\\&
\begin{pmatrix}   0 \\ y^3 \\ (\frac{3}{4}a_0^2-1)y^2z \end{pmatrix},
\begin{pmatrix}   0 \\ y^3z \\ 0 \end{pmatrix},
\begin{pmatrix}   0 \\ xy^2 \\ 0 \end{pmatrix},
\begin{pmatrix}   0 \\ xy^3 \\ 0 \end{pmatrix},
\begin{pmatrix}   0 \\ 0 \\ \frac{1}{2}a_0y^2z+z^2 \end{pmatrix},
\begin{pmatrix}   0 \\ 0 \\ y^3z \end{pmatrix},
\begin{pmatrix}   0 \\ 0 \\ xy^3 \end{pmatrix}
\}:=W.
\end{align*}
We see that $I_1\subset (z^2,y^3,yz,xy)+(f,J(f))$, $I_1^2\subset (z^2,y^3,yz,xy)^2+(f,J(f))=(y^3z,f,J(f))$.\\
If $a_0^2=\frac{40}{11}$, similarly by solving the equation $Hess(f)\cdot (\alpha_1,\alpha_2,\alpha_3)^T=0$ in $R$, we obtain that 
\begin{align*}
\begin{pmatrix}
\alpha_1\\
\alpha_2\\
\alpha_3
\end{pmatrix}&\in W+\mathbb{C}<\begin{pmatrix}
    \frac{a_0}{4}y^2+z \\ \frac{1507}{240}a_0a_1xy+\frac{3}{4}x \\
    -\frac{143}{4}a_1xy^2-\frac{11}{8}a_0xy
\end{pmatrix}>.\end{align*}
We see that $I_1\subset (z^2,y^3,yz,xy,\frac{a_0}{4}y^2+z)+(f,J(f))$, $I_1^2\subset (z^2,y^3,yz,xy,\frac{a_0}{4}y^2+z)^2+(f,J(f))=(y^3z,xy^3,(\frac{a_0}{4}y^2+z)^2,f,J(f))$.
Now in both cases, we define $(\beta_{ij})\in M_3(R)$ by $$\beta_{11}=(3a_0^2-10)y^2z+a_0z^2,$$
$$\beta_{12}=\beta_{21}=-\frac{88a_0a_1}{9(a_0^2-4)}xy^3+\frac{10}{3}xy^2,$$
$$\beta_{13}=\beta_{31}=0,$$
$$\beta_{22}=\frac{20}{9}y^3+\frac{196a_0a_1}{135(a_0^2-4)}z^2,$$
$$\beta_{23}=\beta_{32}=(\frac{196a_0^3a_1}{45(a_0^2-4)}-\frac{392a_0a_1}{27(a_0^2-4)})y^3z+\frac{40}{9}y^2z,$$
$$\beta_{33}=-\frac{40a_0}{9}y^3z.$$
Then the lifting $(\widetilde{\beta_{ij}})\in M_3(\mathcal{O}_3)$ satisfies$$
\begin{pmatrix}
\widetilde{\beta_{11}} & \widetilde{\beta_{12}} & \widetilde{\beta_{13}}\\
\widetilde{\beta_{21}} & \widetilde{\beta_{22}} & \widetilde{\beta_{23}}\\
\widetilde{\beta_{31}} & \widetilde{\beta_{32}} & \widetilde{\beta_{33}}
\end{pmatrix}
\begin{pmatrix}
f_{x}\\
f_{y}\\
f_{z}
\end{pmatrix}
\equiv
\begin{pmatrix}
0\\
0\\
0
\end{pmatrix} \ \ \ (  \textup{mod} \ (f,J(f)^2)  ),$$ and in case of $a_0^2\ne \frac{40}{11}$, $\widetilde{\beta_{11}}\notin (y^3z, f,J(f))$, so $\beta_{11}\notin I_1^2 \mod J(f)$; in case of $a_0^2=\frac{40}{11}$, $\widetilde{\beta_{11}}\notin (y^3z,xy^3,(\frac{a_0}{4}y^2+z)^2,f, J(f))$, so $\beta_{11}\notin I_1^2 \mod J(f)$.
By proposition \ref{prop2.14}, the Nakai Conjecture holds.
\\

4.$S_{1,p}$ singularity

For an $S_{1,p}$ singularity $(V(f),0)$ defined by $f=x^2z+yz^2+x^2y^2+(a_0+a_1y)y^{5+p}$ for some $p>0,a_0\ne 0$, $f_x=2xz+2xy^2$, $f_y=z^2+2x^2y+(5+p)a_0y^{4+p}+(6+p)a_1y^{5+p}$, $f_z=x^2+2yz$. As $\frac{3}{2}xf_x+yf_y+2zf_z-5f=(pa_0+(p+1)a_1y)y^{p+5}$, and $pa_0+(p+1)a_1y$ is invertible in $R=\mathbb{C}\{x,y,z\}/(f,J(f))$, $y^{p+5}=0$ in $R$. Then by calculation, $R=\mathbb{C}\{x,y,z\}/(f,J(f))=\mathbb{C}\{x,y,z\}/(y^{p+5},xz+xy^2,x^2+2yz,yz^2,z^2+2x^2y+(5+p)a_0y^{4+p})=\mathbb{C}\{x,y\}/(x^2y^2,x^2z,xz^2,x^2yz,y^3z,yz^2,y^{p+5},xz+xy^2,x^2+2yz,z^2+2x^2y+(5+p)a_0y^{4+p})$, it has a $\mathbb{C}$-basis: $\{1,z,z^2,y,y^2,\cdots,y^{p+3},yz,y^2z,x,xz,xy,xyz\}$, $dim_{\mathbb{C}}R=12+p$.

By solving the equation $Hess(f)\cdot (\alpha_1,\alpha_2,\alpha_3)^{T}=0$ in $R$, we obtain that 
\begin{align*}
\begin{pmatrix}
\alpha_1\\
\alpha_2\\
\alpha_3
\end{pmatrix}&\in Span_{\mathbb{C}}\{
\begin{pmatrix} z^2 \\ 0 \\ 0  \end{pmatrix},
\begin{pmatrix} yz \\ 0 \\-xz  \end{pmatrix},
\begin{pmatrix} y^2z \\ 0 \\ 0  \end{pmatrix},
\begin{pmatrix} y^{p+3} \\ 0 \\ 0  \end{pmatrix},
\begin{pmatrix} xz \\ 0 \\ 0  \end{pmatrix},
\begin{pmatrix} xy \\ \frac{2}{3}y^2 \\ \frac{4}{3}yz  \end{pmatrix},
\begin{pmatrix} xyz \\ 0 \\ 0  \end{pmatrix},\\&
\begin{pmatrix} 0 \\ -\frac{2}{3}y^2+z \\ (p+5)a_0y^{p+3}-\frac{10}{3}yz \end{pmatrix},
\begin{pmatrix} 0 \\ z^2 \\ 0 \end{pmatrix},
\begin{pmatrix} 0 \\ yz \\ 0  \end{pmatrix},
\begin{pmatrix} 0 \\ y^2z \\ 0  \end{pmatrix},
\begin{pmatrix} 0 \\ y^3 \\ 0 \end{pmatrix},
\begin{pmatrix} 0 \\ y^4 \\ 0  \end{pmatrix},\cdots,
\begin{pmatrix} 0 \\ y^{p+3} \\ 0 \end{pmatrix},\\&
\begin{pmatrix} 0 \\ xz \\ 0  \end{pmatrix},
\begin{pmatrix} 0 \\ xy \\ -xz  \end{pmatrix},
\begin{pmatrix} 0 \\ xyz \\ 0  \end{pmatrix},
\begin{pmatrix} 0 \\ 0 \\ z^2  \end{pmatrix},
\begin{pmatrix} 0 \\ 0 \\ y^2z  \end{pmatrix},
\begin{pmatrix} 0 \\ 0 \\ xyz  \end{pmatrix}
\}.
\end{align*}
We see that $I_1\subset (z^2,yz,y^{p+3},xz,xy)+(f,J(f))$, $I_1^2\subset (z^2,yz,y^{p+3},xz,xy)^2+(f,J(f))=(f,J(f))$.

Now we define $(\beta_{ij})\in M_3(R)$ by $$\beta_{11}=xyz,$$
$$\beta_{12}=\beta_{21}=\frac{10(p+3)}{9(p+5)}y^2z-\frac{p}{9(p+5)}x^2,$$
$$\beta_{13}=\beta_{31}=0,$$
$$\beta_{22}=\frac{2}{9}xz,$$
$$\beta_{23}=\beta_{32}=-\frac{4}{9}xyz,$$
$$\beta_{33}=0.$$
Then the lifting $(\widetilde{\beta_{ij}})\in M_3(\mathcal{O}_3)$ satisfies$$
\begin{pmatrix}
\widetilde{\beta_{11}} & \widetilde{\beta_{12}} & \widetilde{\beta_{13}}\\
\widetilde{\beta_{21}} & \widetilde{\beta_{22}} & \widetilde{\beta_{23}}\\
\widetilde{\beta_{31}} & \widetilde{\beta_{32}} & \widetilde{\beta_{33}}
\end{pmatrix}
\begin{pmatrix}
f_{x}\\
f_{y}\\
f_{z}
\end{pmatrix}
\equiv
\begin{pmatrix}
0\\
0\\
0
\end{pmatrix} \ \ \ (  \textup{mod} \ (f,J(f)^2)  ),$$ and $\widetilde{\beta_{11}}\notin (f,J(f))$, so $\beta_{11}\notin I_1^2 \mod J(f)$.
By proposition \ref{prop2.14}, the Nakai Conjecture holds.
\\

5.$S_{1,2q-1}^{\#}$ singularity

For a $S_{1,2q-1}^{\#}$ singularity $(V(f),0)$ defined by $f=x^2z+yz^2+zy^3+(a_0+a_1y)xy^{3+q}$ for some $q>0,a_0\ne 0$, $f_x=2xz+(a_0+a_1y)y^{3+q}$, $f_y=z^2+3y^2z+(q+3)a_0xy^{q+2}+(q+4)a_1xy^{q+3}$, $f_z=x^2+2yz+y^3$. As $\frac{3}{10}xf_x+\frac{1}{5}yf_y+\frac{2}{5}zf_z-f=(\frac{2q-1}{10}a_0+\frac{2q+1}{10}a_1y)xy^{q+3}$, $\frac{2q-1}{10}a_0+\frac{2q+1}{10}a_1y$ is invertible in $R=\mathbb{C}\{x,y,z\}/(f,J(f))$, $xy^{q+3}=0$ in $R$. And in $R$, $0=y^2f_x-\frac{2}{3}xf_y+\frac{2(q+3)a_0}{3}y^{q+2}f_z=(\frac{2q+9}{3}a_0+a_1y)y^{q+5}$, as $\frac{2q+9}{3}a_0+a_1y$ is invertible in $R$, $y^{q+5}=0$ in $R$. Then $R=\mathbb{C}\{x,y,z\}/(y^{q+5}, xy^{q+3},x^2z,yz^2,2xz+(a_0+a_1y)y^{3+q},z^2+3y^2z+(q+3)a_0xy^{q+2},x^2+2yz+y^3)$, it has a $\mathbb{C}$-basis: $\{1,z,z^2,y,yz,y^2,y^2z,y^3,y^4,\cdots,y^{q+4},x,xy,xy^2,\cdots,xy^{q+1}\}$, $dim_{\mathbb{C}} R=2q+11$.

(1) If $q\ge 2$, by solving the equation $Hess(f)\cdot (\alpha_1,\alpha_2,\alpha_3)^{T}=0$ in $R$, we obtain that 
\begin{align*}
\begin{pmatrix}
\alpha_1\\
\alpha_2\\
\alpha_3
\end{pmatrix}&\in Span_{\mathbb{C}}\{
\begin{pmatrix}  z^2 \\ 0 \\ 0\end{pmatrix},
\begin{pmatrix}  yz \\ 0 \\ \frac{1}{2}a_0y^{q+3}\end{pmatrix},
\begin{pmatrix}  y^2z \\ 0 \\ 0 \end{pmatrix},
\begin{pmatrix}  y^3 \\ -\frac{2}{3}xy \\ -\frac{1}{3}a_0y^{q+3}\end{pmatrix},
\begin{pmatrix} y^4 \\ -\frac{2}{3}xy^2 \\ 0 \end{pmatrix},
\begin{pmatrix} y^5 \\ -\frac{2}{3}xy^3 \\ 0\end{pmatrix}, \cdots,\\&
\begin{pmatrix}  y^{q+2} \\ -\frac{2}{3}xy^q \\ 0\end{pmatrix},
\begin{pmatrix}  y^{q+3} \\ 0 \\ 0 \end{pmatrix},
\begin{pmatrix}  y^{q+4} \\ 0 \\ 0 \end{pmatrix},
\begin{pmatrix}  xy \\ -\frac{2q-1}{9}a_0xy^q+\frac{2}{3}y^2 \\ \frac{2q-1}{6}a_0xy^{q+1}+\frac{4}{3}yz \end{pmatrix},
\begin{pmatrix} xy^2 \\ \frac{2}{3}y^3 \\ 0 \end{pmatrix},
\begin{pmatrix}  xy^3 \\ \frac{2}{3}y^4 \\ 0 \end{pmatrix},\cdots,
\begin{pmatrix}  xy^q \\ \frac{2}{3}y^{q+1} \\ 0 \end{pmatrix},\\&
\begin{pmatrix}  xy^{q+1} \\ 0 \\ y^{q+3} \end{pmatrix},
\begin{pmatrix}  0 \\ \frac{4q+15}{6}a_0xy^q+z \\ -\frac{3}{4}a_0xy^{q+1}+\frac{3}{2}yz \end{pmatrix},
\begin{pmatrix}  0 \\ z^2 \\ 0 \end{pmatrix},
\begin{pmatrix}  0 \\ yz \\ 0 \end{pmatrix},
\begin{pmatrix}  0 \\ y^2z \\ 0 \end{pmatrix},
\begin{pmatrix}  0 \\ y^{q+2} \\ -\frac{3}{2}y^{q+3} \end{pmatrix},
\begin{pmatrix}  0 \\ y^{q+3} \\ 0 \end{pmatrix},
\begin{pmatrix}  0 \\ y^{q+4} \\ 0 \end{pmatrix},\\&
\begin{pmatrix}  0 \\ xy^{q+1} \\ 0 \end{pmatrix},
\begin{pmatrix}  0 \\ 0 \\ z^2\end{pmatrix},
\begin{pmatrix}  0 \\ 0 \\ y^2z \end{pmatrix},
\begin{pmatrix}  0 \\ 0 \\ y^{q+4} \end{pmatrix}
\}.
\end{align*}
We see that $I_1\subset (z^2,yz,y^3,xy)+(f,J(f))$, $I_1^2\subset (z^2,yz,y^3,xy)^2+(f,J(f))=(xy^4,y^5,f,J(f))$. 

Now we define $(\beta_{ij})\in M_3(R)$ by $$\beta_{11}=\frac{6}{2q+11}y^4+y^2z,$$
$$\beta_{12}=\beta_{21}=-\frac{4}{2q+11}xy^2-\frac{2q-1}{3(2q+11)}a_0y^{q+3},$$
$$\beta_{13}=\beta_{31}=\frac{2q+3}{2q+11}a_0y^{q+4},$$
\begin{align*}
    \beta_{22}=&-\frac{(20q^3+76q^2+47q-45)a_0}{54(q+2)(2q+11)}xy^{q+1}-\frac{8}{3(2q+11)}y^3-\frac{(20q^2+80q+3)a_1}{18(2q+11)(q+1)a_0}y^2z\\&+\frac{20q^2+32q-21}{18(2q+11)(q+2)}yz,
\end{align*}
$$\beta_{23}=\beta_{32}=-\frac{4(2q+3)}{3(2q+11)}y^2z-\frac{4q^2+4q-3}{9(2q+11)(q+2)}z^2,$$
$$\beta_{33}=0.$$
Then the lifting $(\widetilde{\beta_{ij}})\in M_3(\mathcal{O}_3)$ satisfies$$
\begin{pmatrix}
\widetilde{\beta_{11}} & \widetilde{\beta_{12}} & \widetilde{\beta_{13}}\\
\widetilde{\beta_{21}} & \widetilde{\beta_{22}} & \widetilde{\beta_{23}}\\
\widetilde{\beta_{31}} & \widetilde{\beta_{32}} & \widetilde{\beta_{33}}
\end{pmatrix}
\begin{pmatrix}
f_{x}\\
f_{y}\\
f_{z}
\end{pmatrix}
\equiv
\begin{pmatrix}
0\\
0\\
0
\end{pmatrix} \ \ \ (  \textup{mod} \ (f,J(f)^2)  ),$$ and $\widetilde{\beta_{11}}\notin (xy^4,y^5,f,J(f))$, so $\beta_{11}\notin I_1^2 \mod J(f)$.
By proposition \ref{prop2.14}, the Nakai Conjecture holds.

(2) If $q=1$, similarly 
by solving the equation $Hess(f)\cdot (\alpha_1,\alpha_2,\alpha_3)^{T}=0$ in $R$, we obtain that 
\begin{align*}
\begin{pmatrix}
\alpha_1\\
\alpha_2\\
\alpha_3
\end{pmatrix}&\in Span_{\mathbb{C}}\{
\begin{pmatrix} z^2 \\ 0 \\ 0 \end{pmatrix},
\begin{pmatrix} yz \\ 0 \\ \frac{1}{2}a_0y^4  \end{pmatrix},
\begin{pmatrix} y^2z \\ 0 \\ 0   \end{pmatrix},
\begin{pmatrix} y^3 \\ -\frac{2}{3}xy \\-\frac{1}{3}a_0y^4 \end{pmatrix},
\begin{pmatrix} y^4 \\ 0 \\ 0   \end{pmatrix},
\begin{pmatrix} y^5 \\ 0 \\ 0   \end{pmatrix},
\begin{pmatrix} xy \\ -\frac{1}{9}a_0xy+\frac{2}{3}y^2 \\ \frac{1}{6}a_0xy^2-\frac{1}{18}a_0^2y^4+\frac{4}{3}yz   \end{pmatrix},\\&
\begin{pmatrix}  xy^2 \\ 0 \\ y^4  \end{pmatrix},
\begin{pmatrix}  0 \\ \frac{19}{6}a_0xy+z \\ -\frac{3}{4}a_0xy^2+\frac{19}{12}a_0^2y^4+\frac{3}{2}yz  \end{pmatrix},
\begin{pmatrix}  0 \\ z^2 \\ 0  \end{pmatrix},
\begin{pmatrix}  0 \\ yz \\ 0  \end{pmatrix},
\begin{pmatrix}  0 \\ y^2z \\ 0  \end{pmatrix},
\begin{pmatrix}  0 \\ y^3 \\ -\frac{3}{2}y^4  \end{pmatrix},
\begin{pmatrix}  0 \\ y^4 \\ 0  \end{pmatrix},
\begin{pmatrix}  0 \\ y^5 \\ 0  \end{pmatrix},\\&
\begin{pmatrix}  0 \\ xy^2 \\ 0  \end{pmatrix},
\begin{pmatrix}  0 \\ 0 \\ z^2  \end{pmatrix},
\begin{pmatrix}  0 \\ 0 \\ y^2z  \end{pmatrix},
\begin{pmatrix}  0 \\ 0 \\ y^5  \end{pmatrix}
\}.
\end{align*}
We see that $I_1\subset (z^2,yz,y^3,xy)+(f,J(f))$, $I_1^2\subset (z^2,yz,y^3,xy)^2+(f,J(f))=(xy^4,y^5,f,J(f))$. 

Now we define $(\beta_{ij})\in M_3(R)$ by $$\beta_{11}=\frac{6}{13}y^4+y^2z,$$
$$\beta_{12}=\beta_{21}=-\frac{4}{13}xy^2-\frac{1}{39}a_0y^4+\frac{2167}{11232}a_0y^2z+\frac{2887}{101088}a_0z^2,$$
$$\beta_{13}=\beta_{31}=\frac{5}{13}a_0y^{5},$$
\begin{align*}
    \beta_{22}=&-\frac{49}{1053}a_0xy^{2}-\frac{18805}{151632}a_0^2y^4-\frac{8}{39}y^3-(\frac{42059a_0^2}{303264}+\frac{103a_1}{468a_0})y^2z+\frac{31}{702}yz,
\end{align*}
$$\beta_{23}=\beta_{32}=\frac{1721}{16848}a_0^2y^5-\frac{20}{39}y^2z-\frac{5}{351}z^2,$$
$$\beta_{33}=0.$$
Then the lifting $(\widetilde{\beta_{ij}})\in M_3(\mathcal{O}_3)$ satisfies$$
\begin{pmatrix}
\widetilde{\beta_{11}} & \widetilde{\beta_{12}} & \widetilde{\beta_{13}}\\
\widetilde{\beta_{21}} & \widetilde{\beta_{22}} & \widetilde{\beta_{23}}\\
\widetilde{\beta_{31}} & \widetilde{\beta_{32}} & \widetilde{\beta_{33}}
\end{pmatrix}
\begin{pmatrix}
f_{x}\\
f_{y}\\
f_{z}
\end{pmatrix}
\equiv
\begin{pmatrix}
0\\
0\\
0
\end{pmatrix} \ \ \ (  \textup{mod} \ (f,J(f)^2)  ),$$ and $\widetilde{\beta_{11}}\notin (xy^4,y^5,f,J(f))$, so $\beta_{11}\notin I_1^2 \mod J(f)$.
By proposition \ref{prop2.14}, the Nakai Conjecture holds.
\\

6.$S_{1,2q}^{\#}$ singularity

For a $S_{1,2q}^{\#}$ singularity $(V(f),0)$ defined by $f=x^2z+yz^2+zy^3+(a_0+a_1y)x^2y^{2+q}$ for some $q>0, a_0\ne 0$, $f_x=2xz+2x(a_0+a_1y)y^{q+2}$, $f_y=z^2+3y^2z+(q+2)a_0x^2y^{q+1}+(q+3)a_1x^2y^{q+2}$, $f_z=x^2+2yz+y^3$. As $\frac{3}{10}xf_x+\frac{1}{5}yf_y+\frac{2}{5}zf_z-f=(\frac{q}{5}a_0+\frac{q+1}{5}a_1y)x^2y^{q+2}$, $\frac{q}{5}a_0+\frac{q+1}{5}a_1y$ is invertible in $R=\mathbb{C}\{x,y,z\}/(f,J(f))$, $x^2y^{q+2}=0$ in $R$. Then $R=\mathbb{C}\{x,y,z\}/(x^2y^{q+2},x^2z,yz^2,y^3z,xz^2,xz+x(a_0+a_1y)y^{q+2},z^2+3y^2z+(q+2)a_0x^2y^{q+1},x^2+2yz+y^3)$, by calculation, it has a $\mathbb{C}$-basis: $\{1,z,z^2,y,yz,y^2,y^2z,y^3,y^4,\cdots,y^{q+3},x,xz,xy,xy^2,\cdots,xy^{q+1}\}$, $dim_{\mathbb{C}}R=12+2q$.

By solving the equation $Hess(f)\cdot (\alpha_1,\alpha_2,\alpha_3)^{T}=0$ in $R$, we obtain that 
\begin{align*}
\begin{pmatrix}
\alpha_1\\
\alpha_2\\
\alpha_3
\end{pmatrix}&\in Span_{\mathbb{C}}\{
\begin{pmatrix}  z^2 \\ 0 \\ 0  \end{pmatrix},
\begin{pmatrix}  yz \\ 0 \\ -xz  \end{pmatrix},
\begin{pmatrix}  y^2z \\ 0 \\ 0  \end{pmatrix},
\begin{pmatrix}  y^3 \\ -\frac{2}{3}xy \\ \frac{2}{3}xz \end{pmatrix},
\begin{pmatrix}  y^4 \\ -\frac{2}{3}xy^2 \\ 0  \end{pmatrix},
\begin{pmatrix}  y^5 \\ -\frac{2}{3}xy^3 \\ 0  \end{pmatrix},\cdots,
\begin{pmatrix}  y^{q+2} \\ -\frac{2}{3}xy^q \\ 0  \end{pmatrix},\\&
\begin{pmatrix}  y^{q+3} \\ 0 \\ \frac{1}{a_0}xz  \end{pmatrix},
\begin{pmatrix}  xz \\ 0 \\ 0  \end{pmatrix},
\begin{pmatrix}  xy \\ \frac{4q}{9}a_0y^{q+2}+\frac{2}{3}y^2 \\ -\frac{2q}{3}a_0y^{q+3}+\frac{4}{3}yz \end{pmatrix},
\begin{pmatrix}  xyz \\ 0 \\ 0  \end{pmatrix},
\begin{pmatrix}  xy^2 \\ \frac{2}{3}y^3 \\ 0  \end{pmatrix},
\begin{pmatrix}  xy^3 \\ \frac{2}{3}y^4 \\ 0 \end{pmatrix},\cdots,
\begin{pmatrix}  xy^{q+1} \\ \frac{2}{3}y^{q+2} \\ 0  \end{pmatrix},\\&
\begin{pmatrix}  0 \\ -\frac{2q+7}{3}a_0y^{q+2}+z \\ \frac{3}{2}a_0y^{q+3}+\frac{3}{2}yz  \end{pmatrix},
\begin{pmatrix}  0 \\ z^2 \\ 0  \end{pmatrix},
\begin{pmatrix}  0 \\ yz \\ 0  \end{pmatrix},
\begin{pmatrix}  0 \\ y^2z \\ 0  \end{pmatrix},
\begin{pmatrix}  0 \\ y^{q+3} \\ 0  \end{pmatrix},
\begin{pmatrix}  0 \\ xz \\ 0  \end{pmatrix},
\begin{pmatrix}  0 \\ xyz \\ 0  \end{pmatrix},
\begin{pmatrix}  0 \\ xy^{q+1} \\ \frac{3}{2a_0}xz \end{pmatrix},\\&
\begin{pmatrix} 0 \\ 0 \\ z^2 \end{pmatrix},
\begin{pmatrix}  0 \\ 0 \\ y^2z \end{pmatrix},
\begin{pmatrix}  0 \\ 0 \\ xyz \end{pmatrix}
\}.
\end{align*}
We see that $I_1\subset (z^2,yz,y^3,xz,xy)+(f,J(f))$, and $I_1^2\subset (z^2,yz,y^3,xz,xy)^2+(f,J(f))=(xy^4,y^5,f,J(f))$.

Now we define $(\beta_{ij})\in M_3(R)$ by $$\beta_{11}=-\frac{9}{q+6}y^4+z^2,$$
$$\beta_{12}=\beta_{21}=\frac{6}{q+6}xy^2+\frac{2(q^2+5q+12)}{3(q+6)}xz,$$
$$\beta_{13}=\beta_{31}=\frac{6(q+2)}{q+6}xyz,$$
\begin{align*}
    \beta_{22}=&-\frac{2(9q^3+58q^2+137q+120)}{9(q+6)(2q+5)}a_0y^{q+3}+\frac{4}{q+6}y^3+\frac{2(q+2)(3q^2+15q+20)a_1}{(q+6)(2q+5)(q+3)a_0}y^2z\\&-\frac{2q(5q+13)}{3(q+6)(2q+5)}yz,
\end{align*}
$$\beta_{23}=\beta_{32}=\frac{4(q+2)}{q+6}y^2z+\frac{4q(q+2)}{3(q+6)(2q+5)}z^2,$$
$$\beta_{33}=0.$$
Then the lifting $(\widetilde{\beta_{ij}})\in M_3(\mathcal{O}_3)$ satisfies$$
\begin{pmatrix}
\widetilde{\beta_{11}} & \widetilde{\beta_{12}} & \widetilde{\beta_{13}}\\
\widetilde{\beta_{21}} & \widetilde{\beta_{22}} & \widetilde{\beta_{23}}\\
\widetilde{\beta_{31}} & \widetilde{\beta_{32}} & \widetilde{\beta_{33}}
\end{pmatrix}
\begin{pmatrix}
f_{x}\\
f_{y}\\
f_{z}
\end{pmatrix}
\equiv
\begin{pmatrix}
0\\
0\\
0
\end{pmatrix} \ \ \ (  \textup{mod} \ (f,J(f)^2)  ),$$ and $\widetilde{\beta_{11}}\notin (xy^4,y^5,f,J(f))$, so $\beta_{11}\notin I_1^2 \mod J(f)$.
By proposition \ref{prop2.14}, the Nakai Conjecture holds.
\\

7.$U_{1,0}$ singularity

For a $U_{1,0}$ singularity $(V(f),0)$ defined by $f=x^3+xz^2+xy^3+(a_0+a_1y)y^3z$ for some $a_0,a_1$ with $a_0(a_0^2+1)\ne 0$, $f_x=3x^2+z^2+y^3$, $f_y=3xy^2+3a_0y^2z+4a_1y^3z$, $f_z=2xz+(a_0+a_1y)y^3$. If $a_1=0$, then it reduces to the weighted homogeneous trinomial case, and $(V(f)\cap \{z=0\},0)$ is an isolated singularity in $\mathbb{A}^2_{x,y}$, then the Nakai Conjecture holds by proposition \ref{prop3.1}, now we only need to consider the case of $a_1\ne 0$. By calculation, $R=\mathbb{C}\{x,y,z\}/(f,J(f))=\mathbb{C}\{x,y,z\}/(y^4z,xy^3+a_0y^3z,xz^2+\frac{a_0}{2}y^3z,x^3+xz^2,3x^2+z^2+y^3,3xy^2+3a_0y^2z+4a_1y^3z,2xz+(a_0+a_1y)y^3$, it has a $\mathbb{C}$-basis: $\{1,y,y^2,y^3,y^4,x,xy,xy^2,xy^3,x^2,x^2y,z,yz\}$, $dim_{\mathbb{C}}R=13$.

If $a_0^2\ne \frac{1}{3}$, by solving the equation $Hess(f)\cdot (\alpha_1,\alpha_2,\alpha_3)^T=0$ in $R$, we obtain that 
\begin{align*}
\begin{pmatrix}
\alpha_1\\
\alpha_2\\
\alpha_3
\end{pmatrix}&\in Span_{\mathbb{C}}\{
\begin{pmatrix} y^3 \\ 0 \\ \frac{6a_0^2-2}{a_0(a_0^2+1)}x^2+\frac{4a_0}{a_0^2+1}y^3 \end{pmatrix},
\begin{pmatrix} y^4 \\ 0 \\ 0  \end{pmatrix},
\begin{pmatrix} xy \\ \frac{2}{3}y^2 \\ yz  \end{pmatrix},
\begin{pmatrix} xy^2 \\ 0 \\ 0  \end{pmatrix},
\begin{pmatrix} xy^3 \\ 0 \\ 0  \end{pmatrix},
\begin{pmatrix} x^2 \\ 0 \\ -\frac{4a_0}{a_0^2+1}x^2+\frac{a_0(a_0^2-3)}{2(a_0^2+1)}y^3  \end{pmatrix},\\&
\begin{pmatrix} x^2y \\ 0 \\ 0  \end{pmatrix},
\begin{pmatrix} 0 \\ y^3 \\ 0  \end{pmatrix},
\begin{pmatrix} 0 \\ y^4 \\ 0  \end{pmatrix},
\begin{pmatrix} 0 \\ xy \\ \frac{6a_0}{a_0^2+1}x^2-\frac{3a_0(a_0^2-1)}{2(a_0^2+1)}y^3  \end{pmatrix},
\begin{pmatrix} 0 \\ xy^2 \\ 0  \end{pmatrix},
\begin{pmatrix} 0 \\ xy^3 \\ 0  \end{pmatrix},
\begin{pmatrix} 0 \\ x^2 \\ 0  \end{pmatrix},
\begin{pmatrix} 0 \\ x^2y \\ 0  \end{pmatrix},\\&
\begin{pmatrix} 0 \\ yz \\ -\frac{6}{a_0^2+1}x^2+\frac{3(a_0^2-1)}{2(a_0^2+1)}y^3  \end{pmatrix},
\begin{pmatrix} 0 \\ 0 \\ y^4  \end{pmatrix},
\begin{pmatrix} 0 \\ 0 \\ xy^2  \end{pmatrix},
\begin{pmatrix} 0 \\ 0 \\ xy^3  \end{pmatrix},
\begin{pmatrix} 0 \\ 0 \\ x^2y  \end{pmatrix}
\}:=W.
\end{align*}
If $a_0^2=\frac{1}{3}$, similarly by solving the equation $Hess(f)\cdot (\alpha_1,\alpha_2,\alpha_3)^T=0$ in $R$, we obtain that 
\begin{align*}
\begin{pmatrix}
\alpha_1\\
\alpha_2\\
\alpha_3
\end{pmatrix}&\in
W+\mathbb{C}<\begin{pmatrix}
    y^2 \\ \frac{3}{2}x+\frac{3}{2}a_0z \\ 21a_1x^2+\frac{7}{2}a_1y^3+3a_0y^2
\end{pmatrix},
\begin{pmatrix}
    yz \\ 2a_0y^2 \\ 3xy
\end{pmatrix}>
\end{align*}
We see that in both cases, $I_1\subset (y^2,yz,xy,x^2)+(f,J(f))$, $I_1^2\subset (y^2,yz,xy,x^2)^2+(f,J(f))=(xy^3,y^4,f,J(f))$.

Now we define $(\beta_{ij})\in M_3(R)$ by $$\beta_{11}=x^2y,$$
$$
\beta_{12}=\beta_{21}=\frac{2}{3}xy^2,$$
$$\beta_{13}=\beta_{31}=-\frac{1}{2}a_0y^4,$$
$$\beta_{22}=\frac{8a_1}{27a_0}x^2y+\frac{4}{9}y^3,$$
$$\beta_{23}=\beta_{32}=(-\frac{2}{9}a_1+\frac{8a_1}{9a_0^2})xy^3-\frac{2}{3a_0}xy^2,$$
$$\beta_{33}=-3x^2y-y^4.$$
Then in both cases, the lifting $(\widetilde{\beta_{ij}})\in M_3(\mathcal{O}_3)$ satisfies$$
\begin{pmatrix}
\widetilde{\beta_{11}} & \widetilde{\beta_{12}} & \widetilde{\beta_{13}}\\
\widetilde{\beta_{21}} & \widetilde{\beta_{22}} & \widetilde{\beta_{23}}\\
\widetilde{\beta_{31}} & \widetilde{\beta_{32}} & \widetilde{\beta_{33}}
\end{pmatrix}
\begin{pmatrix}
f_{x}\\
f_{y}\\
f_{z}
\end{pmatrix}
\equiv
\begin{pmatrix}
0\\
0\\
0
\end{pmatrix} \ \ \ (  \textup{mod} \ (f,J(f)^2)  ),$$ and $\widetilde{\beta_{11}}\notin (xy^3,y^4,f,J(f))$, so $\beta_{11}\notin I_1^2 \mod J(f)$.
By proposition \ref{prop2.14}, the Nakai Conjecture holds.
\\

8.$U_{1,2q-1}$ singularity

For a $U_{1,2q-1}$ singularity $(V(f),0)$ defined by $f=x^3+xz^2+xy^3+(a_0+a_1y)y^{1+q}z^2$ for some $q>0, a_0\ne 0$, $f_x=3x^2+z^2+y^3$, $f_y=3xy^2+(q+1)a_0y^qz^2+(q+2)a_1y^{q+1}z^2$, $f_z=2xz+2(a_0+a_1y)y^{q+1}z$. As $\frac{1}{3}xf_x+\frac{2}{9}yf_y+\frac{1}{3}zf_z-f=(\frac{2q-1}{9}a_0+\frac{2q+1}{9}a_1y)y^{q+1}z^2$, $\frac{2q-1}{9}a_0+\frac{2q+1}{9}a_1y$ is invertible in $R=\mathbb{C}\{x,y,z\}/(f,J(f))$, $y^{q+1}z^2=0$ in $R$. Then $R=\mathbb{C}\{x,y,z\}/(y^{q+1}z^2,x^3,xz^2,xy^3,3x^2+z^2+y^3,3xy^2+(q+1)a_0y^qz^2,xz+(a_0+a_1y)y^{q+1}z)$, by calculation, it has a $\mathbb{C}$-basis:$\{1,y,y^2,\cdots,y^{q+3},z,yz,y^2z,\cdots,y^{q+2}z,x,xy,x^2y,x^2\}$, $dim_{\mathbb{C}}R=2q+11$.

By solving the equation $Hess(f)\cdot (\alpha_1,\alpha_2,\alpha_3)^{T}=0$ in $R$, we obtain that 
\begin{align*}
\begin{pmatrix}
\alpha_1\\
\alpha_2\\
\alpha_3
\end{pmatrix}&\in Span_{\mathbb{C}}\{
\begin{pmatrix} -\frac{3}{(2q-1)a_0}xy+y^{q+2} \\  -\frac{2}{(2q-1)a_0}y^2 \\ -\frac{3}{(2q-1)a_0}yz \end{pmatrix},
\begin{pmatrix} y^{q+3} \\ 0 \\ 0  \end{pmatrix},
\begin{pmatrix} x^2 \\ 0 \\ 0  \end{pmatrix},
\begin{pmatrix} x^2y \\ 0 \\ 0  \end{pmatrix},
\begin{pmatrix} y^{q+1}z \\ 0 \\ 0  \end{pmatrix},
\begin{pmatrix} y^{q+2}z \\ 0 \\ 0  \end{pmatrix},\\&
\begin{pmatrix} 0 \\ y^3 \\ \frac{3}{2}y^2z  \end{pmatrix},
\begin{pmatrix} 0 \\ y^4 \\ \frac{3}{2}y^3z  \end{pmatrix},\cdots,
\begin{pmatrix} 0 \\ y^{q+3} \\ \frac{3}{2}y^{q+2}z  \end{pmatrix},
\begin{pmatrix} 0 \\ xy \\ 0  \end{pmatrix},
\begin{pmatrix} 0 \\ x^2 \\ 0  \end{pmatrix},
\begin{pmatrix} 0 \\ x^2y \\ 0  \end{pmatrix},
\begin{pmatrix} 0 \\ yz \\ -\frac{3}{2}y^3  \end{pmatrix},
\begin{pmatrix} 0 \\ y^2z \\ -\frac{3}{2}y^4  \end{pmatrix},\cdots,\\&
\begin{pmatrix} 0 \\ y^qz\\  -\frac{3}{2}y^{q+2} \end{pmatrix},
\begin{pmatrix} 0 \\ y^{q+1}z \\ 0  \end{pmatrix},
\begin{pmatrix} 0 \\ y^{q+2}z \\ 0  \end{pmatrix},
\begin{pmatrix} 0 \\ 0 \\ y^{q+3}  \end{pmatrix},
\begin{pmatrix} 0 \\ 0 \\ x^2  \end{pmatrix},
\begin{pmatrix} 0 \\ 0 \\ x^2y  \end{pmatrix},
\begin{pmatrix} 0 \\ 0 \\ y^{q+1}z  \end{pmatrix},
\begin{pmatrix} 0 \\ 0 \\ y^{q+2}z\end{pmatrix}
\}.
\end{align*}
We see that $I_1\subset (-\frac{3}{(2q-1)a_0}xy+y^{q+2},y^{q+3},x^2,y^{q+1}z)+(f,J(f))$, and $I_1^2\subset (-\frac{3}{(2q-1)a_0}xy+y^{q+2},y^{q+3},x^2,y^{q+1}z)^2+(f,J(f))=(f,J(f))$.

Now for $q=1$, we define $(\beta_{ij})\in M_3(R)$ by $$\beta_{11}=-\frac{3}{2a_0}xy^2+y^4,$$
$$\beta_{12}=\beta_{21}=-\frac{11}{5}xy^2,$$
$$\beta_{13}=\beta_{31}=\frac{33}{10}a_0y^3z,$$
$$\beta_{22}=\frac{11}{30}x^2-\frac{11}{10}y^3,$$
$$\beta_{23}=\beta_{32}=\frac{11}{15}a_0^2y^3z-\frac{33}{20}y^2z,$$
$$\beta_{33}=\frac{33}{8a_0}xy^2-\frac{11}{40}y^4;$$
and for $q\ge 2$, we define $(\beta_{ij})\in M_3(R)$ by
$$\beta_{11}=x^2y,$$
$$\beta_{12}=\beta_{21}=\frac{(4q+7)(q+1)}{9(q+4)}a_0y^{q+3},$$
$$\beta_{13}=\beta_{31}=-\frac{4q+7}{2(q+4)}a_0y^{q+2}z,$$
$$\beta_{22}=-\frac{(4q+7)(2q-1)}{9(q+4)(q+1)}x^2+\frac{4q+7}{3(q+4)(q+1)}y^3,$$
$$\beta_{23}=\beta_{32}=\frac{4q+7}{2(q+4)(q+1)}y^2z,$$
$$\beta_{33}=-\frac{4q+7}{2(q+1)}x^2y-\frac{3(4q+7)}{4(q+1)(q+4)}y^4.$$
Then the lifting $(\widetilde{\beta_{ij}})\in M_3(\mathcal{O}_3)$ satisfies$$
\begin{pmatrix}
\widetilde{\beta_{11}} & \widetilde{\beta_{12}} & \widetilde{\beta_{13}}\\
\widetilde{\beta_{21}} & \widetilde{\beta_{22}} & \widetilde{\beta_{23}}\\
\widetilde{\beta_{31}} & \widetilde{\beta_{32}} & \widetilde{\beta_{33}}
\end{pmatrix}
\begin{pmatrix}
f_{x}\\
f_{y}\\
f_{z}
\end{pmatrix}
\equiv
\begin{pmatrix}
0\\
0\\
0
\end{pmatrix} \ \ \ (  \textup{mod} \ (f,J(f)^2)  ),$$ and $\widetilde{\beta_{11}}\notin (f,J(f))$, so $\beta_{11}\notin I_1^2 \mod J(f)$.
By proposition \ref{prop2.14}, the Nakai Conjecture holds.
\\

9.$U_{1,2q}$ singularity

For a $U_{1,2q}$ singularity $(V(f),0)$ defined by $f=x^3+xz^2+xy^3+(a_0+a_1y)y^{3+q}z$ for some $q>0, a_0\ne 0$, $f_x=3x^2+z^2+y^3$, $f_y=3xy^2+(q+3)a_0y^{q+2}z+(q+4)a_1y^{q+3}z$, $f_z=2xz+(a_0+a_1y)y^{q+3}$. As $\frac{1}{3}xf_x+\frac{2}{9}yf_y+\frac{1}{3}zf_z-f=(\frac{2q}{9}a_0+\frac{2(q+1)}{9}a_1y)y^{q+3}z$, $\frac{2q}{9}a_0+\frac{2(q+1)}{9}a_1y$ is invertible in $R=\mathbb{C}\{x,y,z\}/(f,J(f))$, $y^{q+3}z=0$ in $R$. And in $R$, $x^3=xz^2=xy^3=0$, $y^{q+6}=y^{q+3}(-3x^2-z^2)=0$, $0=3y^2f_z-2zf_y=3a_0y^{q+5}-2(q+3)a_0y^{q+2}z^2=3a_0y^{q+5}+2(q+3)a_0y^{q+2}(3x^2+y^3)=(2q+9)a_0y^{q+5}$, so $y^{q+5}=0$ in $R$. Then $R=\mathbb{C}\{x,y,z\}/(y^{q+5}, y^{q+3}z,x^3,xz^2,xy^3,3x^2+z^2+y^3,3xy^2+(q+3)a_0y^{q+2}z,2xz+(a_0+a_1y)y^{q+3})$, it has a $\mathbb{C}$-basis:$\{1,y,y^2,\cdots,y^{q+4},z,yz,y^2z,\cdots,y^{q+1}z, x,xy,\\ x^2y,xy^2,x^2\}$, $dim_{\mathbb{C}} R=2q+12$.

By solving the equation $Hess(f)\cdot (\alpha_1,\alpha_2,\alpha_3)^{T}=0$ in $R$, we obtain that 
\begin{align*}
\begin{pmatrix}
\alpha_1\\
\alpha_2\\
\alpha_3
\end{pmatrix}&\in Span_{\mathbb{C}}\{
\begin{pmatrix} y^{q+3} \\ 0 \\ 0  \end{pmatrix},
\begin{pmatrix} y^{q+4} \\ 0 \\ 0  \end{pmatrix},
\begin{pmatrix} xy+\frac{q}{3}a_0y^{q+1}z \\ \frac{2}{3}y^2 \\ yz \end{pmatrix},
\begin{pmatrix} xy^2 \\ 0 \\ 0  \end{pmatrix},
\begin{pmatrix} x^2 \\ 0 \\ 0  \end{pmatrix},
\begin{pmatrix} x^2y \\ 0 \\ 0  \end{pmatrix},
\begin{pmatrix} 0 \\ y^3 \\ \frac{3}{2}y^2z \end{pmatrix},\\&
\begin{pmatrix} 0 \\ y^4 \\ \frac{3}{2}y^3z  \end{pmatrix},\cdots,
\begin{pmatrix} 0 \\ y^{q+2} \\ \frac{3}{2}y^{q+1}z  \end{pmatrix},
\begin{pmatrix} 0 \\ y^{q+3} \\ 0  \end{pmatrix},
\begin{pmatrix} 0 \\ y^{q+4} \\ 0  \end{pmatrix},
\begin{pmatrix} 0 \\ xy \\ 0  \end{pmatrix},
\begin{pmatrix} 0 \\ xy^2 \\ 0  \end{pmatrix},
\begin{pmatrix} 0 \\ x^2 \\ 0  \end{pmatrix},
\begin{pmatrix} 0 \\ x^2y \\ 0  \end{pmatrix},
\begin{pmatrix} 0 \\ yz \\ -\frac{3}{2}y^3  \end{pmatrix},\\&
\begin{pmatrix} 0 \\ y^2z \\ -\frac{3}{2}y^4  \end{pmatrix},\cdots,
\begin{pmatrix} 0 \\ y^qz \\ -\frac{3}{2}y^{q+2}  \end{pmatrix},
\begin{pmatrix} 0 \\ y^{q+1}z \\ 0  \end{pmatrix},
\begin{pmatrix} 0 \\ 0 \\ y^{q+3}  \end{pmatrix},
\begin{pmatrix} 0 \\ 0 \\ y^{q+4}  \end{pmatrix},
\begin{pmatrix} 0 \\ 0 \\ xy^2  \end{pmatrix},
\begin{pmatrix} 0 \\ 0 \\ x^2  \end{pmatrix},
\begin{pmatrix} 0 \\ 0 \\ x^2y  \end{pmatrix}
\}.
\end{align*}
We see that $I_1\subset (y^{q+3},xy+\frac{q}{3}a_0y^{q+1}z,xy^2,x^2)+(f,J(f))$, $I_1^2\subset (y^{q+3},xy+\frac{q}{3}a_0y^{q+1}z,xy^2,x^2)^2+(f,J(f))=(f,J(f)) $.

Now we define $(\beta_{ij})\in M_3(R)$ by $$\beta_{11}=x^2y,$$
$$\beta_{12}=\beta_{21}=\frac{2(4q+9)}{3(2q+9)}xy^2,$$
$$\beta_{13}=\beta_{31}=-\frac{4q+9}{2(2q+9)}a_0y^{q+4},$$
$$\beta_{22}=-\frac{8q(4q+9)}{9(2q+3)(2q+9)}x^2+\frac{4(4q+9)}{3(2q+3)(2q+9)}y^3,$$
$$\beta_{23}=\beta_{32}=\frac{2(4q+9)}{(2q+3)(2q+9)}y^2z,$$
$$\beta_{33}=-\frac{4q+9}{2q+3}x^2y-\frac{3(4q+9)}{(2q+3)(2q+9)}y^4.$$
Then the lifting $(\widetilde{\beta_{ij}})\in M_3(\mathcal{O}_3)$ satisfies$$
\begin{pmatrix}
\widetilde{\beta_{11}} & \widetilde{\beta_{12}} & \widetilde{\beta_{13}}\\
\widetilde{\beta_{21}} & \widetilde{\beta_{22}} & \widetilde{\beta_{23}}\\
\widetilde{\beta_{31}} & \widetilde{\beta_{32}} & \widetilde{\beta_{33}}
\end{pmatrix}
\begin{pmatrix}
f_{x}\\
f_{y}\\
f_{z}
\end{pmatrix}
\equiv
\begin{pmatrix}
0\\
0\\
0
\end{pmatrix} \ \ \ (  \textup{mod} \ (f,J(f)^2)  ),$$ and $\widetilde{\beta_{11}}\notin (f,J(f))$, so $\beta_{11}\notin I_1^2 \mod J(f)$.
By proposition \ref{prop2.14}, the Nakai Conjecture holds.
\\

\subsection{Cases of exceptional families of bimodal isolated hypersurface singularities}
\ 
\newline
\indent
Now we treat the exceptional families of unimodal isolated hypersurface singularities. In table \ref{table-3}, we see that if $a_0=a_1=0$, then they are reduced to weighted homogeneous fewnomial singularities (see \ref{thm2.15}), so we only consider the case of $a_0=0  \ \&  \ a_1\ne 0$ and the case of $a_0\ne 0$ in table \ref{table-3}.

1. $E_{18}$ singularity

For a $E_{18}$ singularity $(V(f),0)$ defined by $f=x^3+y^{10}+(a_0+a_1y)xy^7$ for some $a_0,a_1$, we treat the case of $a_0=0  \ \&  \ a_1\ne 0$ and the case of $a_0\ne 0$ separately.

(1) If $a_0\ne 0$, by calculation, the Tjurina algebra $R=\mathcal{O}_2/(f,f_x,f_y)\simeq \mathbb{C}\{x,y\}/(x^3,xy^7,y^{10},3x^2+(a_0+a_1y)y^7,10y^9+7a_0xy^6)$, it has a $\mathbb{C}$-basis: $\{1,y,y^2,\cdots, y^9, x,xy,xy^2,\cdots,xy^5\}$, $dim_{\mathbb{C}} R=16$. 
By solving the equation $Hess(f)\cdot (\alpha_1,\alpha_2)^{T}=0$ in $R$, we obtain that 
\begin{align*}
\begin{pmatrix}
\alpha_1\\
\alpha_2
\end{pmatrix}&\in Span_{\mathbb{C}}\{
\begin{pmatrix} a_0y^6 \\ -\frac{6}{7}x \end{pmatrix}, 
\begin{pmatrix} y^7 \\ 0 \end{pmatrix}, 
\begin{pmatrix} y^8 \\ 0 \end{pmatrix}, 
\begin{pmatrix} y^9 \\ 0 \end{pmatrix}, 
\begin{pmatrix} 35xy \\ a_1x+10y^2 \end{pmatrix}, 
\begin{pmatrix} xy^2 \\ -\frac{1}{5}a_0x \end{pmatrix}, 
\begin{pmatrix} xy^3 \\ 0 \end{pmatrix}, 
\begin{pmatrix} xy^4 \\ 0 \end{pmatrix}, \\&
\begin{pmatrix} xy^5 \\ 0 \end{pmatrix}, 
\begin{pmatrix} 0 \\ a_0x+\frac{10}{7}y^3 \end{pmatrix}, 
\begin{pmatrix} 0 \\ y^4 \end{pmatrix}, 
\begin{pmatrix} 0 \\ y^5 \end{pmatrix}, \cdots,
\begin{pmatrix} 0 \\ y^9 \end{pmatrix}, 
\begin{pmatrix} 0 \\ xy \end{pmatrix}, 
\begin{pmatrix} 0 \\ xy^2 \end{pmatrix},\cdots, 
\begin{pmatrix} 0 \\ xy^5 \end{pmatrix}
\}.
\end{align*}
We see that the ideal $I_1\subset (y^6,xy)+(f,J(f))$, $I_1^2\subset (y^6,xy)^2+(f,J(f))=(y^9, f,J(f))$. 

Now we define $(\beta_{ij})\in M_2(R)$ by
$$\beta_{11}= -\frac{189}{884}a_0^2xy^5 + (a_0^4 - \frac{4995}{442}a_1)y^9 + \frac{2835}{442}a_0y^8, $$
 $$\beta_{12}=\beta_{21}= \frac{22113}{7480}a_0a_1^2xy^5 - \frac{83349}{97240}a_0^3xy^3 - \frac{25515}{4862}xy^2 +
\frac{49329}{9724}a_1^2y^8 + \frac{369117}{97240}a_0a_1y^7 - \frac{11907}{97240}a_0^2y^6, $$
 $$\beta_{22}= -\frac{5103}{1496}a_1xy + \frac{5103}{97240}a_0x - \frac{5103}{48620}a_1^3y^7 + \frac{39123}{48620}a_0a_1^2y^6 +\frac{1701}{48620}a_0^2a_1y^5 - \frac{11907}{48620}a_0^3y^4 - \frac{15309}{9724}y^3. $$
Then the lifting $(\widetilde{\beta_{ij}})\in M_2(\mathcal{O}_2)$ satisfies$$
\begin{pmatrix}
\widetilde{\beta_{11}} & \widetilde{\beta_{12}} \\
\widetilde{\beta_{21}} & \widetilde{\beta_{22}}
\end{pmatrix}
\begin{pmatrix}
f_{x}\\
f_{y}
\end{pmatrix}
\equiv
\begin{pmatrix}
0\\
0
\end{pmatrix} \ \ \ (  \textup{mod} \ (f,J(f)^2)  ),$$ and $\widetilde{\beta_{11}}\notin (y^9, f,J(f))$, so $\beta_{11}\notin I_1^2 \mod J(f)$.
By proposition \ref{prop2.14}, the Nakai Conjecture holds.

(2) If $a_0=0 \ \& \ a_1\ne 0$, by calculation, the Tjurina algebra 
$R=\mathcal{O}_2/(f,f_x,f_y)\simeq \mathbb{C}\{x,y\}/(x^3,xy^8,\\ y^{10},3x^2+a_1y^8,10y^9+8a_1xy^7)$, it has a $\mathbb{C}$-basis: $\{1,y,y^2,\cdots, y^9, x,xy,xy^2,\cdots,xy^6\}$, $dim_{\mathbb{C}} R=17$.
By solving the equation $Hess(f)\cdot (\alpha_1,\alpha_2)^{T}=0$ in $R$, we obtain that 
\begin{align*}
\begin{pmatrix}
\alpha_1\\
\alpha_2
\end{pmatrix}&\in Span_{\mathbb{C}}\{
\begin{pmatrix} a_1y^7 \\ -\frac{3}{4}x \end{pmatrix}, 
\begin{pmatrix} y^8 \\ 0 \end{pmatrix}, 
\begin{pmatrix} y^9 \\ 0 \end{pmatrix}, 
\begin{pmatrix} xy \\ -\frac{1}{5}a_1x \end{pmatrix}, 
\begin{pmatrix} xy^2 \\ 0 \end{pmatrix}, 
\begin{pmatrix} xy^3 \\ 0 \end{pmatrix}, 
\begin{pmatrix} xy^4 \\ 0 \end{pmatrix}, 
\begin{pmatrix} xy^5 \\ 0 \end{pmatrix}, \\&
\begin{pmatrix} xy^6 \\ 0 \end{pmatrix},
\begin{pmatrix} 0 \\ \frac{4}{5}a_1x+y^2 \end{pmatrix}, 
\begin{pmatrix} 0 \\ y^3 \end{pmatrix}, 
\begin{pmatrix} 0 \\ y^4 \end{pmatrix}, \cdots,
\begin{pmatrix} 0 \\ y^9 \end{pmatrix}, 
\begin{pmatrix} 0 \\ xy \end{pmatrix}, 
\begin{pmatrix} 0 \\ xy^2 \end{pmatrix},\cdots, 
\begin{pmatrix} 0 \\ xy^6 \end{pmatrix}
\}.
\end{align*}
We see that the ideal $I_1\subset (y^7,xy)+(f,J(f))$, $I_1^2\subset (y^7,xy)^2+(f,J(f))=(f,J(f))$. 

Now we define $(\beta_{ij})\in M_2(R)$ by
$$\beta_{11}=y^9, $$
$$\beta_{12}=\beta_{21}=-\frac{6}{25}a_1y^8, $$
$$\beta_{22}=\frac{27}{175}xy+\frac{6}{875}a_1^2y^7. $$
Then the lifting $(\widetilde{\beta_{ij}})\in M_2(\mathcal{O}_2)$ satisfies$$
\begin{pmatrix}
\widetilde{\beta_{11}} & \widetilde{\beta_{12}} \\
\widetilde{\beta_{21}} & \widetilde{\beta_{22}}
\end{pmatrix}
\begin{pmatrix}
f_{x}\\
f_{y}
\end{pmatrix}
\equiv
\begin{pmatrix}
0\\
0
\end{pmatrix} \ \ \ (  \textup{mod} \ (f,J(f)^2)  ),$$ and $\widetilde{\beta_{11}}\notin (f,J(f))$, so $\beta_{11}\notin I_1^2 \mod J(f)$.
By proposition \ref{prop2.14}, the Nakai Conjecture holds.
\\

2. $E_{19}$ singularity

For a $E_{19}$ singularity $(V(f),0)$ defined by $f=x^3+xy^7+(a_0+a_1y)y^{11}$ for some $a_0,a_1$, we treat the case of $a_0=0  \ \&  \ a_1\ne 0$ and the case of $a_0\ne 0$ separately.

(1) If $a_0\ne 0$, by calculation, the Tjurina algebra $R=\mathcal{O}_2/(f,f_x,f_y)\simeq \mathbb{C}\{x,y\}/(x^3,xy^7,y^{11},3x^2+y^7,7xy^6+11a_0y^{10})$, it has a $\mathbb{C}$-basis: $\{1,y,y^2,\cdots, y^{10}, x,xy,xy^2,\cdots,xy^5\}$, $dim_{\mathbb{C}} R=17$. 
By solving the equation $Hess(f)\cdot (\alpha_1,\alpha_2)^{T}=0$ in $R$, we obtain that 
\begin{align*}
\begin{pmatrix}
\alpha_1\\
\alpha_2
\end{pmatrix}&\in Span_{\mathbb{C}}\{
\begin{pmatrix} y^6 \\ -\frac{6}{7}x \end{pmatrix}, 
\begin{pmatrix} y^7 \\ 0 \end{pmatrix}, 
\begin{pmatrix} y^8 \\ 0 \end{pmatrix}, 
\begin{pmatrix} y^9 \\ 0 \end{pmatrix}, 
\begin{pmatrix} y^{10} \\ 0 \end{pmatrix}, 
\begin{pmatrix} xy \\ \frac{2}{7}y^2 \end{pmatrix}, 
\begin{pmatrix} xy^2 \\ \frac{2}{7}y^3 \end{pmatrix}, 
\begin{pmatrix} xy^3 \\ \frac{2}{7}y^4 \end{pmatrix}, 
\begin{pmatrix} xy^4 \\ 0 \end{pmatrix}, 
\begin{pmatrix} xy^5 \\ 0 \end{pmatrix}, \\&
\begin{pmatrix} 0 \\ \frac{7}{11}x+a_0y^4 \end{pmatrix}, 
\begin{pmatrix} 0 \\ y^5 \end{pmatrix}, 
\begin{pmatrix} 0 \\ y^6 \end{pmatrix}, \cdots,
\begin{pmatrix} 0 \\ y^{10} \end{pmatrix}, 
\begin{pmatrix} 0 \\ xy \end{pmatrix}, 
\begin{pmatrix} 0 \\ xy^2 \end{pmatrix},\cdots, 
\begin{pmatrix} 0 \\ xy^5 \end{pmatrix}
\}.
\end{align*}
We see that the ideal $I_1\subset (y^6,xy)+(f,J(f))$, $I_1^2\subset (y^6,xy)^2+(f,J(f))=(y^9, f,J(f))$. 

Now we define $(\beta_{ij})\in M_2(R)$ by
$$\beta_{11}=\frac{158123}{30}a_0xy^5 + 173558a_0^2y^9 + \frac{1106861}{45}y^8, $$
\begin{align*}
    \beta_{12}=\beta_{21}=&-\frac{7109740953}{1000}a_1^2xy^5 - \frac{668764929}{100}a_0a_1xy^4 - \frac{1561861}{10}a_0^2xy^3 -\frac{316246}{15}xy^2 \\& + (a_0^3 + \frac{67767147}{50}a_1)y^7 + \frac{45178}{15}a_0y^6, 
\end{align*}
\begin{align*}
    \beta_{22}=&-\frac{729}{7000}a_1^3xy^4 + \frac{81}{1000}a_0a_1^2xy^3 -\frac{9}{100}a_0^2a_1xy^2 + (-\frac{9}{10}a_0^3 - \frac{28946253}{25}a_1)xy \\ &- \frac{6454}{5}a_0x -\frac{888717609}{500}a_1^2y^6 - \frac{45500799}{25}a_0a_1y^5 - \frac{212981}{5}a_0^2y^4 - \frac{90356}{15}y^3. 
\end{align*}
Then the lifting $(\widetilde{\beta_{ij}})\in M_2(\mathcal{O}_2)$ satisfies$$
\begin{pmatrix}
\widetilde{\beta_{11}} & \widetilde{\beta_{12}} \\
\widetilde{\beta_{21}} & \widetilde{\beta_{22}}
\end{pmatrix}
\begin{pmatrix}
f_{x}\\
f_{y}
\end{pmatrix}
\equiv
\begin{pmatrix}
0\\
0
\end{pmatrix} \ \ \ (  \textup{mod} \ (f,J(f)^2)  ),$$ and $\widetilde{\beta_{11}}\notin (y^9, f,J(f))$, so $\beta_{11}\notin I_1^2 \mod J(f)$.
By proposition \ref{prop2.14}, the Nakai Conjecture holds.

(2) If $a_0=0 \ \& \ a_1\ne 0$, by calculation, the Tjurina algebra 
$R=\mathcal{O}_2/(f,f_x,f_y)\simeq \mathbb{C}\{x,y\}/(x^3,xy^7,\\ y^{12},3x^2+y^7,7xy^6+12a_1y^{11})$, it has a $\mathbb{C}$-basis: $\{1,y,y^2,\cdots, y^{11}, x,xy,xy^2,\cdots,xy^5\}$, $dim_{\mathbb{C}} R=18$.
By solving the equation $Hess(f)\cdot (\alpha_1,\alpha_2)^{T}=0$ in $R$, we obtain that 
\begin{align*}
\begin{pmatrix}
\alpha_1\\
\alpha_2
\end{pmatrix}&\in Span_{\mathbb{C}}\{
\begin{pmatrix} y^6 \\ -\frac{6}{7}x \end{pmatrix}, 
\begin{pmatrix} y^7 \\ 0 \end{pmatrix}, 
\begin{pmatrix} y^8 \\ 0 \end{pmatrix}, 
\begin{pmatrix} y^9 \\ 0 \end{pmatrix}, 
\begin{pmatrix} y^{10} \\ 0 \end{pmatrix}, 
\begin{pmatrix} y^{11} \\ 0 \end{pmatrix}, 
\begin{pmatrix} xy \\ \frac{2}{7}y^2 \end{pmatrix}, 
\begin{pmatrix} xy^2 \\ \frac{2}{7}y^3 \end{pmatrix}, 
\begin{pmatrix} xy^3 \\ \frac{2}{7}y^4 \end{pmatrix}, 
\begin{pmatrix} xy^4 \\ \frac{2}{7}y^5 \end{pmatrix},\\& 
\begin{pmatrix} xy^5 \\ 0 \end{pmatrix}, 
\begin{pmatrix} 0 \\ \frac{7}{12}x+a_1y^5 \end{pmatrix}, 
\begin{pmatrix} 0 \\ y^6 \end{pmatrix}, 
\begin{pmatrix} 0 \\ y^7 \end{pmatrix}, \cdots,
\begin{pmatrix} 0 \\ y^{11} \end{pmatrix}, 
\begin{pmatrix} 0 \\ xy \end{pmatrix}, 
\begin{pmatrix} 0 \\ xy^2 \end{pmatrix},\cdots, 
\begin{pmatrix} 0 \\ xy^5 \end{pmatrix}
\}.
\end{align*}
We see that the ideal $I_1\subset (y^6,xy)+(f,J(f))$, $I_1^2\subset (y^6,xy)^2+(f,J(f))=(y^9, f,J(f))$. 

Now we define $(\beta_{ij})\in M_2(R)$ by
$$\beta_{11}=y^8, $$
$$\beta_{12}=\beta_{21}=\frac{81}{686}a_1^2xy^5-\frac{6}{7}xy^2, $$
$$\beta_{22}=-\frac{729}{686}a_1^3xy^4+\frac{54}{343}a_1xy+\frac{81}{343}a_1^2y^6-\frac{12}{49}y^3. $$
Then the lifting $(\widetilde{\beta_{ij}})\in M_2(\mathcal{O}_2)$ satisfies$$
\begin{pmatrix}
\widetilde{\beta_{11}} & \widetilde{\beta_{12}} \\
\widetilde{\beta_{21}} & \widetilde{\beta_{22}}
\end{pmatrix}
\begin{pmatrix}
f_{x}\\
f_{y}
\end{pmatrix}
\equiv
\begin{pmatrix}
0\\
0
\end{pmatrix} \ \ \ (  \textup{mod} \ (f,J(f)^2)  ),$$ and $\widetilde{\beta_{11}}\notin (y^9, f,J(f))$, so $\beta_{11}\notin I_1^2 \mod J(f)$.
By proposition \ref{prop2.14}, the Nakai Conjecture holds.
\\

3. $E_{20}$ singularity

For a $E_{20}$ singularity $(V(f),0)$ defined by $f=x^3+y^{11}+(a_0+a_1y)xy^8$ for some $a_0,a_1$, we treat the case of $a_0=0  \ \&  \ a_1\ne 0$ and the case of $a_0\ne 0$ separately.

(1) If $a_0\ne 0$, by calculation, the Tjurina algebra $R=\mathcal{O}_2/(f,f_x,f_y)\simeq \mathbb{C}\{x,y\}/(x^3,xy^8,y^{11},3x^2+(a_0+a_1y)y^8,11y^{10}+8a_0xy^7)$, it has a $\mathbb{C}$-basis: $\{1,y,y^2,\cdots, y^{10}, x,xy,xy^2,\cdots,xy^5, xy^6\}$, $dim_{\mathbb{C}} R=18$. 
By solving the equation $Hess(f)\cdot (\alpha_1,\alpha_2)^{T}=0$ in $R$, we obtain that 
\begin{align*}
\begin{pmatrix}
\alpha_1\\
\alpha_2
\end{pmatrix}&\in Span_{\mathbb{C}}\{
\begin{pmatrix} a_0y^7 \\ -\frac{3}{4}x \end{pmatrix}, 
\begin{pmatrix} y^8 \\ 0 \end{pmatrix}, 
\begin{pmatrix} y^9 \\ 0 \end{pmatrix}, 
\begin{pmatrix} y^{10} \\ 0 \end{pmatrix}, 
\begin{pmatrix} 44xy \\ a_1x+11y^2 \end{pmatrix}, 
\begin{pmatrix} xy^2 \\ -\frac{2}{11}a_0x \end{pmatrix}, 
\begin{pmatrix} xy^3 \\ 0 \end{pmatrix}, 
\begin{pmatrix} xy^4 \\ 0 \end{pmatrix},\\& 
\begin{pmatrix} xy^5 \\ 0 \end{pmatrix}, 
\begin{pmatrix} xy^6 \\ 0 \end{pmatrix},
\begin{pmatrix} 0 \\ \frac{8}{11}a_0x+y^3 \end{pmatrix}, 
\begin{pmatrix} 0 \\ y^4 \end{pmatrix}, 
\begin{pmatrix} 0 \\ y^5 \end{pmatrix}, \cdots,
\begin{pmatrix} 0 \\ y^{10} \end{pmatrix}, 
\begin{pmatrix} 0 \\ xy \end{pmatrix}, 
\begin{pmatrix} 0 \\ xy^2 \end{pmatrix},\cdots, 
\begin{pmatrix} 0 \\ xy^6 \end{pmatrix}
\}.
\end{align*}
We see that the ideal $I_1\subset (y^7,xy)+(f,J(f))$, $I_1^2\subset (y^7,xy)^2+(f,J(f))=(y^{10}, f,J(f))$. 

Now we define $(\beta_{ij})\in M_2(R)$ by
$$\beta_{11}= \frac{143}{288}a_0^2xy^6 + \frac{64493}{4608}a_1y^{10}- \frac{1573}{192}a_0y^9, $$
$$\beta_{12}=\beta_{21}=a_0^3xy^4 + \frac{363}{64}xy^2 - \frac{371}{64}a_1^2y^9 - \frac{35}{8}a_0a_1y^8 + \frac{1}{4}a_0^2y^7, $$
$$\beta_{22}=\frac{219}{64}a_1xy - \frac{3}{32}a_0x + \frac{85}{704}a_1^3y^8 +
\frac{79}{704}a_0a_1^2y^7 - \frac{1}{32}a_0^2a_1y^6 + \frac{1}{4}a_0^3y^5 + \frac{99}{64}y^3. $$
Then the lifting $(\widetilde{\beta_{ij}})\in M_2(\mathcal{O}_2)$ satisfies$$
\begin{pmatrix}
\widetilde{\beta_{11}} & \widetilde{\beta_{12}} \\
\widetilde{\beta_{21}} & \widetilde{\beta_{22}}
\end{pmatrix}
\begin{pmatrix}
f_{x}\\
f_{y}
\end{pmatrix}
\equiv
\begin{pmatrix}
0\\
0
\end{pmatrix} \ \ \ (  \textup{mod} \ (f,J(f)^2)  ),$$ and $\widetilde{\beta_{11}}\notin (y^{10}, f,J(f))$, so $\beta_{11}\notin I_1^2 \mod J(f)$.
By proposition \ref{prop2.14}, the Nakai Conjecture holds.

(2) If $a_0=0 \ \& \ a_1\ne 0$, by calculation, the Tjurina algebra $R=\mathcal{O}_2/(f,f_x,f_y)\simeq \mathbb{C}\{x,y\}/(x^3,xy^9,\\y^{11},3x^2+a_1y^9,11y^{10}+9a_1xy^8)$, it has a $\mathbb{C}$-basis: $\{1,y,y^2,\cdots, y^{10}, x,xy,xy^2,\cdots,xy^6, xy^7\}$, $dim_{\mathbb{C}} R=19$. 
By solving the equation $Hess(f)\cdot (\alpha_1,\alpha_2)^{T}=0$ in $R$, we obtain that 
\begin{align*}
\begin{pmatrix}
\alpha_1\\
\alpha_2
\end{pmatrix}&\in Span_{\mathbb{C}}\{
\begin{pmatrix} a_1y^8 \\ -\frac{2}{3}x \end{pmatrix}, 
\begin{pmatrix} y^9 \\ 0 \end{pmatrix}, 
\begin{pmatrix} y^{10} \\ 0 \end{pmatrix}, 
\begin{pmatrix} xy \\ -\frac{2}{11}a_1x \end{pmatrix}, 
\begin{pmatrix} xy^2 \\ 0 \end{pmatrix}, 
\begin{pmatrix} xy^3 \\ 0 \end{pmatrix},\cdots, 
\begin{pmatrix} xy^7 \\ 0 \end{pmatrix},\\& 
\begin{pmatrix} 0 \\ \frac{9}{11}a_1x+y^2 \end{pmatrix}, 
\begin{pmatrix} 0 \\ y^3 \end{pmatrix}, 
\begin{pmatrix} 0 \\ y^4 \end{pmatrix}, \cdots,
\begin{pmatrix} 0 \\ y^{10} \end{pmatrix}, 
\begin{pmatrix} 0 \\ xy \end{pmatrix}, 
\begin{pmatrix} 0 \\ xy^2 \end{pmatrix},\cdots, 
\begin{pmatrix} 0 \\ xy^7 \end{pmatrix}
\}.
\end{align*}
We see that the ideal $I_1\subset (y^8,xy)+(f,J(f))$, $I_1^2\subset (y^8,xy)^2+(f,J(f))=(f,J(f))$. 

Now we define $(\beta_{ij})\in M_2(R)$ by
$$\beta_{11}=y^{10}, $$
$$\beta_{12}=\beta_{21}=-\frac{27}{121}a_1y^9, $$
$$\beta_{22}=\frac{243}{1936}xy+\frac{135}{21296}a_1^2y^8. $$
Then the lifting $(\widetilde{\beta_{ij}})\in M_2(\mathcal{O}_2)$ satisfies$$
\begin{pmatrix}
\widetilde{\beta_{11}} & \widetilde{\beta_{12}} \\
\widetilde{\beta_{21}} & \widetilde{\beta_{22}}
\end{pmatrix}
\begin{pmatrix}
f_{x}\\
f_{y}
\end{pmatrix}
\equiv
\begin{pmatrix}
0\\
0
\end{pmatrix} \ \ \ (  \textup{mod} \ (f,J(f)^2)  ),$$ and $\widetilde{\beta_{11}}\notin (f,J(f))$, so $\beta_{11}\notin I_1^2 \mod J(f)$.
By proposition \ref{prop2.14}, the Nakai Conjecture holds.
\\

4. $Z_{17}$ singularity

For a $Z_{17}$ singularity $(V(f),0)$ defined by $f=x^3y+y^8+(a_0+a_1y)xy^6$ for some $a_0,a_1$, we treat the case of $a_0=0  \ \&  \ a_1\ne 0$ and the case of $a_0\ne 0$ separately.

(1) If $a_0\ne 0$, by calculation, the Tjurina algebra $R=\mathcal{O}_2/(f,f_x,f_y)\simeq \mathbb{C}\{x,y\}/(x^3y,xy^6,y^8,\\3x^2y+(a_0+a_1y)y^6,x^3+8y^7+6a_0xy^5)$, it has a $\mathbb{C}$-basis: $\{1,y,y^2,\cdots, y^7, x,xy,xy^2,\cdots,xy^5, x^2\}$, $dim_{\mathbb{C}} R=15$. 
By solving the equation $Hess(f)\cdot (\alpha_1,\alpha_2)^{T}=0$ in $R$, we obtain that 
\begin{align*}
\begin{pmatrix}
\alpha_1\\
\alpha_2
\end{pmatrix}&\in Span_{\mathbb{C}}\{
\begin{pmatrix} a_0^2y^4 \\ \frac{1}{2}a_0x+\frac{12}{5}y^2 \end{pmatrix}, 
\begin{pmatrix} y^5 \\ 0 \end{pmatrix}, 
\begin{pmatrix} y^6 \\ 0 \end{pmatrix}, 
\begin{pmatrix} y^7 \\ 0 \end{pmatrix}, 
\begin{pmatrix} xy \\ \frac{2}{5}y^2 \end{pmatrix}, 
\begin{pmatrix} xy^2 \\ 0 \end{pmatrix}, 
\begin{pmatrix} xy^3 \\ 0 \end{pmatrix}, 
\begin{pmatrix} xy^4 \\ 0 \end{pmatrix},
\begin{pmatrix} xy^5 \\ 0 \end{pmatrix},\\& 
\begin{pmatrix} x^2 \\ 0 \end{pmatrix}, 
\begin{pmatrix} 0 \\ y^3 \end{pmatrix}, 
\begin{pmatrix} 0 \\ y^4 \end{pmatrix}, \cdots,
\begin{pmatrix} 0 \\ y^7 \end{pmatrix}, 
\begin{pmatrix} 0 \\ xy \end{pmatrix}, 
\begin{pmatrix} 0 \\ xy^2 \end{pmatrix},\cdots, 
\begin{pmatrix} 0 \\ xy^5 \end{pmatrix},
\begin{pmatrix} 0 \\ x^2 \end{pmatrix}
\}.
\end{align*}
We see that the ideal $I_1\subset (y^4,xy,x^2)+(f,J(f))$, $I_1^2\subset (y^4,xy,x^2)^2+(f,J(f))=(xy^5,y^7, f,J(f))$. 

Now we define $(\beta_{ij})\in M_2(R)$ by
$$\beta_{11}= -\frac{135}{452}a_0^2xy^4 + (a_0^4 - \frac{4725}{452}a_1)y^7 + \frac{2835}{452}a_0y^6, $$
$$\beta_{12}=\beta_{21}= -\frac{30375}{25312}a_0^3xy^3 - \frac{25515}{3616}xy^2 + \frac{22275}{3164}a_1^2y^7 +
\frac{66825}{12656}a_0a_1y^6 - \frac{6075}{25312}a_0^2y^5, $$ 
$$\beta_{22}= -\frac{171315}{25312}a_1xy^2 + \frac{3645}{25312}a_0xy + \frac{1215}{12656}a_0^2a_1y^5 - \frac{6075}{12656}a_0^3y^4 - \frac{10935}{3616}y^3. $$
Then the lifting $(\widetilde{\beta_{ij}})\in M_2(\mathcal{O}_2)$ satisfies$$
\begin{pmatrix}
\widetilde{\beta_{11}} & \widetilde{\beta_{12}} \\
\widetilde{\beta_{21}} & \widetilde{\beta_{22}}
\end{pmatrix}
\begin{pmatrix}
f_{x}\\
f_{y}
\end{pmatrix}
\equiv
\begin{pmatrix}
0\\
0
\end{pmatrix} \ \ \ (  \textup{mod} \ (f,J(f)^2)  ),$$ and $\widetilde{\beta_{11}}\notin (xy^5, y^7, f,J(f))$, so $\beta_{11}\notin I_1^2 \mod J(f)$.
By proposition \ref{prop2.14}, the Nakai Conjecture holds.

(2) If $a_0=0 \ \& \ a_1\ne 0$, by calculation, the Tjurina algebra $R=\mathcal{O}_2/(f,f_x,f_y)\simeq \mathbb{C}\{x,y\}/(x^3y,\\xy^7,y^8,3x^2y+a_1y^7,x^3+8y^7+7a_1xy^6)$, it has a $\mathbb{C}$-basis: $\{1,y,y^2,\cdots, y^7, x,xy,xy^2,\cdots,xy^6, x^2\}$, $dim_{\mathbb{C}} R=16$. 
By solving the equation $Hess(f)\cdot (\alpha_1,\alpha_2)^{T}=0$ in $R$, we obtain that 
\begin{align*}
\begin{pmatrix}
\alpha_1\\
\alpha_2
\end{pmatrix}&\in Span_{\mathbb{C}}\{
\begin{pmatrix} y^6 \\ 0 \end{pmatrix}, 
\begin{pmatrix} y^7 \\ 0 \end{pmatrix}, 
\begin{pmatrix} xy \\ 0 \end{pmatrix}, 
\begin{pmatrix} xy^2 \\ 0 \end{pmatrix}, \cdots,
\begin{pmatrix} xy^6 \\ 0 \end{pmatrix}, 
\begin{pmatrix} x^2 \\ 0 \end{pmatrix},
\begin{pmatrix} 0 \\ y^2 \end{pmatrix},\\& 
\begin{pmatrix} 0 \\ y^3 \end{pmatrix}, \cdots,
\begin{pmatrix} 0 \\ y^7 \end{pmatrix}, 
\begin{pmatrix} 0 \\ xy \end{pmatrix}, 
\begin{pmatrix} 0 \\ xy^2 \end{pmatrix},\cdots, 
\begin{pmatrix} 0 \\ xy^6 \end{pmatrix},
\begin{pmatrix} 0 \\ x^2 \end{pmatrix}
\}.
\end{align*}
We see that the ideal $I_1\subset (y^6,xy,x^2)+(f,J(f))$, $I_1^2\subset (y^6,xy,x^2)^2+(f,J(f))=(f,J(f))$. 

Now we define $(\beta_{ij})\in M_2(R)$ by
$$\beta_{11}=xy^6, $$
$$\beta_{12}=\beta_{21}=\frac{3}{7}y^7, $$
$$\beta_{22}=-\frac{54}{343}a_1y^7. $$
Then the lifting $(\widetilde{\beta_{ij}})\in M_2(\mathcal{O}_2)$ satisfies$$
\begin{pmatrix}
\widetilde{\beta_{11}} & \widetilde{\beta_{12}} \\
\widetilde{\beta_{21}} & \widetilde{\beta_{22}}
\end{pmatrix}
\begin{pmatrix}
f_{x}\\
f_{y}
\end{pmatrix}
\equiv
\begin{pmatrix}
0\\
0
\end{pmatrix} \ \ \ (  \textup{mod} \ (f,J(f)^2)  ),$$ and $\widetilde{\beta_{11}}\notin (f,J(f))$, so $\beta_{11}\notin I_1^2 \mod J(f)$.
By proposition \ref{prop2.14}, the Nakai Conjecture holds.
\\

5. $Z_{18}$ singularity

For a $Z_{18}$ singularity $(V(f),0)$ defined by $f=x^3y+xy^6+(a_0+a_1y)y^9$ for some $a_0,a_1$, we treat the case of $a_0=0  \ \&  \ a_1\ne 0$ and the case of $a_0\ne 0$ separately.

(1) If $a_0\ne 0$, by calculation, the Tjurina algebra $R=\mathcal{O}_2/(f,f_x,f_y)\simeq \mathbb{C}\{x,y\}/(x^3y,xy^6,y^9,\\3x^2y+y^6,x^3+6xy^5+9a_0y^8)$, it has a $\mathbb{C}$-basis: $\{1,y,y^2,\cdots, y^8, x,xy,xy^2,\cdots,xy^5, x^2\}$, $dim_{\mathbb{C}} R=16$. 
By solving the equation $Hess(f)\cdot (\alpha_1,\alpha_2)^{T}=0$ in $R$, we obtain that 
\begin{align*}
\begin{pmatrix}
\alpha_1\\
\alpha_2
\end{pmatrix}&\in Span_{\mathbb{C}}\{
\begin{pmatrix} y^4 \\ \frac{1}{2}x+\frac{27}{10}a_0y^3 \end{pmatrix}, 
\begin{pmatrix} y^5 \\ 0 \end{pmatrix}, 
\begin{pmatrix} y^6 \\ 0 \end{pmatrix}, 
\begin{pmatrix} y^7 \\ 0 \end{pmatrix}, 
\begin{pmatrix} y^8 \\ 0 \end{pmatrix}, 
\begin{pmatrix} xy \\ \frac{2}{5}y^2 \end{pmatrix}, 
\begin{pmatrix} xy^2 \\ \frac{2}{5}y^3 \end{pmatrix}, 
\begin{pmatrix} xy^3 \\ 0 \end{pmatrix}, 
\begin{pmatrix} xy^4 \\ 0 \end{pmatrix},\\&
\begin{pmatrix} xy^5 \\ 0 \end{pmatrix},
\begin{pmatrix} x^2 \\ 0 \end{pmatrix}, 
\begin{pmatrix} 0 \\ y^4 \end{pmatrix}, 
\begin{pmatrix} 0 \\ y^5 \end{pmatrix}, \cdots,
\begin{pmatrix} 0 \\ y^8 \end{pmatrix}, 
\begin{pmatrix} 0 \\ xy \end{pmatrix}, 
\begin{pmatrix} 0 \\ xy^2 \end{pmatrix},\cdots, 
\begin{pmatrix} 0 \\ xy^5 \end{pmatrix},
\begin{pmatrix} 0 \\ x^2 \end{pmatrix}
\}.
\end{align*}
We see that the ideal $I_1\subset (y^4,xy,x^2)+(f,J(f))$, $I_1^2\subset (y^4,xy,x^2)^2+(f,J(f))=(xy^5,y^7, f,J(f))$. 

Now we define $(\beta_{ij})\in M_2(R)$ by
$$\beta_{11}=-\frac{274374125}{6272}a_1xy^5 + \frac{29875}{21}a_0xy^4 + 34000a_0^2y^7 +\frac{298750}{63}y^6, $$
$$\beta_{12}=\beta_{21}=-\frac{305915}{7}a_0^2xy^3 - \frac{119500}{21}xy^2 + (a_0^3 -\frac{3346025}{196}a_1)y^6 + \frac{23900}{21}a_0y^5, $$ 
\begin{align*}
    \beta_{22}=&-\frac{9}{49}a_0^2a_1xy^3 + (-\frac{9}{7}a_0^3 + \frac{2208375}{98}a_1)xy^2 - \frac{4780}{7}a_0xy - \frac{7950051}{3332}a_1^2y^6\\& + \frac{401517}{196}a_0a_1y^5 - \frac{114718}{7}a_0^2y^4 -\frac{47800}{21}y^3.
\end{align*}
Then the lifting $(\widetilde{\beta_{ij}})\in M_2(\mathcal{O}_2)$ satisfies$$
\begin{pmatrix}
\widetilde{\beta_{11}} & \widetilde{\beta_{12}} \\
\widetilde{\beta_{21}} & \widetilde{\beta_{22}}
\end{pmatrix}
\begin{pmatrix}
f_{x}\\
f_{y}
\end{pmatrix}
\equiv
\begin{pmatrix}
0\\
0
\end{pmatrix} \ \ \ (  \textup{mod} \ (f,J(f)^2)  ),$$ and $\widetilde{\beta_{11}}\notin (xy^5,y^7, f,J(f))$, so $\beta_{11}\notin I_1^2 \mod J(f)$.
By proposition \ref{prop2.14}, the Nakai Conjecture holds.

(2) If $a_0=0 \ \& \ a_1\ne 0$, by calculation, the Tjurina algebra $R=\mathcal{O}_2/(f,f_x,f_y)\simeq \mathbb{C}\{x,y\}/(x^3y,\\xy^6,y^{10},3x^2y+y^6,x^3+6xy^5+10a_1y^9)$, it has a $\mathbb{C}$-basis: $\{1,y,y^2,\cdots, y^9, x,xy,xy^2,\cdots,xy^5, x^2\}$, $dim_{\mathbb{C}} R=17$. 
By solving the equation $Hess(f)\cdot (\alpha_1,\alpha_2)^{T}=0$ in $R$, we obtain that 
\begin{align*}
\begin{pmatrix}
\alpha_1\\
\alpha_2
\end{pmatrix}&\in Span_{\mathbb{C}}\{
\begin{pmatrix} y^4 \\ \frac{1}{2}x+3a_1y^4 \end{pmatrix}, 
\begin{pmatrix} y^5 \\ 0 \end{pmatrix}, 
\begin{pmatrix} y^6 \\ 0 \end{pmatrix}, \cdots, 
\begin{pmatrix} y^9 \\ 0 \end{pmatrix}, 
\begin{pmatrix} xy \\ \frac{2}{5}y^2 \end{pmatrix}, 
\begin{pmatrix} xy^2 \\ \frac{2}{5}y^3 \end{pmatrix}, 
\begin{pmatrix} xy^3 \\ \frac{2}{5}y^4 \end{pmatrix}, 
\begin{pmatrix} xy^4 \\ 0 \end{pmatrix},\\&
\begin{pmatrix} xy^5 \\ 0 \end{pmatrix},
\begin{pmatrix} x^2 \\ 0 \end{pmatrix}, 
\begin{pmatrix} 0 \\ y^5 \end{pmatrix}, 
\begin{pmatrix} 0 \\ y^6 \end{pmatrix}, \cdots,
\begin{pmatrix} 0 \\ y^9 \end{pmatrix}, 
\begin{pmatrix} 0 \\ xy \end{pmatrix}, 
\begin{pmatrix} 0 \\ xy^2 \end{pmatrix},\cdots, 
\begin{pmatrix} 0 \\ xy^5 \end{pmatrix},
\begin{pmatrix} 0 \\ x^2 \end{pmatrix}
\}.
\end{align*}
We see that the ideal $I_1\subset (y^4,xy,x^2)+(f,J(f))$, $I_1^2\subset (y^4,xy,x^2)^2+(f,J(f))=(xy^5,y^7, f,J(f))$. 

Now we define $(\beta_{ij})\in M_2(R)$ by
$$\beta_{11}= -\frac{10}{3}x^2 + a_1xy^4, $$
$$\beta_{12}=\beta_{21}= -\frac{63}{25}a_1^2xy^4 - \frac{4}{3}xy + \frac{4}{5}a_1y^5, $$ 
$$\beta_{22}=-\frac{2349}{625}a_1^3xy^4- \frac{12}{25}a_1xy - \frac{18}{125}a_1^2y^5 - \frac{8}{15}y^2. $$
Then the lifting $(\widetilde{\beta_{ij}})\in M_2(\mathcal{O}_2)$ satisfies$$
\begin{pmatrix}
\widetilde{\beta_{11}} & \widetilde{\beta_{12}} \\
\widetilde{\beta_{21}} & \widetilde{\beta_{22}}
\end{pmatrix}
\begin{pmatrix}
f_{x}\\
f_{y}
\end{pmatrix}
\equiv
\begin{pmatrix}
0\\
0
\end{pmatrix} \ \ \ (  \textup{mod} \ (f,J(f)^2)  ),$$ and $\widetilde{\beta_{11}}\notin (xy^5,y^7, f,J(f))$, so $\beta_{11}\notin I_1^2 \mod J(f)$.
By proposition \ref{prop2.14}, the Nakai Conjecture holds.
\\

6. $Z_{19}$ singularity

For a $Z_{19}$ singularity $(V(f),0)$ defined by $f=x^3y+y^9+(a_0+a_1y)xy^7$ for some $a_0,a_1$, we treat the case of $a_0=0  \ \&  \ a_1\ne 0$ and the case of $a_0\ne 0$ separately.

(1) If $a_0\ne 0$, by calculation, the Tjurina algebra $R=\mathcal{O}_2/(f,f_x,f_y)\simeq \mathbb{C}\{x,y\}/(x^3y,xy^7,y^9,\\3x^2y+(a_0+a_1y)y^7,x^3+9y^8+7a_0xy^6)$, it has a $\mathbb{C}$-basis: $\{1,y,y^2,\cdots, y^8, x,xy,xy^2,\cdots,xy^6, x^2\}$, $dim_{\mathbb{C}} R=17$. 
By solving the equation $Hess(f)\cdot (\alpha_1,\alpha_2)^{T}=0$ in $R$, we obtain that 
\begin{align*}
\begin{pmatrix}
\alpha_1\\
\alpha_2
\end{pmatrix}&\in Span_{\mathbb{C}}\{
\begin{pmatrix} a_0^2y^5 \\ \frac{3}{7}a_0x+\frac{27}{14}y^2 \end{pmatrix}, 
\begin{pmatrix} y^6 \\ 0 \end{pmatrix}, 
\begin{pmatrix} y^7 \\ 0 \end{pmatrix}, 
\begin{pmatrix} y^8 \\ 0 \end{pmatrix}, 
\begin{pmatrix} xy \\ \frac{1}{3}y^2 \end{pmatrix}, 
\begin{pmatrix} xy^2 \\ 0 \end{pmatrix}, 
\begin{pmatrix} xy^3 \\ 0 \end{pmatrix},\cdots, 
\begin{pmatrix} xy^6 \\ 0 \end{pmatrix},\\&
\begin{pmatrix} x^2 \\ 0 \end{pmatrix}, 
\begin{pmatrix} 0 \\ y^3 \end{pmatrix}, 
\begin{pmatrix} 0 \\ y^4 \end{pmatrix}, \cdots,
\begin{pmatrix} 0 \\ y^8 \end{pmatrix}, 
\begin{pmatrix} 0 \\ xy \end{pmatrix}, 
\begin{pmatrix} 0 \\ xy^2 \end{pmatrix},\cdots, 
\begin{pmatrix} 0 \\ xy^6 \end{pmatrix},
\begin{pmatrix} 0 \\ x^2 \end{pmatrix}
\}.
\end{align*}
We see that the ideal $I_1\subset (y^5,xy,x^2)+(f,J(f))$, $I_1^2\subset (y^5,xy,x^2)^2+(f,J(f))=(xy^6,y^8, f,J(f))$. 

Now we define $(\beta_{ij})\in M_2(R)$ by
$$\beta_{11}=\frac{40}{81}a_0^2xy^5 + \frac{4240}{243}a_1y^8 - \frac{160}{27}a_0y^7, $$
$$\beta_{12}=\beta_{21}=a_0^3xy^4 +\frac{16}{3}xy^2 - \frac{277}{36}a_1^2y^8 - \frac{59}{9}a_0a_1y^7 + \frac{1}{3}a_0^2y^6, $$ 
$$\beta_{22}=\frac{79}{12}a_1xy^2 - \frac{1}{6}a_0xy + \frac{1}{3}a_0^3y^5 + 2y^3. $$
Then the lifting $(\widetilde{\beta_{ij}})\in M_2(\mathcal{O}_2)$ satisfies$$
\begin{pmatrix}
\widetilde{\beta_{11}} & \widetilde{\beta_{12}} \\
\widetilde{\beta_{21}} & \widetilde{\beta_{22}}
\end{pmatrix}
\begin{pmatrix}
f_{x}\\
f_{y}
\end{pmatrix}
\equiv
\begin{pmatrix}
0\\
0
\end{pmatrix} \ \ \ (  \textup{mod} \ (f,J(f)^2)  ),$$ and $\widetilde{\beta_{11}}\notin (xy^6,y^8, f,J(f))$, so $\beta_{11}\notin I_1^2 \mod J(f)$.
By proposition \ref{prop2.14}, the Nakai Conjecture holds.

(2) If $a_0=0 \ \& \ a_1\ne 0$, by calculation, the Tjurina algebra $R=\mathcal{O}_2/(f,f_x,f_y)\simeq \mathbb{C}\{x,y\}/(x^3y,\\xy^8,y^9,3x^2y+a_1y^8,x^3+9y^8+8a_1xy^7)$, it has a $\mathbb{C}$-basis: $\{1,y,y^2,\cdots, y^8, x,xy,xy^2,\cdots,xy^7, x^2\}$, $dim_{\mathbb{C}} R=18$. 
By solving the equation $Hess(f)\cdot (\alpha_1,\alpha_2)^{T}=0$ in $R$, we obtain that 
\begin{align*}
\begin{pmatrix}
\alpha_1\\
\alpha_2
\end{pmatrix}&\in Span_{\mathbb{C}}\{ 
\begin{pmatrix} y^7 \\ 0 \end{pmatrix}, 
\begin{pmatrix} y^8 \\ 0 \end{pmatrix}, 
\begin{pmatrix} xy \\ 0 \end{pmatrix}, 
\begin{pmatrix} xy^2 \\ 0 \end{pmatrix},\cdots, 
\begin{pmatrix} xy^7 \\ 0 \end{pmatrix},
\begin{pmatrix} x^2 \\ 0 \end{pmatrix}, 
\begin{pmatrix} 0 \\ y^2 \end{pmatrix}, 
\begin{pmatrix} 0 \\ y^3 \end{pmatrix}, \cdots,\\&
\begin{pmatrix} 0 \\ y^8 \end{pmatrix}, 
\begin{pmatrix} 0 \\ xy \end{pmatrix}, 
\begin{pmatrix} 0 \\ xy^2 \end{pmatrix},\cdots, 
\begin{pmatrix} 0 \\ xy^7 \end{pmatrix},
\begin{pmatrix} 0 \\ x^2 \end{pmatrix}
\}.
\end{align*}
We see that the ideal $I_1\subset (y^7,xy,x^2)+(f,J(f))$, $I_1^2\subset (y^7,xy,x^2)^2+(f,J(f))=(f,J(f))$. 

Now we define $(\beta_{ij})\in M_2(R)$ by
$$\beta_{11}=xy^7, $$
$$\beta_{12}=\beta_{21}=\frac{3}{8}y^8, $$
$$\beta_{22}=-\frac{63}{512}a_1y^8. $$
Then the lifting $(\widetilde{\beta_{ij}})\in M_2(\mathcal{O}_2)$ satisfies$$
\begin{pmatrix}
\widetilde{\beta_{11}} & \widetilde{\beta_{12}} \\
\widetilde{\beta_{21}} & \widetilde{\beta_{22}}
\end{pmatrix}
\begin{pmatrix}
f_{x}\\
f_{y}
\end{pmatrix}
\equiv
\begin{pmatrix}
0\\
0
\end{pmatrix} \ \ \ (  \textup{mod} \ (f,J(f)^2)  ),$$ and $\widetilde{\beta_{11}}\notin (f,J(f))$, so $\beta_{11}\notin I_1^2 \mod J(f)$.
By proposition \ref{prop2.14}, the Nakai Conjecture holds.
\\

7. $W_{17}$ singularity

For a $W_{17}$ singularity $(V(f),0)$ defined by $f=x^4+xy^5+(a_0+a_1y)y^7$ for some $a_0,a_1$, we treat the case of $a_0=0  \ \&  \ a_1\ne 0$ and the case of $a_0\ne 0$ separately.

(1) If $a_0\ne 0$, by calculation, the Tjurina algebra $R=\mathcal{O}_2/(f,f_x,f_y)\simeq \mathbb{C}\{x,y\}/(x^4,xy^5,y^7,\\4x^3+y^5,5xy^4+7a_0y^6)$, it has a $\mathbb{C}$-basis: $\{1,y,y^2,\cdots, y^6, x,xy,xy^2,xy^3,x^2,x^2y,x^2y^2,x^2y^3\}$, $dim_{\mathbb{C}} R=15$. 
By solving the equation $Hess(f)\cdot (\alpha_1,\alpha_2)^{T}=0$ in $R$, we obtain that 
\begin{align*}
\begin{pmatrix}
\alpha_1\\
\alpha_2
\end{pmatrix}&\in Span_{\mathbb{C}}\{
\begin{pmatrix} y^4 \\ 0 \end{pmatrix}, 
\begin{pmatrix} y^5 \\ 0 \end{pmatrix}, 
\begin{pmatrix} y^6 \\ 0 \end{pmatrix}, 
\begin{pmatrix} xy \\ \frac{3}{5}y^2 \end{pmatrix}, 
\begin{pmatrix} xy^2 \\ 0 \end{pmatrix}, 
\begin{pmatrix} xy^3 \\ 0 \end{pmatrix}, 
\begin{pmatrix} x^2 \\ 0 \end{pmatrix}, 
\begin{pmatrix} x^2y \\ 0 \end{pmatrix}, 
\begin{pmatrix} x^2y^2 \\ 0 \end{pmatrix}, 
\begin{pmatrix} x^2y^3 \\ 0 \end{pmatrix}, \\&
\begin{pmatrix} 0 \\ y^3 \end{pmatrix}, 
\begin{pmatrix} 0 \\ y^4 \end{pmatrix}, 
\begin{pmatrix} 0 \\ y^5 \end{pmatrix}, 
\begin{pmatrix} 0 \\ y^6 \end{pmatrix}, 
\begin{pmatrix} 0 \\ xy \end{pmatrix}, 
\begin{pmatrix} 0 \\ xy^2 \end{pmatrix}, 
\begin{pmatrix} 0 \\ xy^3 \end{pmatrix},
\begin{pmatrix} 0 \\ x^2 \end{pmatrix},
\begin{pmatrix} 0 \\ x^2y \end{pmatrix},
\begin{pmatrix} 0 \\ x^2y^2 \end{pmatrix},
\begin{pmatrix} 0 \\ x^2y^3 \end{pmatrix}
\}.
\end{align*}
We see that the ideal $I_1\subset (y^4,xy,x^2)+(f,J(f))$, $I_1^2\subset (y^4,xy,x^2)^2+(f,J(f))=(x^2y^2,y^6,f,J(f))$. 

Now we define $(\beta_{ij})\in M_2(R)$ by
$$\beta_{11}=a_0^3x^2y^2 - \frac{375}{3724}x^2y + \frac{25}{1862}a_0xy^3 +\frac{25}{3724}a_0^2 y^5, $$
$$\beta_{12}=\beta_{21}= \frac{2896}{931}a_0a_1x^2y^2 - \frac{61}{931}a_0^2x^2y+\frac{4}{7}a_0^3xy^3 - \frac{225}{3724}xy^2 - \frac{495}{931}a_1y^5 + \frac{45}{3724}a_0y^4, $$
$$\beta_{22}= \frac{180}{133}a_1x^2y-\frac{9}{931}a_0x^2 + \frac{1716}{931}a_0a_1xy^3 - \frac{24}{931}a_0^2xy^2 + \frac{16}{49}a_0^3y^4 -
\frac{135}{3724}y^3. $$
Then the lifting $(\widetilde{\beta_{ij}})\in M_2(\mathcal{O}_2)$ satisfies$$
\begin{pmatrix}
\widetilde{\beta_{11}} & \widetilde{\beta_{12}} \\
\widetilde{\beta_{21}} & \widetilde{\beta_{22}}
\end{pmatrix}
\begin{pmatrix}
f_{x}\\
f_{y}
\end{pmatrix}
\equiv
\begin{pmatrix}
0\\
0
\end{pmatrix} \ \ \ (  \textup{mod} \ (f,J(f)^2)  ),$$ and $\widetilde{\beta_{11}}\notin (x^2y^2, y^6, f,J(f))$, so $\beta_{11}\notin I_1^2 \mod J(f)$.
By proposition \ref{prop2.14}, the Nakai Conjecture holds.

(2) If $a_0=0 \ \& \ a_1\ne 0$, by calculation, the Tjurina algebra $R=\mathcal{O}_2/(f,f_x,f_y)\simeq \mathbb{C}\{x,y\}/(x^4,\\xy^5,y^8,4x^3+y^5,5xy^4+8a_1y^7)$, it has a $\mathbb{C}$-basis: $\{1,y,\cdots, y^7, x,xy,xy^2,xy^3,x^2,x^2y,x^2y^2,x^2y^3\}$, $dim_{\mathbb{C}} R=16$. 
By solving the equation $Hess(f)\cdot (\alpha_1,\alpha_2)^{T}=0$ in $R$, we obtain that 
\begin{align*}
\begin{pmatrix}
\alpha_1\\
\alpha_2
\end{pmatrix}&\in Span_{\mathbb{C}}\{
\begin{pmatrix} y^4 \\ 0 \end{pmatrix}, 
\begin{pmatrix} y^5 \\ 0 \end{pmatrix}, 
\begin{pmatrix} y^6 \\ 0 \end{pmatrix}, 
\begin{pmatrix} y^7 \\ 0 \end{pmatrix}, 
\begin{pmatrix} xy \\ \frac{3}{5}y^2 \end{pmatrix}, 
\begin{pmatrix} xy^2 \\ \frac{3}{5}y^3 \end{pmatrix}, 
\begin{pmatrix} xy^3 \\ 0 \end{pmatrix}, 
\begin{pmatrix} x^2 \\ 0 \end{pmatrix}, 
\begin{pmatrix} x^2y \\ 0 \end{pmatrix}, \\&
\begin{pmatrix} x^2y^2 \\ 0 \end{pmatrix}, 
\begin{pmatrix} x^2y^3 \\ 0 \end{pmatrix}, 
\begin{pmatrix} 0 \\ y^4 \end{pmatrix}, 
\begin{pmatrix} 0 \\ y^5 \end{pmatrix}, 
\begin{pmatrix} 0 \\ y^6 \end{pmatrix}, 
\begin{pmatrix} 0 \\ y^7 \end{pmatrix}, 
\begin{pmatrix} 0 \\ xy \end{pmatrix}, 
\begin{pmatrix} 0 \\ xy^2 \end{pmatrix}, 
\begin{pmatrix} 0 \\ xy^3 \end{pmatrix},
\begin{pmatrix} 0 \\ x^2 \end{pmatrix},\\&
\begin{pmatrix} 0 \\ x^2y \end{pmatrix},
\begin{pmatrix} 0 \\ x^2y^2 \end{pmatrix},
\begin{pmatrix} 0 \\ x^2y^3 \end{pmatrix}
\}.
\end{align*}
We see that the ideal $I_1\subset (y^4,xy,x^2)+(f,J(f))$, $I_1^2\subset (y^4,xy,x^2)^2+(f,J(f))=(x^2y^2,y^6,f,J(f))$. 

Now we define $(\beta_{ij})\in M_2(R)$ by
$$\beta_{11}=y^5, $$
$$\beta_{12}=\beta_{21}=-\frac{12}{5}x^2y, $$
$$\beta_{22}=-\frac{1024}{625}a_1^2x^2y^3-\frac{36}{25}xy^2. $$
Then the lifting $(\widetilde{\beta_{ij}})\in M_2(\mathcal{O}_2)$ satisfies$$
\begin{pmatrix}
\widetilde{\beta_{11}} & \widetilde{\beta_{12}} \\
\widetilde{\beta_{21}} & \widetilde{\beta_{22}}
\end{pmatrix}
\begin{pmatrix}
f_{x}\\
f_{y}
\end{pmatrix}
\equiv
\begin{pmatrix}
0\\
0
\end{pmatrix} \ \ \ (  \textup{mod} \ (f,J(f)^2)  ),$$ and $\widetilde{\beta_{11}}\notin (x^2y^2, y^6, f,J(f))$, so $\beta_{11}\notin I_1^2 \mod J(f)$.
By proposition \ref{prop2.14}, the Nakai Conjecture holds.
\\

8. $W_{18}$ singularity

For a $W_{18}$ singularity $(V(f),0)$ defined by $f=x^4+y^7+(a_0+a_1y)x^2y^4$ for some $a_0,a_1$, we treat the case of $a_0=0  \ \&  \ a_1\ne 0$ and the case of $a_0\ne 0$ separately.

(1) If $a_0\ne 0$, by calculation, the Tjurina algebra $R=\mathcal{O}_2/(f,f_x,f_y)\simeq \mathbb{C}\{x,y\}/(x^4,y^7,x^2y^4, 4x^3+2(a_0+a_1y)xy^4,7y^6+4a_0x^2y^3)$, it has a $\mathbb{C}$-basis: $\{1,y,y^2,\cdots, y^6, x,xy,xy^2,\cdots, xy^5,x^2,x^2y,x^2y^2\}$, $dim_{\mathbb{C}} R=16$. 
By solving the equation $Hess(f)\cdot (\alpha_1,\alpha_2)^{T}=0$ in $R$, we obtain that 
\begin{align*}
\begin{pmatrix}
\alpha_1\\
\alpha_2
\end{pmatrix}&\in Span_{\mathbb{C}}\{
\begin{pmatrix} y^4 \\ 0 \end{pmatrix}, 
\begin{pmatrix} y^5 \\ 0 \end{pmatrix}, 
\begin{pmatrix} y^6 \\ 0 \end{pmatrix}, 
\begin{pmatrix} xy \\ \frac{1}{2}y^2 \end{pmatrix}, 
\begin{pmatrix} xy^2 \\ 0 \end{pmatrix}, 
\begin{pmatrix} xy^3 \\ 0 \end{pmatrix}, 
\begin{pmatrix} xy^4 \\ 0 \end{pmatrix}, 
\begin{pmatrix} xy^5 \\ 0 \end{pmatrix},
\begin{pmatrix} x^2 \\ 0 \end{pmatrix}, 
\begin{pmatrix} x^2y \\ 0 \end{pmatrix}, \\&
\begin{pmatrix} x^2y^2 \\ 0 \end{pmatrix}, 
\begin{pmatrix} 0 \\ y^3 \end{pmatrix}, 
\begin{pmatrix} 0 \\ y^4 \end{pmatrix}, 
\begin{pmatrix} 0 \\ y^5 \end{pmatrix}, 
\begin{pmatrix} 0 \\ y^6 \end{pmatrix}, 
\begin{pmatrix} 0 \\ xy \end{pmatrix}, 
\begin{pmatrix} 0 \\ xy^2 \end{pmatrix}, \cdots,
\begin{pmatrix} 0 \\ xy^5 \end{pmatrix},
\begin{pmatrix} 0 \\ x^2 \end{pmatrix},
\begin{pmatrix} 0 \\ x^2y \end{pmatrix},
\begin{pmatrix} 0 \\ x^2y^2 \end{pmatrix}
\}.
\end{align*}
We see that the ideal $I_1\subset (y^4,xy,x^2)+(f,J(f))$, $I_1^2\subset (y^4,xy,x^2)^2+(f,J(f))=(xy^5,x^2y^2,f,J(f))$. 

Now we define $(\beta_{ij})\in M_2(R)$ by
$$\beta_{11}=a_0^2x^2y^2 - \frac{7}{2}x^2y - \frac{1}{8}a_1y^6 + \frac{1}{4}a_0y^5, $$
$$\beta_{12}=\beta_{21}= \frac{1}{2}a_0^2xy^3 - 2xy^2, $$
$$\beta_{22}=-\frac{5}{49}a_1x^2y + \frac{4}{49}a_0x^2 +\frac{114}{343}a_0^2y^4 - \frac{8}{7}y^3. $$
Then the lifting $(\widetilde{\beta_{ij}})\in M_2(\mathcal{O}_2)$ satisfies$$
\begin{pmatrix}
\widetilde{\beta_{11}} & \widetilde{\beta_{12}} \\
\widetilde{\beta_{21}} & \widetilde{\beta_{22}}
\end{pmatrix}
\begin{pmatrix}
f_{x}\\
f_{y}
\end{pmatrix}
\equiv
\begin{pmatrix}
0\\
0
\end{pmatrix} \ \ \ (  \textup{mod} \ (f,J(f)^2)  ),$$ and $\widetilde{\beta_{11}}\notin (xy^5, x^2y^2, f,J(f))$, so $\beta_{11}\notin I_1^2 \mod J(f)$.
By proposition \ref{prop2.14}, the Nakai Conjecture holds.

(2) If $a_0=0 \ \& \ a_1\ne 0$, by calculation, the Tjurina algebra $R=\mathcal{O}_2/(f,f_x,f_y)\simeq \mathbb{C}\{x,y\}/(x^4,\\y^7,x^2y^5, 4x^3+2a_1xy^5,7y^6+5a_1x^2y^4)$, it has a $\mathbb{C}$-basis: $\{1,y,y^2,\cdots, y^6, x,xy,xy^2,\cdots, xy^5,x^2,\\x^2y,x^2y^2,x^2y^3\}$, $dim_{\mathbb{C}} R=17$. 
By solving the equation $Hess(f)\cdot (\alpha_1,\alpha_2)^{T}=0$ in $R$, we obtain that 
\begin{align*}
\begin{pmatrix}
\alpha_1\\
\alpha_2
\end{pmatrix}&\in Span_{\mathbb{C}}\{ 
\begin{pmatrix} y^5 \\ 0 \end{pmatrix}, 
\begin{pmatrix} y^6 \\ 0 \end{pmatrix}, 
\begin{pmatrix} xy \\ 0 \end{pmatrix}, 
\begin{pmatrix} xy^2 \\ 0 \end{pmatrix}, \cdots,
\begin{pmatrix} xy^5 \\ 0 \end{pmatrix},
\begin{pmatrix} x^2 \\ 0 \end{pmatrix}, 
\begin{pmatrix} x^2y \\ 0 \end{pmatrix}, 
\begin{pmatrix} x^2y^2 \\ 0 \end{pmatrix}, 
\begin{pmatrix} x^2y^3 \\ 0 \end{pmatrix},\\&
\begin{pmatrix} 0 \\ y^2 \end{pmatrix}, 
\begin{pmatrix} 0 \\ y^3 \end{pmatrix},\cdots,  
\begin{pmatrix} 0 \\ y^6 \end{pmatrix}, 
\begin{pmatrix} 0 \\ xy \end{pmatrix}, 
\begin{pmatrix} 0 \\ xy^2 \end{pmatrix}, \cdots,
\begin{pmatrix} 0 \\ xy^5 \end{pmatrix},
\begin{pmatrix} 0 \\ x^2 \end{pmatrix},
\begin{pmatrix} 0 \\ x^2y \end{pmatrix},
\begin{pmatrix} 0 \\ x^2y^2 \end{pmatrix},
\begin{pmatrix} 0 \\ x^2y^3 \end{pmatrix}
\}.
\end{align*}
We see that the ideal $I_1\subset (y^5,xy,x^2)+(f,J(f))$, $I_1^2\subset (y^5,xy,x^2)^2+(f,J(f))=(x^2y^2,f,J(f))$. 

Now we define $(\beta_{ij})\in M_2(R)$ by
$$\beta_{11}= -\frac{14}{3}x^2 + a_1y^5, $$
$$\beta_{12}=\beta_{21}= -\frac{2}{7}a_1^2xy^4 - \frac{8}{3}xy, $$
$$\beta_{22}= \frac{16}{49}a_1x^2 + \frac{24}{343}a_1^2 y^5 -
\frac{32}{21}y^2. $$ 
Then the lifting $(\widetilde{\beta_{ij}})\in M_2(\mathcal{O}_2)$ satisfies$$
\begin{pmatrix}
\widetilde{\beta_{11}} & \widetilde{\beta_{12}} \\
\widetilde{\beta_{21}} & \widetilde{\beta_{22}}
\end{pmatrix}
\begin{pmatrix}
f_{x}\\
f_{y}
\end{pmatrix}
\equiv
\begin{pmatrix}
0\\
0
\end{pmatrix} \ \ \ (  \textup{mod} \ (f,J(f)^2)  ),$$ and $\widetilde{\beta_{11}}\notin (x^2y^2, f,J(f))$, so $\beta_{11}\notin I_1^2 \mod J(f)$.
By proposition \ref{prop2.14}, the Nakai Conjecture holds.
\\

9. $Q_{16}$ singularity

For a $Q_{16}$ singularity $(V(f),0)$ defined by $f=x^3+yz^2+y^7+(a_0+a_1y)xy^5$ for some $a_0,a_1$, we treat the case of $a_0=0  \ \&  \ a_1\ne 0$ and the case of $a_0\ne 0$ separately.

(1) If $a_0\ne 0$, by calculation, the Tjurina algebra $R=\mathcal{O}_3/(f,f_x,f_y,f_z)\simeq \mathbb{C}\{x,y,z\}/(x^3,xy^5,\\y^7,yz,3x^2+(a_0+a_1y)y^5, z^2+7y^6+5a_0xy^4)$, it has a $\mathbb{C}$-basis: $\{1,z,z^2,y,y^2,\cdots, y^6, x,xz,xy,xy^2,\\xy^3\}$, $dim_{\mathbb{C}} R=14$. 
By solving the equation $Hess(f)\cdot (\alpha_1,\alpha_2,\alpha_3)^{T}=0$ in $R$, we obtain that 
\begin{align*}
\begin{pmatrix}
\alpha_1\\
\alpha_2\\
\alpha_3
\end{pmatrix}&\in Span_{\mathbb{C}}\{ 
\begin{pmatrix} z^2 \\ 0 \\ 0 \end{pmatrix},
\begin{pmatrix} y^5 \\ 0 \\ 0 \end{pmatrix},
\begin{pmatrix} y^6 \\ 0 \\ 0 \end{pmatrix},
\begin{pmatrix} xz \\ 0 \\ 0 \end{pmatrix},
\begin{pmatrix} xy \\ \frac{2}{5}y^2 \\ 0 \end{pmatrix},
\begin{pmatrix} xy^2 \\ 0 \\ 0 \end{pmatrix},
\begin{pmatrix} xy^3 \\ 0 \\ 0 \end{pmatrix},
\begin{pmatrix} 0 \\ z \\ 5a_0xy^3+7y^5 \end{pmatrix},
\begin{pmatrix} 0 \\ z^2 \\ 0 \end{pmatrix},\\&
\begin{pmatrix} 0 \\ y^3 \\ 0 \end{pmatrix},
\begin{pmatrix} 0 \\ y^4 \\ 0 \end{pmatrix},
\begin{pmatrix} 0 \\ y^5 \\ 0 \end{pmatrix},
\begin{pmatrix} 0 \\ y^6 \\ 0 \end{pmatrix},
\begin{pmatrix} 0 \\ xz \\ 0 \end{pmatrix},
\begin{pmatrix} 0 \\ xy \\ 0 \end{pmatrix},
\begin{pmatrix} 0 \\ xy^2 \\ 0 \end{pmatrix},
\begin{pmatrix} 0 \\ xy^3 \\ 0 \end{pmatrix},
\begin{pmatrix} 0 \\ 0 \\ z^2 \end{pmatrix},
\begin{pmatrix} 0 \\ 0 \\ y^6 \end{pmatrix},
\begin{pmatrix} 0 \\ 0 \\ xz \end{pmatrix}
\}.
\end{align*}
We see that the ideal $I_1\subset (z^2,y^5,xz,xy)+(f,J(f))$, $I_1^2\subset (z^2,y^5,xz,xy)^2+(f,J(f))=(f,J(f))$.

Now we define $(\beta_{ij})\in M_3(R)$ by
$$\beta_{11}=-\frac{5432}{975}a_0y^6 - \frac{56}{975}a_0z^2, $$
$$\beta_{12}=\beta_{21}=a_0^3xy^3 + \frac{378}{65}xy^2 -\frac{115}{13}a_0a_1y^6 + \frac{14}{65}a_0^2y^5, $$
$$\beta_{13}=\beta_{31}=0, $$ 
$$\beta_{22}=\frac{693}{65}a_1xy^2 - \frac{9}{65}a_0xy + \frac{2}{5}a_0^3y^4 + \frac{162}{65}y^3, $$
$$\beta_{23}=\beta_{32}= 0, $$
$$\beta_{33}=0. $$
Then the lifting $(\widetilde{\beta_{ij}})\in M_3(\mathcal{O}_3)$ satisfies$$
\begin{pmatrix}
\widetilde{\beta_{11}} & \widetilde{\beta_{12}} & \widetilde{\beta_{13}}\\
\widetilde{\beta_{21}} & \widetilde{\beta_{22}} & \widetilde{\beta_{23}}\\
\widetilde{\beta_{31}} & \widetilde{\beta_{32}} & \widetilde{\beta_{33}}
\end{pmatrix}
\begin{pmatrix}
f_{x}\\
f_{y}\\
f_{z}
\end{pmatrix}
\equiv
\begin{pmatrix}
0\\
0\\
0
\end{pmatrix} \ \ \ (  \textup{mod} \ (f,J(f)^2)  ),$$ and $\widetilde{\beta_{11}}\notin (f,J(f))$, so $\beta_{11}\notin I_1^2 \mod J(f)$.
By proposition \ref{prop2.14}, the Nakai Conjecture holds.

(2) If $a_0=0 \ \& \ a_1\ne 0$, by calculation, the Tjurina algebra $R=\mathcal{O}_3/(f,f_x,f_y,f_z)\simeq \mathbb{C}\{x,y,z\}/(x^3,xy^6,y^7,yz,3x^2+a_1y^6, z^2+7y^6+6a_1xy^5)$, it has a $\mathbb{C}$-basis: $\{1,z,z^2,y,y^2,\cdots, y^6, x,\\xz,xy,xy^2,xy^3,xy^4\}$, $dim_{\mathbb{C}} R=15$. 
By solving the equation $Hess(f)\cdot (\alpha_1,\alpha_2,\alpha_3)^{T}=0$ in $R$, we obtain that 
\begin{align*}
\begin{pmatrix}
\alpha_1\\
\alpha_2\\
\alpha_3
\end{pmatrix}&\in Span_{\mathbb{C}}\{ 
\begin{pmatrix} z^2 \\ 0 \\ 0 \end{pmatrix},
\begin{pmatrix} y^6 \\ 0 \\ 0 \end{pmatrix},
\begin{pmatrix} xz \\ 0 \\ 0 \end{pmatrix},
\begin{pmatrix} xy \\ 0 \\ 0 \end{pmatrix},
\begin{pmatrix} xy^2 \\ 0 \\ 0 \end{pmatrix},
\begin{pmatrix} xy^3 \\ 0 \\ 0 \end{pmatrix},
\begin{pmatrix} xy^4 \\ 0 \\ 0 \end{pmatrix},
\begin{pmatrix} 0 \\ z \\ 6a_1xy^4+7y^5 \end{pmatrix},
\begin{pmatrix} 0 \\ z^2 \\ 0 \end{pmatrix},\\&
\begin{pmatrix} 0 \\ y^2 \\ 0 \end{pmatrix},
\begin{pmatrix} 0 \\ y^3 \\ 0 \end{pmatrix},\cdots,
\begin{pmatrix} 0 \\ y^6 \\ 0 \end{pmatrix},
\begin{pmatrix} 0 \\ xz \\ 0 \end{pmatrix},
\begin{pmatrix} 0 \\ xy \\ 0 \end{pmatrix},
\begin{pmatrix} 0 \\ xy^2 \\ 0 \end{pmatrix},
\begin{pmatrix} 0 \\ xy^3 \\ 0 \end{pmatrix},
\begin{pmatrix} 0 \\ xy^4 \\ 0 \end{pmatrix},
\begin{pmatrix} 0 \\ 0 \\ z^2 \end{pmatrix},
\begin{pmatrix} 0 \\ 0 \\ y^6 \end{pmatrix},
\begin{pmatrix} 0 \\ 0 \\ xz \end{pmatrix}
\}.
\end{align*}
We see that the ideal $I_1\subset (z^2,y^6,xz,xy)+(f,J(f))$, $I_1^2\subset (z^2,y^6,xz,xy)^2+(f,J(f))=(f,J(f))$.

Now we define $(\beta_{ij})\in M_3(R)$ by
$$\beta_{11}=34a_1y^6 + a_1z^2, $$
$$\beta_{12}=\beta_{21}=-\frac{243}{11}xy - \frac{306}{77}a_1^2y^6, $$
$$\beta_{13}=\beta_{31}= -\frac{729}{11}xz, $$ 
$$\beta_{22}=\frac{162}{77}a_1xy + \frac{288}{539}a_1^3y^6
- \frac{729}{77}y^2, $$
$$\beta_{23}=\beta_{32}=0, $$
$$\beta_{33}=-\frac{243}{11}y^6 - \frac{972}{11}z^2. $$
Then the lifting $(\widetilde{\beta_{ij}})\in M_3(\mathcal{O}_3)$ satisfies$$
\begin{pmatrix}
\widetilde{\beta_{11}} & \widetilde{\beta_{12}} & \widetilde{\beta_{13}}\\
\widetilde{\beta_{21}} & \widetilde{\beta_{22}} & \widetilde{\beta_{23}}\\
\widetilde{\beta_{31}} & \widetilde{\beta_{32}} & \widetilde{\beta_{33}}
\end{pmatrix}
\begin{pmatrix}
f_{x}\\
f_{y}\\
f_{z}
\end{pmatrix}
\equiv
\begin{pmatrix}
0\\
0\\
0
\end{pmatrix} \ \ \ (  \textup{mod} \ (f,J(f)^2)  ),$$ and $\widetilde{\beta_{11}}\notin (f,J(f))$, so $\beta_{11}\notin I_1^2 \mod J(f)$.
By proposition \ref{prop2.14}, the Nakai Conjecture holds.
\\

10. $Q_{17}$ singularity

For a $Q_{17}$ singularity $(V(f),0)$ defined by $f=x^3+yz^2+xy^5+(a_0+a_1y)y^8$ for some $a_0,a_1$, we treat the case of $a_0=0  \ \&  \ a_1\ne 0$ and the case of $a_0\ne 0$ separately.

(1) If $a_0\ne 0$, by calculation, the Tjurina algebra $R=\mathcal{O}_3/(f,f_x,f_y,f_z)\simeq \mathbb{C}\{x,y,z\}/(x^3,xy^5,y^8,\\yz,3x^2+y^5, z^2+5xy^4+8a_0y^7)$, it has a $\mathbb{C}$-basis: $\{1,z,z^2,y,y^2,\cdots, y^7, x,xz,xy,xy^2,xy^3\}$, $dim_{\mathbb{C}} R=15$. 
By solving the equation $Hess(f)\cdot (\alpha_1,\alpha_2,\alpha_3)^{T}=0$ in $R$, we obtain that 
\begin{align*}
\begin{pmatrix}
\alpha_1\\
\alpha_2\\
\alpha_3
\end{pmatrix}&\in Span_{\mathbb{C}}\{ 
\begin{pmatrix} z^2 \\ 0 \\ 0 \end{pmatrix},
\begin{pmatrix} y^5 \\ 0 \\ 0 \end{pmatrix},
\begin{pmatrix} y^6 \\ 0 \\ 0 \end{pmatrix},
\begin{pmatrix} y^7 \\ 0 \\ 0 \end{pmatrix},
\begin{pmatrix} xz \\ 0 \\ 0 \end{pmatrix},
\begin{pmatrix} xy \\ \frac{2}{5}y^2 \\ 0 \end{pmatrix},
\begin{pmatrix} xy^2 \\ \frac{2}{5}y^3 \\ 0 \end{pmatrix},
\begin{pmatrix} xy^3 \\ 0 \\ 0 \end{pmatrix},
\begin{pmatrix} 0 \\ z \\ 5xy^3+8a_0y^6 \end{pmatrix},\\&
\begin{pmatrix} 0 \\ z^2 \\ 0 \end{pmatrix},
\begin{pmatrix} 0 \\ y^4 \\ 0 \end{pmatrix},
\begin{pmatrix} 0 \\ y^5 \\ 0 \end{pmatrix},\cdots,
\begin{pmatrix} 0 \\ y^7 \\ 0 \end{pmatrix},
\begin{pmatrix} 0 \\ xz \\ 0 \end{pmatrix},
\begin{pmatrix} 0 \\ xy \\ 0 \end{pmatrix},
\begin{pmatrix} 0 \\ xy^2 \\ 0 \end{pmatrix},
\begin{pmatrix} 0 \\ xy^3 \\ 0 \end{pmatrix},
\begin{pmatrix} 0 \\ 0 \\ z^2 \end{pmatrix},
\begin{pmatrix} 0 \\ 0 \\ y^7 \end{pmatrix},
\begin{pmatrix} 0 \\ 0 \\ xz \end{pmatrix}
\}.
\end{align*}
We see that the ideal $I_1\subset (z^2,y^5,xz,xy)+(f,J(f))$, $I_1^2\subset (z^2,y^5,xz,xy)^2+(f,J(f))=(y^7, f,J(f))$.

Now we define $(\beta_{ij})\in M_3(R)$ by
 $$\beta_{11}=a_0^2y^7 + \frac{1300}{8691}y^6 - \frac{30}{2897}a_0z^2, $$
 $$\beta_{12}=\beta_{21}= -\frac{4032}{2897}a_0^2 xy^3 - \frac{520}{2897}xy^2 - \frac{1569}{2897}a_1y^6 + \frac{112}{2897}a_0y^5, $$
 $$\beta_{13}=\beta_{31}=0, $$ 
 $$\beta_{22}=\frac{2070}{2897}a_1xy^2 - \frac{72}{2897}a_0xy +\frac{189}{2897}a_0a_1y^5- \frac{1512}{2897}a_0^2y^4 - \frac{208}{2897}y^3, $$
 $$\beta_{23}=\beta_{32}=0, $$
 $$\beta_{33}=0. $$
 Then the lifting $(\widetilde{\beta_{ij}})\in M_3(\mathcal{O}_3)$ satisfies$$
\begin{pmatrix}
\widetilde{\beta_{11}} & \widetilde{\beta_{12}} & \widetilde{\beta_{13}}\\
\widetilde{\beta_{21}} & \widetilde{\beta_{22}} & \widetilde{\beta_{23}}\\
\widetilde{\beta_{31}} & \widetilde{\beta_{32}} & \widetilde{\beta_{33}}
\end{pmatrix}
\begin{pmatrix}
f_{x}\\
f_{y}\\
f_{z}
\end{pmatrix}
\equiv
\begin{pmatrix}
0\\
0\\
0
\end{pmatrix} \ \ \ (  \textup{mod} \ (f,J(f)^2)  ),$$ and $\widetilde{\beta_{11}}\notin (y^7, f,J(f))$, so $\beta_{11}\notin I_1^2 \mod J(f)$.
By proposition \ref{prop2.14}, the Nakai Conjecture holds.

(2) If $a_0=0 \ \& \ a_1\ne 0$, by calculation, the Tjurina algebra $R=\mathcal{O}_3/(f,f_x,f_y,f_z)\simeq \mathbb{C}\{x,y,z\}/(x^3,xy^5,y^9,yz,3x^2+y^5, z^2+5xy^4+9a_1y^8)$, it has a $\mathbb{C}$-basis: $\{1,z,z^2,y,y^2,\cdots, y^8, x,\\xz,xy,xy^2,xy^3\}$, $dim_{\mathbb{C}} R=16$.  
By solving the equation $Hess(f)\cdot (\alpha_1,\alpha_2,\alpha_3)^{T}=0$ in $R$, we obtain that
\begin{align*}
\begin{pmatrix}
\alpha_1\\
\alpha_2\\
\alpha_3
\end{pmatrix}&\in Span_{\mathbb{C}}\{ 
\begin{pmatrix} z^2 \\ 0 \\ 0 \end{pmatrix},
\begin{pmatrix} y^5 \\ 0 \\ 0 \end{pmatrix},
\begin{pmatrix} y^6 \\ 0 \\ 0 \end{pmatrix},
\begin{pmatrix} y^7 \\ 0 \\ 0 \end{pmatrix},
\begin{pmatrix} y^8 \\ 0 \\ 0 \end{pmatrix},
\begin{pmatrix} xz \\ 0 \\ 0 \end{pmatrix},
\begin{pmatrix} xy \\ \frac{2}{5}y^2 \\ 0 \end{pmatrix},
\begin{pmatrix} xy^2 \\ \frac{2}{5}y^3 \\ 0 \end{pmatrix},
\begin{pmatrix} xy^3 \\ \frac{2}{5}y^4 \\ 0 \end{pmatrix},\\&
\begin{pmatrix} 0 \\ z \\ 5xy^3+9a_1y^7 \end{pmatrix},
\begin{pmatrix} 0 \\ z^2 \\ 0 \end{pmatrix},
\begin{pmatrix} 0 \\ y^5 \\ 0 \end{pmatrix},
\begin{pmatrix} 0 \\ y^6 \\ 0 \end{pmatrix},
\begin{pmatrix} 0 \\ y^7 \\ 0 \end{pmatrix},
\begin{pmatrix} 0 \\ y^8 \\ 0 \end{pmatrix},
\begin{pmatrix} 0 \\ xz \\ 0 \end{pmatrix},
\begin{pmatrix} 0 \\ xy \\ 0 \end{pmatrix},
\begin{pmatrix} 0 \\ xy^2 \\ 0 \end{pmatrix},\\&
\begin{pmatrix} 0 \\ xy^3 \\ 0 \end{pmatrix},
\begin{pmatrix} 0 \\ 0 \\ z^2 \end{pmatrix},
\begin{pmatrix} 0 \\ 0 \\ y^8 \end{pmatrix},
\begin{pmatrix} 0 \\ 0 \\ xz \end{pmatrix}
\}.
\end{align*}
We see that the ideal $I_1\subset (z^2,y^5,xz,xy)+(f,J(f))$, $I_1^2\subset (z^2,y^5,xz,xy)^2+(f,J(f))=(y^7, f,J(f))$.

Now we define $(\beta_{ij})\in M_3(R)$ by
$$\beta_{11}=-\frac{130}{27}y^5 + a_1z^2, $$
$$\beta_{12}=\beta_{21}=\frac{52}{9}xy - \frac{56}{15}a_1y^5 -
\frac{549}{125}a_1^2z^2, $$ 
$$\beta_{13}=\beta_{31}= \frac{169}{9}xz, $$ 
$$\beta_{22}= \frac{12}{5}a_1xy +
\frac{618}{125}a_1^2y^5 + \frac{104}{45}y^2 + \frac{21627}{3125}a_1^3z^2, $$
$$\beta_{23}=\beta_{32}= 0, $$
$$\beta_{33}=-\frac{169}{30}a_1y^8 + \frac{2197}{90}z^2. $$
 Then the lifting $(\widetilde{\beta_{ij}})\in M_3(\mathcal{O}_3)$ satisfies$$
\begin{pmatrix}
\widetilde{\beta_{11}} & \widetilde{\beta_{12}} & \widetilde{\beta_{13}}\\
\widetilde{\beta_{21}} & \widetilde{\beta_{22}} & \widetilde{\beta_{23}}\\
\widetilde{\beta_{31}} & \widetilde{\beta_{32}} & \widetilde{\beta_{33}}
\end{pmatrix}
\begin{pmatrix}
f_{x}\\
f_{y}\\
f_{z}
\end{pmatrix}
\equiv
\begin{pmatrix}
0\\
0\\
0
\end{pmatrix} \ \ \ (  \textup{mod} \ (f,J(f)^2)  ),$$ and $\widetilde{\beta_{11}}\notin (y^7, f,J(f))$, so $\beta_{11}\notin I_1^2 \mod J(f)$.
By proposition \ref{prop2.14}, the Nakai Conjecture holds.
\\

11. $Q_{18}$ singularity

For a $Q_{18}$ singularity $(V(f),0)$ defined by $f=x^3+yz^2+y^8+(a_0+a_1y)xy^6$ for some $a_0,a_1$, we treat the case of $a_0=0  \ \&  \ a_1\ne 0$ and the case of $a_0\ne 0$ separately.

(1) If $a_0\ne 0$, by calculation, the Tjurina algebra $R=\mathcal{O}_3/(f,f_x,f_y,f_z)\simeq \mathbb{C}\{x,y,z\}/(x^3,xy^6,y^8,\\yz,3x^2+(a_0+a_1y)y^6, z^2+8y^7+6a_0xy^5)$, it has a $\mathbb{C}$-basis: $\{1,z,z^2,y,y^2,\cdots, y^7, x,xz,xy,xy^2,xy^3,\\xy^4\}$, $dim_{\mathbb{C}} R=16$. 
By solving the equation $Hess(f)\cdot (\alpha_1,\alpha_2,\alpha_3)^{T}=0$ in $R$, we obtain that 
\begin{align*}
\begin{pmatrix}
\alpha_1\\
\alpha_2\\
\alpha_3
\end{pmatrix}&\in Span_{\mathbb{C}}\{ 
\begin{pmatrix} z^2 \\ 0 \\ 0 \end{pmatrix},
\begin{pmatrix} y^6 \\ 0 \\ 0 \end{pmatrix},
\begin{pmatrix} y^7 \\ 0 \\ 0 \end{pmatrix},
\begin{pmatrix} xz \\ 0 \\ 0 \end{pmatrix},
\begin{pmatrix} xy \\ \frac{1}{3}y^2 \\ 0 \end{pmatrix},
\begin{pmatrix} xy^2 \\ 0 \\ 0 \end{pmatrix},
\begin{pmatrix} xy^3 \\ 0 \\ 0 \end{pmatrix},
\begin{pmatrix} xy^4 \\ 0 \\ 0 \end{pmatrix},\\&
\begin{pmatrix} 0 \\ z \\ 6a_0xy^4+8y^6 \end{pmatrix},
\begin{pmatrix} 0 \\ z^2 \\ 0 \end{pmatrix},
\begin{pmatrix} 0 \\ y^3 \\ 0 \end{pmatrix},
\begin{pmatrix} 0 \\ y^4 \\ 0 \end{pmatrix},
\cdots,
\begin{pmatrix} 0 \\ y^7 \\ 0 \end{pmatrix},
\begin{pmatrix} 0 \\ xz \\ 0 \end{pmatrix},
\begin{pmatrix} 0 \\ xy \\ 0 \end{pmatrix},
\begin{pmatrix} 0 \\ xy^2 \\ 0 \end{pmatrix},\\&
\begin{pmatrix} 0 \\ xy^3 \\ 0 \end{pmatrix},
\begin{pmatrix} 0 \\ xy^4 \\ 0 \end{pmatrix},
\begin{pmatrix} 0 \\ 0 \\ z^2 \end{pmatrix},
\begin{pmatrix} 0 \\ 0 \\ y^7 \end{pmatrix},
\begin{pmatrix} 0 \\ 0 \\ xz \end{pmatrix}
\}.
\end{align*}
We see that the ideal $I_1\subset (z^2,y^6,xz,xy)+(f,J(f))$, $I_1^2\subset (z^2,y^6,xz,xy)^2+(f,J(f))=(f,J(f))$.

Now we define $(\beta_{ij})\in M_3(R)$ by
$$\beta_{11}= \frac{1420}{9}a_0y^7 + \frac{20}{9}a_0z^2, $$
$$\beta_{12}=\beta_{21}=-126 xy^2 + \frac{625}{4}a_0a_1y^7 -
\frac{17}{2}a_0^2y^6, $$
$$\beta_{13}=\beta_{31}=0, $$
$$\beta_{22}= -\frac{315}{2}a_1 xy^2 +
\frac{9}{2}a_0xy + a_0^3y^5 - \frac{189}{4}y^3, $$
$$\beta_{23}=\beta_{32}=0, $$
$$\beta_{33}=0. $$
Then the lifting $(\widetilde{\beta_{ij}})\in M_3(\mathcal{O}_3)$ satisfies$$
\begin{pmatrix}
\widetilde{\beta_{11}} & \widetilde{\beta_{12}} & \widetilde{\beta_{13}}\\
\widetilde{\beta_{21}} & \widetilde{\beta_{22}} & \widetilde{\beta_{23}}\\
\widetilde{\beta_{31}} & \widetilde{\beta_{32}} & \widetilde{\beta_{33}}
\end{pmatrix}
\begin{pmatrix}
f_{x}\\
f_{y}\\
f_{z}
\end{pmatrix}
\equiv
\begin{pmatrix}
0\\
0\\
0
\end{pmatrix} \ \ \ (  \textup{mod} \ (f,J(f)^2)  ),$$ and $\widetilde{\beta_{11}}\notin (f,J(f))$, so $\beta_{11}\notin I_1^2 \mod J(f)$.
By proposition \ref{prop2.14}, the Nakai Conjecture holds.

(2) If $a_0=0 \ \& \ a_1\ne 0$, by calculation, the Tjurina algebra $R=\mathcal{O}_3/(f,f_x,f_y,f_z)\simeq \mathbb{C}\{x,y,z\}/(x^3,xy^7,y^8,yz,3x^2+a_1y^7, z^2+8y^7+7a_1xy^6)$, it has a $\mathbb{C}$-basis: $\{1,z,z^2,y,y^2,\cdots, y^7, x,\\xz,xy,xy^2,\cdots,xy^5\}$, $dim_{\mathbb{C}} R=17$. 
By solving the equation $Hess(f)\cdot (\alpha_1,\alpha_2,\alpha_3)^{T}=0$ in $R$, we obtain that 
\begin{align*}
\begin{pmatrix}
\alpha_1\\
\alpha_2\\
\alpha_3
\end{pmatrix}&\in Span_{\mathbb{C}}\{ 
\begin{pmatrix} z^2 \\ 0 \\ 0 \end{pmatrix},
\begin{pmatrix} y^7 \\ 0 \\ 0 \end{pmatrix},
\begin{pmatrix} xz \\ 0 \\ 0 \end{pmatrix},
\begin{pmatrix} xy \\ 0 \\ 0 \end{pmatrix},
\begin{pmatrix} xy^2 \\ 0 \\ 0 \end{pmatrix},
\cdots,
\begin{pmatrix} xy^5 \\ 0 \\ 0 \end{pmatrix},
\begin{pmatrix} 0 \\ z \\ 7a_1xy^5+8y^6 \end{pmatrix},\\&
\begin{pmatrix} 0 \\ z^2 \\ 0 \end{pmatrix},
\begin{pmatrix} 0 \\ y^2 \\ 0 \end{pmatrix},
\begin{pmatrix} 0 \\ y^3 \\ 0 \end{pmatrix},
\cdots,
\begin{pmatrix} 0 \\ y^7 \\ 0 \end{pmatrix},
\begin{pmatrix} 0 \\ xz \\ 0 \end{pmatrix},
\begin{pmatrix} 0 \\ xy \\ 0 \end{pmatrix},
\begin{pmatrix} 0 \\ xy^2 \\ 0 \end{pmatrix},
\cdots,
\begin{pmatrix} 0 \\ xy^5 \\ 0 \end{pmatrix},
\begin{pmatrix} 0 \\ 0 \\ z^2 \end{pmatrix},\\&
\begin{pmatrix} 0 \\ 0 \\ y^7 \end{pmatrix},
\begin{pmatrix} 0 \\ 0 \\ xz \end{pmatrix}
\}.
\end{align*}
We see that the ideal $I_1\subset (z^2,y^7,xz,xy)+(f,J(f))$, $I_1^2\subset (z^2,y^7,xz,xy)^2+(f,J(f))=(f,J(f))$.

Now we define $(\beta_{ij})\in M_3(R)$ by
$$\beta_{11}= \frac{187}{5}a_1y^7 + a_1z^2, $$
$$\beta_{12}=\beta_{21}=-\frac{1323}{65}xy - \frac{105}{26}a_1^2y^7, $$
$$\beta_{13}=\beta_{31}=-\frac{9261}{130}xz, $$  
$$\beta_{22}= \frac{189}{104}a_1xy +
\frac{399}{832}a_1^3y^7 - \frac{3969}{520}y^2, $$
$$\beta_{23}=\beta_{32}=0, $$
$$\beta_{33}= -\frac{1323}{52}y^7 - \frac{25137}{260}z^2. $$
Then the lifting $(\widetilde{\beta_{ij}})\in M_3(\mathcal{O}_3)$ satisfies$$
\begin{pmatrix}
\widetilde{\beta_{11}} & \widetilde{\beta_{12}} & \widetilde{\beta_{13}}\\
\widetilde{\beta_{21}} & \widetilde{\beta_{22}} & \widetilde{\beta_{23}}\\
\widetilde{\beta_{31}} & \widetilde{\beta_{32}} & \widetilde{\beta_{33}}
\end{pmatrix}
\begin{pmatrix}
f_{x}\\
f_{y}\\
f_{z}
\end{pmatrix}
\equiv
\begin{pmatrix}
0\\
0\\
0
\end{pmatrix} \ \ \ (  \textup{mod} \ (f,J(f)^2)  ),$$ and $\widetilde{\beta_{11}}\notin (f,J(f))$, so $\beta_{11}\notin I_1^2 \mod J(f)$.
By proposition \ref{prop2.14}, the Nakai Conjecture holds.
\\

12. $S_{16}$ singularity

For a $S_{16}$ singularity $(V(f),0)$ defined by $f=x^2z+yz^2+xy^4+(a_0+a_1y)y^6$ for some $a_0,a_1$, we treat the case of $a_0=0  \ \&  \ a_1\ne 0$ and the case of $a_0\ne 0$ separately.

(1) If $a_0\ne 0$, by calculation, the Tjurina algebra $R=\mathcal{O}_3/(f,f_x,f_y,f_z)\simeq \mathbb{C}\{x,y,z\}/(x^2z,yz^2,\\xy^4,y^6,2xz+y^4, z^2+4xy^3+6a_0y^5,x^2+2yz)$, it has a $\mathbb{C}$-basis: $\{1,z,z^2,y,yz,y^2,y^2z,y^3,y^3z,y^4,y^5,\\x,xy,xy^2\}$, $dim_{\mathbb{C}} R=14$. 
By solving the equation $Hess(f)\cdot (\alpha_1,\alpha_2,\alpha_3)^{T}=0$ in $R$, we obtain that 
\begin{align*}
\begin{pmatrix}
\alpha_1\\
\alpha_2\\
\alpha_3
\end{pmatrix}&\in Span_{\mathbb{C}}\{ 
\begin{pmatrix} z^2 \\ 0 \\ 0 \end{pmatrix},
\begin{pmatrix} yz \\ 0 \\ \frac{1}{2}y^4 \end{pmatrix},
\begin{pmatrix} y^2z \\ 0 \\ 0 \end{pmatrix},
\begin{pmatrix} y^3 \\ \frac{1}{2}z \\ xy^2+3a_0y^4 \end{pmatrix},
\begin{pmatrix} y^3z \\ 0 \\ 0 \end{pmatrix},
\begin{pmatrix} y^4 \\ 0 \\ 0 \end{pmatrix},
\begin{pmatrix} y^5 \\ 0 \\ 0 \end{pmatrix},
\begin{pmatrix} xy \\ \frac{3}{5}y^2 \\ \frac{7}{5}yz \end{pmatrix},\\&
\begin{pmatrix} xy^2 \\ 0 \\ 2y^2z \end{pmatrix},
\begin{pmatrix} 0 \\ z^2 \\ 0 \end{pmatrix},
\begin{pmatrix} 0 \\ yz \\ 0 \end{pmatrix},
\begin{pmatrix} 0 \\ y^2z \\ 0 \end{pmatrix},
\begin{pmatrix} 0 \\ y^3 \\ -y^2z \end{pmatrix},
\begin{pmatrix} 0 \\ y^3z \\ 0 \end{pmatrix},
\begin{pmatrix} 0 \\ y^4 \\ 0 \end{pmatrix},
\begin{pmatrix} 0 \\ y^5 \\ 0 \end{pmatrix},
\begin{pmatrix} 0 \\ xy \\ \frac{1}{2}y^4 \end{pmatrix},
\begin{pmatrix} 0 \\ xy^2 \\ 0 \end{pmatrix},\\&
\begin{pmatrix} 0 \\ 0 \\ z^2 \end{pmatrix},
\begin{pmatrix} 0 \\ 0 \\ y^3z \end{pmatrix},
\begin{pmatrix} 0 \\ 0 \\ y^5 \end{pmatrix}
\}.
\end{align*}
We see that the ideal $I_1\subset (z^2,yz,y^3,xy)+(f,J(f))$, $I_1^2\subset (z^2,yz,y^3,xy)^2+(f,J(f))=(x^2y^2,f,J(f))=(y^3z,f,J(f))$.

Now we define $(\beta_{ij})\in M_3(R)$ by
$$\beta_{11}=\frac{25}{168}y^5+\frac{415}{252}a_0y^3z, $$
$$\beta_{12}=\beta_{21}=-\frac{5}{28}y^2z+\frac{1}{8}a_0z^2, $$
$$\beta_{13}=\beta_{31}=0, $$
$$\beta_{22}=\frac{3}{56}xy^2-\frac{17}{56}a_0y^4+a_0^2y^2z, $$
$$\beta_{23}=\beta_{32}=-\frac{1}{16}y^5-\frac{39}{56}a_0y^3z, $$
$$\beta_{33}=0. $$
Then the lifting $(\widetilde{\beta_{ij}})\in M_3(\mathcal{O}_3)$ satisfies$$
\begin{pmatrix}
\widetilde{\beta_{11}} & \widetilde{\beta_{12}} & \widetilde{\beta_{13}}\\
\widetilde{\beta_{21}} & \widetilde{\beta_{22}} & \widetilde{\beta_{23}}\\
\widetilde{\beta_{31}} & \widetilde{\beta_{32}} & \widetilde{\beta_{33}}
\end{pmatrix}
\begin{pmatrix}
f_{x}\\
f_{y}\\
f_{z}
\end{pmatrix}
\equiv
\begin{pmatrix}
0\\
0\\
0
\end{pmatrix} \ \ \ (  \textup{mod} \ (f,J(f)^2)  ),$$ and $\widetilde{\beta_{11}}\notin (y^3z,f,J(f))$, so $\beta_{11}\notin I_1^2 \mod J(f)$.
By proposition \ref{prop2.14}, the Nakai Conjecture holds.

(2) If $a_0=0 \ \& \ a_1\ne 0$, by calculation, the Tjurina algebra $R=\mathcal{O}_3/(f,f_x,f_y,f_z)\simeq \mathbb{C}\{x,y,z\}/(x^2z,yz^2,xy^4,y^7,2xz+y^4, z^2+4xy^3+7a_1y^6,x^2+2yz)$, it has a $\mathbb{C}$-basis: $\{1,z,z^2,y,yz,y^2,\\y^2z,y^3,y^3z,y^4,y^5,y^6,x,xy,xy^2\}$, $dim_{\mathbb{C}} R=15$. 
By solving the equation $Hess(f)\cdot (\alpha_1,\alpha_2,\alpha_3)^{T}=0$ in $R$, we obtain that
\begin{align*}
\begin{pmatrix}
\alpha_1\\
\alpha_2\\
\alpha_3
\end{pmatrix}&\in Span_{\mathbb{C}}\{ 
\begin{pmatrix} z^2 \\ 0 \\ 0 \end{pmatrix},
\begin{pmatrix} yz \\ 0 \\ \frac{1}{2}y^4 \end{pmatrix},
\begin{pmatrix} y^2z \\ 0 \\ \frac{1}{2}y^5 \end{pmatrix},
\begin{pmatrix} y^3 \\ \frac{1}{2}z \\ xy^2+\frac{7}{2}a_1y^5 \end{pmatrix},
\begin{pmatrix} y^3z \\ 0 \\ 0 \end{pmatrix},
\begin{pmatrix} y^4 \\ 0 \\ 0 \end{pmatrix},
\begin{pmatrix} y^5 \\ 0 \\ 0 \end{pmatrix},
\begin{pmatrix} y^6 \\ 0 \\ 0 \end{pmatrix},\\&
\begin{pmatrix} xy \\ \frac{3}{5}y^2 \\ \frac{7}{5}yz \end{pmatrix},
\begin{pmatrix} xy^2 \\ \frac{3}{5}y^3 \\ \frac{7}{5}y^2z \end{pmatrix},
\begin{pmatrix} 0 \\ z^2 \\ 0 \end{pmatrix},
\begin{pmatrix} 0 \\ yz \\ 0 \end{pmatrix},
\begin{pmatrix} 0 \\ y^2z \\ 0 \end{pmatrix},
\begin{pmatrix} 0 \\ y^3z \\ 0 \end{pmatrix},
\begin{pmatrix} 0 \\ y^4 \\ 0 \end{pmatrix},
\begin{pmatrix} 0 \\ y^5 \\ 0 \end{pmatrix},
\begin{pmatrix} 0 \\ y^6 \\ 0 \end{pmatrix},
\begin{pmatrix} 0 \\ xy \\ \frac{1}{2}y^4 \end{pmatrix},\\&
\begin{pmatrix} 0 \\ xy^2 \\ \frac{1}{2}y^5 \end{pmatrix},
\begin{pmatrix} 0 \\ 0 \\ z^2 \end{pmatrix},
\begin{pmatrix} 0 \\ 0 \\ y^3z \end{pmatrix},
\begin{pmatrix} 0 \\ 0 \\ y^6 \end{pmatrix}
\}.
\end{align*}
We see that the ideal $I_1\subset (z^2,yz,y^3,xy)+(f,J(f))$, $I_1^2\subset (z^2,yz,y^3,xy)^2+(f,J(f))=(y^6,x^2y^2,f,J(f))=(y^6, y^3z,f,J(f))$.

Now we define $(\beta_{ij})\in M_3(R)$ by
$$\beta_{11}=y^2z, $$
$$\beta_{12}=\beta_{21}=-\frac{3}{10}xy^2, $$
$$\beta_{13}=\beta_{31}=\frac{7}{20}y^5, $$
$$\beta_{22}=-\frac{9}{50}y^3+\frac{24}{125}a_1y^2z, $$
$$\beta_{23}=\beta_{32}=-\frac{21}{50}y^2z, $$
$$\beta_{33}=0. $$
Then the lifting $(\widetilde{\beta_{ij}})\in M_3(\mathcal{O}_3)$ satisfies$$
\begin{pmatrix}
\widetilde{\beta_{11}} & \widetilde{\beta_{12}} & \widetilde{\beta_{13}}\\
\widetilde{\beta_{21}} & \widetilde{\beta_{22}} & \widetilde{\beta_{23}}\\
\widetilde{\beta_{31}} & \widetilde{\beta_{32}} & \widetilde{\beta_{33}}
\end{pmatrix}
\begin{pmatrix}
f_{x}\\
f_{y}\\
f_{z}
\end{pmatrix}
\equiv
\begin{pmatrix}
0\\
0\\
0
\end{pmatrix} \ \ \ (  \textup{mod} \ (f,J(f)^2)  ),$$ and $\widetilde{\beta_{11}}\notin (y^6, y^3z,f,J(f))$, so $\beta_{11}\notin I_1^2 \mod J(f)$.
By proposition \ref{prop2.14}, the Nakai Conjecture holds.
\\

13. $S_{17}$ singularity

For a $S_{17}$ singularity $(V(f),0)$ defined by $f=x^2z+yz^2+y^6+(a_0+a_1y)y^4z$ for some $a_0,a_1$, we treat the case of $a_0=0  \ \&  \ a_1\ne 0$ and the case of $a_0\ne 0$ separately.

(1) If $a_0\ne 0$, by calculation, the Tjurina algebra $R=\mathcal{O}_3/(f,f_x,f_y,f_z)\simeq \mathbb{C}\{x,y,z\}/(xz,yz^2,\\y^6,y^4z, z^2+6y^5+4a_0y^3z,x^2+2yz+(a_0+a_1y)y^4)$, it has a $\mathbb{C}$-basis: $\{1,z,z^2,y,yz,y^2,y^2z,y^3,\\y^3z,y^4,x,xy,xy^2,\cdots,xy^4\}$, $dim_{\mathbb{C}} R=15$. 
By solving the equation $Hess(f)\cdot (\alpha_1,\alpha_2,\alpha_3)^{T}=0$ in $R$, we obtain that
\begin{align*}
\begin{pmatrix}
\alpha_1\\
\alpha_2\\
\alpha_3
\end{pmatrix}&\in Span_{\mathbb{C}}\{ 
\begin{pmatrix} z^2 \\ 0 \\ 0 \end{pmatrix},
\begin{pmatrix} yz \\ 0 \\ 0 \end{pmatrix},
\begin{pmatrix} y^2z \\ 0 \\ 0 \end{pmatrix},
\begin{pmatrix} y^3z \\ 0 \\ 0 \end{pmatrix},
\begin{pmatrix} y^4 \\ 0 \\ -xy^3 \end{pmatrix},
\begin{pmatrix} xy \\ \frac{1}{2}y^2 \\ \frac{3}{2}yz \end{pmatrix},
\begin{pmatrix} xy^2 \\ 0 \\ 2y^2z \end{pmatrix},
\begin{pmatrix} xy^3 \\ 0 \\ 0 \end{pmatrix},
\begin{pmatrix} xy^4 \\ 0 \\ 0 \end{pmatrix},\\&
\begin{pmatrix} 0 \\ 3y^2+a_0z \\ 2a_0^2y^2z-3yz \end{pmatrix},
\begin{pmatrix} 0 \\ z^2 \\ 0 \end{pmatrix},
\begin{pmatrix} 0 \\ yz \\ 0 \end{pmatrix},
\begin{pmatrix} 0 \\ y^2z \\ 0 \end{pmatrix},
\begin{pmatrix} 0 \\ y^3 \\ -y^2z \end{pmatrix},
\begin{pmatrix} 0 \\ y^3z \\ 0 \end{pmatrix},
\begin{pmatrix} 0 \\ y^4 \\ 0 \end{pmatrix},
\begin{pmatrix} 0 \\ xy \\ -2a_0xy^3 \end{pmatrix},
\begin{pmatrix} 0 \\ xy^2 \\ 0 \end{pmatrix},\\&
\begin{pmatrix} 0 \\ xy^3 \\ 0 \end{pmatrix},
\begin{pmatrix} 0 \\ xy^4 \\ 0 \end{pmatrix},
\begin{pmatrix} 0 \\ 0 \\ z^2 \end{pmatrix},
\begin{pmatrix} 0 \\ 0 \\ y^3z \end{pmatrix},
\begin{pmatrix} 0 \\ 0 \\ xy^4 \end{pmatrix}
\}.
\end{align*}
We see that the ideal $I_1\subset (z^2,yz,y^4,xy)+(f,J(f))$, $I_1^2\subset (z^2,yz,y^4,xy)^2+(f,J(f))=(x^2y^2,f,J(f))=(y^3z,f,J(f))$.

Now we define $(\beta_{ij})\in M_3(R)$ by
$$\beta_{11}=a_0^2 y^3 z - \frac{735}{527}y^2 z + \frac{70}{527}a_0z^2, $$
$$\beta_{12}=\beta_{21}= -\frac{494}{3689}a_0^2x y^3 + \frac{210}{527}xy^2, $$
$$\beta_{13}=\beta_{31}=0, $$
$$\beta_{22}=-\frac{247}{3689}a_0^2 y^4 +
\frac{120}{527}y^3 + \frac{1902}{3689}a_1y^2z - \frac{74}{3689}a_0yz, $$ 
$$\beta_{23}=\beta_{32}= -\frac{494}{3689}a_0^2 y^3z
+ \frac{300}{527}y^2z + \frac{25}{3689}a_0z^2, $$
$$\beta_{33}= 0. $$
Then the lifting $(\widetilde{\beta_{ij}})\in M_3(\mathcal{O}_3)$ satisfies$$
\begin{pmatrix}
\widetilde{\beta_{11}} & \widetilde{\beta_{12}} & \widetilde{\beta_{13}}\\
\widetilde{\beta_{21}} & \widetilde{\beta_{22}} & \widetilde{\beta_{23}}\\
\widetilde{\beta_{31}} & \widetilde{\beta_{32}} & \widetilde{\beta_{33}}
\end{pmatrix}
\begin{pmatrix}
f_{x}\\
f_{y}\\
f_{z}
\end{pmatrix}
\equiv
\begin{pmatrix}
0\\
0\\
0
\end{pmatrix} \ \ \ (  \textup{mod} \ (f,J(f)^2)  ),$$ and $\widetilde{\beta_{11}}\notin (y^3z,f,J(f))$, so $\beta_{11}\notin I_1^2 \mod J(f)$.
By proposition \ref{prop2.14}, the Nakai Conjecture holds.

(2) If $a_0=0 \ \& \ a_1\ne 0$, by calculation, the Tjurina algebra $R=\mathcal{O}_3/(f,f_x,f_y,f_z)\simeq \mathbb{C}\{x,y,z\}/(xz,yz^2,y^6,y^5z, z^2+6y^5+5a_1y^4z,x^2+2yz+a_1y^5)$, it has a $\mathbb{C}$-basis: $\{1,z,z^2,y,yz,\\y^2,y^2z,y^3,y^3z,y^4,y^4z,x,xy,xy^2,\cdots,xy^4\}$, $dim_{\mathbb{C}} R=16$. 
By solving the equation $Hess(f)\cdot (\alpha_1,\alpha_2,\alpha_3)^{T}=0$ in $R$, we obtain that
\begin{align*}
\begin{pmatrix}
\alpha_1\\
\alpha_2\\
\alpha_3
\end{pmatrix}&\in Span_{\mathbb{C}}\{ 
\begin{pmatrix} z^2 \\ 0 \\ 0 \end{pmatrix},
\begin{pmatrix} yz \\ 0 \\ 0 \end{pmatrix},
\begin{pmatrix} y^2z \\ 0 \\ 0 \end{pmatrix},
\begin{pmatrix} y^3z \\ 0 \\ 0 \end{pmatrix},
\begin{pmatrix} y^4z \\ 0 \\ 0 \end{pmatrix},
\begin{pmatrix} xy \\ 0 \\ 2yz \end{pmatrix},
\begin{pmatrix} xy^2 \\ 0 \\ 2y^2z \end{pmatrix},
\begin{pmatrix} xy^3 \\ 0 \\ 2y^3z \end{pmatrix},
\begin{pmatrix} xy^4 \\ 0 \\ 0 \end{pmatrix},\\&
\begin{pmatrix} 0 \\ z^2 \\ 0 \end{pmatrix},
\begin{pmatrix} 0 \\ yz \\ 0 \end{pmatrix},
\begin{pmatrix} 0 \\ y^2 \\ -yz \end{pmatrix},
\begin{pmatrix} 0 \\ y^2z \\ 0 \end{pmatrix},
\begin{pmatrix} 0 \\ y^3 \\ -y^2z \end{pmatrix},
\begin{pmatrix} 0 \\ y^3z \\ 0 \end{pmatrix},
\begin{pmatrix} 0 \\ y^4 \\ -y^3z \end{pmatrix},
\begin{pmatrix} 0 \\ y^4z \\ 0 \end{pmatrix},
\begin{pmatrix} 0 \\ xy \\ 0 \end{pmatrix},
\begin{pmatrix} 0 \\ xy^2 \\ 0 \end{pmatrix},\\&
\begin{pmatrix} 0 \\ xy^3 \\ 0 \end{pmatrix},
\begin{pmatrix} 0 \\ xy^4 \\ 0 \end{pmatrix},
\begin{pmatrix} 0 \\ 0 \\ z^2 \end{pmatrix},
\begin{pmatrix} 0 \\ 0 \\ y^4z \end{pmatrix},
\begin{pmatrix} 0 \\ 0 \\ xy^4 \end{pmatrix}
\}.
\end{align*}
We see that the ideal $I_1\subset (z^2,yz,xy)+(f,J(f))$, $I_1^2\subset (z^2,yz,xy)^2+(f,J(f))=(x^2y^2,f,J(f))=(y^3z,f,J(f))$.

Now we define $(\beta_{ij})\in M_3(R)$ by
$$\beta_{11}= -\frac{42}{5} yz + a_1z^2, $$
$$\beta_{12}=\beta_{21}=\frac{10}{7}a_1^2 xy^4 + \frac{12}{5}xy, $$
$$\beta_{13}=\beta_{31}=0, $$
$$\beta_{22}=\frac{48}{35}y^2 - \frac{444}{1225}a_1yz -
\frac{4463}{25725}a_1^2z^2, $$ 
$$\beta_{23}=\beta_{32}= \frac{29056}{8575}a_1^2y^4z + \frac{24}{7}yz + \frac{6}{49}a_1z^2, $$
$$\beta_{33}= -\frac{216}{49}a_1y^4z + \frac{60}{7}z^2. $$
Then the lifting $(\widetilde{\beta_{ij}})\in M_3(\mathcal{O}_3)$ satisfies$$
\begin{pmatrix}
\widetilde{\beta_{11}} & \widetilde{\beta_{12}} & \widetilde{\beta_{13}}\\
\widetilde{\beta_{21}} & \widetilde{\beta_{22}} & \widetilde{\beta_{23}}\\
\widetilde{\beta_{31}} & \widetilde{\beta_{32}} & \widetilde{\beta_{33}}
\end{pmatrix}
\begin{pmatrix}
f_{x}\\
f_{y}\\
f_{z}
\end{pmatrix}
\equiv
\begin{pmatrix}
0\\
0\\
0
\end{pmatrix} \ \ \ (  \textup{mod} \ (f,J(f)^2)  ),$$ and $\widetilde{\beta_{11}}\notin (y^3z,f,J(f))$, so $\beta_{11}\notin I_1^2 \mod J(f)$.
By proposition \ref{prop2.14}, the Nakai Conjecture holds.
\\

14. $U_{16}$ singularity

For a $U_{16}$ singularity $(V(f),0)$ defined by $f=x^3+xz^2+y^5+(a_0+a_1y)x^2y^2$ for some $a_0,a_1$, we treat the case of $a_0=0  \ \&  \ a_1\ne 0$ and the case of $a_0\ne 0$ separately.

(1) If $a_0\ne 0$, by calculation, the Tjurina algebra $R=\mathcal{O}_3/(f,f_x,f_y,f_z)\simeq \mathbb{C}\{x,y,z\}/(x^3,xz,y^5,\\x^2y^2,3x^2+z^2+2(a_0+a_1y)xy^2,5y^4+2a_0x^2y)$, it has a $\mathbb{C}$-basis: $\{1,z,z^2,y,yz,yz^2,y^2,y^2z,y^3,y^3z,y^4,\\x,xy,xy^2\}$, $dim_{\mathbb{C}} R=14$. 
By solving the equation $Hess(f)\cdot (\alpha_1,\alpha_2,\alpha_3)^{T}=0$ in $R$, we obtain that
\begin{align*}
\begin{pmatrix}
\alpha_1\\
\alpha_2\\
\alpha_3
\end{pmatrix}&\in Span_{\mathbb{C}}\{ 
\begin{pmatrix} z^2 \\ 0 \\ 0 \end{pmatrix},
\begin{pmatrix} yz^2 \\ 0 \\ 0 \end{pmatrix},
\begin{pmatrix} y^2z \\ 0 \\ 0 \end{pmatrix},
\begin{pmatrix} y^3z \\ 0 \\ 0 \end{pmatrix},
\begin{pmatrix} y^4 \\ 0 \\ 0 \end{pmatrix},
\begin{pmatrix} xy \\ \frac{1}{2}y^2 \\ yz \end{pmatrix},
\begin{pmatrix} xy^2 \\ 0 \\ 0 \end{pmatrix},
\begin{pmatrix} 0 \\ z^2 \\ 0 \end{pmatrix},
\begin{pmatrix} 0 \\ yz \\ 0 \end{pmatrix},\\&
\begin{pmatrix} 0 \\ yz^2 \\ 0 \end{pmatrix},
\begin{pmatrix} 0 \\ y^2z \\ 0 \end{pmatrix},
\begin{pmatrix} 0 \\ y^3 \\ 0 \end{pmatrix},
\begin{pmatrix} 0 \\ y^3z \\ 0 \end{pmatrix},
\begin{pmatrix} 0 \\ y^4 \\ 0 \end{pmatrix},
\begin{pmatrix} 0 \\ xy \\ 0 \end{pmatrix},
\begin{pmatrix} 0 \\ xy^2 \\ 0 \end{pmatrix},
\begin{pmatrix} 0 \\ 0 \\ z^2 \end{pmatrix},
\begin{pmatrix} 0 \\ 0 \\ yz^2 \end{pmatrix},
\begin{pmatrix} 0 \\ 0 \\ y^2z \end{pmatrix},\\&
\begin{pmatrix} 0 \\ 0 \\ y^3z \end{pmatrix},
\begin{pmatrix} 0 \\ 0 \\ y^4 \end{pmatrix},
\begin{pmatrix} 0 \\ 0 \\ xy^2 \end{pmatrix}
\}.
\end{align*}
We see that the ideal $I_2\subset (y^2,z^2,yz,xy)+(f,J(f))$, $I_2^2\subset (y^2,z^2,yz,xy)^2+(f,J(f))=(y^4,y^3z,xy^3,f,J(f))$.

Now we define $(\beta_{ij})\in M_3(R)$ by
$$\beta_{11}=\frac{575}{16}y^4 - \frac{5}{8}a_0yz^2, $$
$$\beta_{12}=\beta_{21}= -\frac{15}{2}a_0xy^2 - \frac{9}{2}a_0^2y^4, $$
$$\beta_{13}=\beta_{31}=0, $$
$$\beta_{22}= a_0^3xy^2 - \frac{9}{2}a_0y^3 -\frac{1}{10}a_0^2z^2, $$ 
$$\beta_{23}=\beta_{32}= -\frac{15}{2}a_0y^2z, $$
$$\beta_{33}= \frac{75}{16}y^4 - \frac{105}{8}a_0yz^2. $$
Then the lifting $(\widetilde{\beta_{ij}})\in M_3(\mathcal{O}_3)$ satisfies$$
\begin{pmatrix}
\widetilde{\beta_{11}} & \widetilde{\beta_{12}} & \widetilde{\beta_{13}}\\
\widetilde{\beta_{21}} & \widetilde{\beta_{22}} & \widetilde{\beta_{23}}\\
\widetilde{\beta_{31}} & \widetilde{\beta_{32}} & \widetilde{\beta_{33}}
\end{pmatrix}
\begin{pmatrix}
f_{x}\\
f_{y}\\
f_{z}
\end{pmatrix}
\equiv
\begin{pmatrix}
0\\
0\\
0
\end{pmatrix} \ \ \ (  \textup{mod} \ (f,J(f)^2)  ),$$ and $\widetilde{\beta_{22}}\notin (y^4,y^3z,xy^3,f,J(f))$, so $\beta_{22}\notin I_2^2 \in J(f)$.
By proposition \ref{prop2.14}, the Nakai Conjecture holds.

(2) If $a_0=0 \ \& \ a_1\ne 0$, by calculation, the Tjurina algebra $R=\mathcal{O}_3/(f,f_x,f_y,f_z)\simeq \mathbb{C}\{x,y,z\}/(x^3,xz,y^5,x^2y^3,3x^2+z^2+2a_1xy^3,5y^4+3a_1x^2y^2)$, it has a $\mathbb{C}$-basis: $\{1,z,z^2,y,yz,yz^2,\\y^2,y^2z,y^2z^2,y^3,y^3z,x,xy,xy^2,xy^3\}$, $dim_{\mathbb{C}} R=15$. 
By solving the equation $Hess(f)\cdot (\alpha_1,\alpha_2,\alpha_3)^{T}\\=0$ in $R$, we obtain that
\begin{align*}
\begin{pmatrix}
\alpha_1\\
\alpha_2\\
\alpha_3
\end{pmatrix}&\in Span_{\mathbb{C}}\{ 
\begin{pmatrix} z^2 \\ 0 \\ 0 \end{pmatrix},
\begin{pmatrix} yz \\ 0 \\ 3xy \end{pmatrix},
\begin{pmatrix} yz^2 \\ 0 \\ 0 \end{pmatrix},
\begin{pmatrix} y^2z \\ 0 \\ 3xy^2 \end{pmatrix},
\begin{pmatrix} y^2z^2 \\ 0 \\ 0 \end{pmatrix},
\begin{pmatrix} y^3z \\ 0 \\ 0 \end{pmatrix},
\begin{pmatrix} xy \\ 0 \\ yz \end{pmatrix},
\begin{pmatrix} xy^2 \\ 0 \\ y^2z \end{pmatrix},
\begin{pmatrix} xy^3 \\ 0 \\ 0 \end{pmatrix},\\&
\begin{pmatrix} 0 \\ z^2 \\ 0 \end{pmatrix},
\begin{pmatrix} 0 \\ yz \\ 0 \end{pmatrix},
\begin{pmatrix} 0 \\ yz^2 \\ 0 \end{pmatrix},
\begin{pmatrix} 0 \\ y^2 \\ 0 \end{pmatrix},
\begin{pmatrix} 0 \\ y^2z \\ 0 \end{pmatrix},
\begin{pmatrix} 0 \\ y^2z^2 \\ 0 \end{pmatrix},
\begin{pmatrix} 0 \\ y^3 \\ 0 \end{pmatrix},
\begin{pmatrix} 0 \\ y^3z \\ 0 \end{pmatrix},
\begin{pmatrix} 0 \\ xy \\ 0 \end{pmatrix},
\begin{pmatrix} 0 \\ xy^2 \\ 0 \end{pmatrix},\\&
\begin{pmatrix} 0 \\ xy^3 \\ 0 \end{pmatrix},
\begin{pmatrix} 0 \\ 0 \\ z^2 \end{pmatrix},
\begin{pmatrix} 0 \\ 0 \\ yz^2 \end{pmatrix},
\begin{pmatrix} 0 \\ 0 \\ y^2z^2 \end{pmatrix},
\begin{pmatrix} 0 \\ 0 \\ y^3z \end{pmatrix},
\begin{pmatrix} 0 \\ 0 \\ xy^3 \end{pmatrix}
\}.
\end{align*}
We see that the ideal $I_2\subset (y^2,z^2,yz,xy)+(f,J(f))$, $I_2^2\subset (y^2,z^2,yz,xy)^2+(f,J(f))=(y^4,y^2z^2,y^3z,xy^3,f,J(f))$.

Now we define $(\beta_{ij})\in M_3(R)$ by
$$\beta_{11}= \frac{16}{5}a_1xy^3 + z^2, $$
$$\beta_{12}=\beta_{21}=-\frac{9}{5}xy, $$
$$\beta_{13}=\beta_{31}= 0, $$
$$\beta_{22}=-\frac{12}{625}a_1^2 xy^3 - \frac{27}{25}y^2 - \frac{12}{125}a_1z^2, $$ 
$$\beta_{23}=\beta_{32}= -\frac{9}{5}yz, $$ 
$$\beta_{33}= \frac{6}{5}a_1xy^3 - 3 z^2. $$ 
Then the lifting $(\widetilde{\beta_{ij}})\in M_3(\mathcal{O}_3)$ satisfies$$
\begin{pmatrix}
\widetilde{\beta_{11}} & \widetilde{\beta_{12}} & \widetilde{\beta_{13}}\\
\widetilde{\beta_{21}} & \widetilde{\beta_{22}} & \widetilde{\beta_{23}}\\
\widetilde{\beta_{31}} & \widetilde{\beta_{32}} & \widetilde{\beta_{33}}
\end{pmatrix}
\begin{pmatrix}
f_{x}\\
f_{y}\\
f_{z}
\end{pmatrix}
\equiv
\begin{pmatrix}
0\\
0\\
0
\end{pmatrix} \ \ \ (  \textup{mod} \ (f,J(f)^2)  ),$$ and $\widetilde{\beta_{22}}\notin (y^4,y^2z^2,y^3z,xy^3,f,J(f))$, so $\beta_{22}\notin I_2^2 \mod J(f)$.
By proposition \ref{prop2.14}, the Nakai Conjecture holds.

From the calculations in section \ref{sec 3} and section \ref{sec 4}, we have finished the proof of the main theorem.

\end{document}